\documentclass[openany,12pt,reqno]{amsart}
\usepackage{amsthm,amsmath,amssymb,enumerate,mathtools}
\newtheorem{theorem}{Theorem}
\newtheorem{claim}{}[theorem]
\newtheorem{lemma}[theorem]{Lemma}
\newtheorem{problem}[theorem]{Problem}
\newtheorem{proposition}[theorem]{Proposition}
\newtheorem{corollary}[theorem]{Corollary}
\newtheorem{conjecture}[theorem]{Conjecture}

\newtheorem*{innerthm}{\innerthmname}
\newcommand{\innerthmname}{}
\newenvironment{statement}[1]
 {\renewcommand{\innerthmname}{#1}\innerthm}
 {\endinnerthm}
 
 \theoremstyle{definition}

\newcommand{\bF}{\mathbb F}
\newcommand{\bR}{\mathbb R}
\newcommand{\bC}{\mathbb C}
\newcommand{\bZ}{\mathbb Z}

\newcommand{\cA}{\mathcal{A}}
\newcommand{\cB}{\mathcal{B}}
\newcommand{\cC}{\mathcal{C}}
\newcommand{\cR}{\mathcal{R}}
\newcommand{\cF}{\mathcal{F}}

\newcommand{\cK}{\mathcal{K}}

\newcommand{\cL}{\mathcal{L}}

\newcommand{\cH}{\mathcal{H}}
\newcommand{\cM}{\mathcal{M}}
\newcommand{\cP}{\mathcal{P}}

\newcommand{\cS}{\mathcal{S}}
\newcommand{\cT}{\mathcal{T}}
\newcommand{\cU}{\mathcal{U}}

\newcommand{\cX}{\mathcal{X}}

\DeclareMathOperator{\si}{si}

\DeclareMathOperator{\cl}{cl}

\DeclareMathOperator{\dist}{dist}

\DeclareMathOperator{\PG}{PG}
\DeclareMathOperator{\DG}{DG}
\DeclareMathOperator{\LG}{LG}

\DeclareMathOperator{\GF}{GF}

\DeclareMathOperator{\rank}{rank}
\newcommand{\elem}{\varepsilon}
\newcommand{\del}{\backslash}
\newcommand{\con}{/}

\usepackage{comment,chngcntr,setspace}

\pdfoutput=1

\newenvironment{alignLetter}{
    \setcounter{equation}{0}
    
    \align
}{
    \endalign
}

\counterwithout{equation}{section}

\newif\ifhideproofs
\ifhideproofs
\usepackage{environ}
\NewEnviron{hide}{}

\fi

\numberwithin{theorem}{subsection}

\author[Geelen]{Jim Geelen}
\author[Nelson]{Peter Nelson}
\author[Walsh]{Zach Walsh}
\title{Excluding a Line from $\bC$-Representable Matroids}

\date{February 2025. }

\newcommand{\ep}{\epsilon}
\newcommand{\cO}{\mathcal O}
\newcommand{\pT}{\mathcal T}
\newcommand{\cY}{\mathcal Y}
\newcommand{\A}{{\bf A}}
\newcommand{\Y}{{\bf Y}}
\newcommand{\Q}{{\bf Q}}
\newcommand{\R}{{\bf R}}
\newcommand{\bQ}{\mathbb{Q}}
\newcommand{\bK}{\mathbb K}

\begin{document}

\begin{abstract} 
For each positive integer $t$ and each sufficiently large integer $r$, we show that the maximum number of elements of a simple, rank-$r$, $\bC$-representable matroid with no $U_{2,t+3}$-minor is $t{r\choose 2}+r$.
We derive this as a consequence of a much more general result concerning matroids on group-labeled graphs.
\end{abstract}

\maketitle


\begin{section}{Introduction} \label{intro}
\pagenumbering{arabic}
	
A \emph{matroid} is a pair $M = (E,r)$, where $E$ is a finite `ground' set and $r : 2^E \to \mathbb{Z}$ is a `rank' function for which $0 \le r(X) \le |X|$ for each $X$, and $r$ is both nondecreasing with respect to set inclusion, and submodular: that is $r(X') \le r(X)$ for each $X' \subseteq X \subseteq E$, and $r(X) + r(Y) \ge r(X \cup Y) + r(X \cap Y)$ for all $X,Y \subseteq E$. These objects were introduced independently by Whitney \cite{Whitney} and Nakasawa \cite{Nakasawa1,Nakasawa2,Nakasawa3} as a combinatorial abstraction of rank functions as they arise in linear algebra. Indeed, the most natural examples of matroids come from matrices; if $A$ is a matrix over a field $\mathbb{F}$ with columns indexed by $E$, then the function $r : 2^E \to \bZ$ defined by $r(X) = \rank(A[X])$ gives a matroid on $E$. Such a matroid is \emph{$\bF$-representable}. 

The \emph{rank} $r(M)$ of a matroid is the rank of its ground set, $r(E)$. The \emph{size} $|M|$ of $M$ is the cardinality of its ground set $E$. This article concerns the natural problem of bounding size above by a function of rank. In his seminal 1993 survey paper \cite{Kung1}, Kung phrased this as follows:

\begin{problem} \label{Kung}
	Let $\cM$ be a class of matroids. For each positive integer $r$, determine the maximum size of a rank-$r$ matroid in $\cM$.
\end{problem}

When $\cM$ is the class of all matroids, this problem is not meaningful, since a rank-$r$ matroid can have arbitrarily many elements. For example, for each integer $a \ge 0$, the \emph{rank-$a$ uniform matroid} $U_{a,b}$, defined by $|E| = b$ and $r(X) = \min(|X|,a)$ for each $X \subseteq E$, exists for all $b \ge a$. 

For many classes, however, the size of the matroids in the class is bounded by a function of their rank. A matroid $M = (E,r)$ is \emph{simple} if $r(X) = |X|$ whenever $|X| \le 2$. If $M$ is represented by a matrix $A$, then simplicity is equivalent to the assertion that the columns of $A$ are nonzero and pairwise non-parallel. For this reason, it is easy to argue that, if $\GF(q)$ denotes the $q$-element finite field, then a simple rank-$r$ $\GF(q)$-representable matroid has size at most $\frac{q^r-1}{q-1}$, since a rank-$r$ matrix over such a field with pairwise non-parallel, nonzero columns has at most $\frac{q^r-1}{q-1}$ columns. A $\GF(q)$-representable matroid attaining this bound is a rank-$r$ \emph{projective geometry}, denoted $\PG(r-1,q)$.

In fact, a similar upper bound applies under a much weaker hypothesis than representability, where we simply forbid a certain substructure from $M$. Given an element $e$ of a matroid $M = (E,r)$, we can \emph{delete} $e$ from $M$ to obtain the matroid $M \del e = (E-\{e\},r')$ where $r'(X) = r(X)$ for each $X \subseteq E - \{e\}$, or \emph{contract} $e$ from $M$ to obtain the matroid $M \con e = (E-\{e\}, r'')$, where $r''(X) = r(X \cup \{e\}) - r(\{e\})$ for each $X \subseteq E-\{e\}$. A matroid obtained from $M$ by a sequence of deletions or contractions is a \emph{minor} of $M$. 
Kung \cite{Kung1} showed that size is bounded above by a function of rank in any simple matroid omitting a fixed rank-$2$ uniform minor.

\begin{theorem}[Kung] \label{l-Kung}
	For each integer $\ell\ge 2$, each simple rank-$r$ matroid $M$ with no $U_{2,\ell+2}$-minor satisfies $|M|\le \frac{\ell^{r}-1}{\ell-1}$. 
\end{theorem}

This upper bound is best possible in the weak sense that any upper bound must be exponential in $r$; if $q$ is a prime power with $q \le \ell$, then the rank-$r$ projective geometry $\PG(r-1,q)$ has no $U_{2,\ell+2}$-minor and has size $\frac{q^r-1}{q-1}$. 
In fact, for large $r$ these projective geometries are the unique tight examples \cite{n-line}:

\begin{theorem}[Geelen, Nelson] \label{line}
	Let $\ell\ge 2$ be an integer, and let $q$ be the largest prime power not exceeding $\ell$.
	Then, for all sufficiently large $r$, each simple rank-$r$ matroid $M$ with no $U_{2,\ell+2}$-minor satisfies $|M|\le \frac{q^r-1}{q-1}$. Moreover, equality holds only if $M$ is isomorphic to $\PG(r-1,q)$. 
\end{theorem}

This result is surprising in that the purely combinatorial hypothesis of excluding a rank-$2$ uniform minor gives an upper bound attained uniquely by nontrivial algebraic objects: projective geometries over finite fields. 

There is a natural additional hypothesis that drastically changes the flavour of the problem: imposing representability over a field of characteristic zero. One can show that nontrivial finite projective geometries are not representable over $\bC$ (or, therefore, $\bR$ or $\bQ$), and thus the bound in Theorem~\ref{line} can be strengthened under the hypothesis that $M$ is $\bC$-representable. In fact, it can be improved rather dramatically. In  \cite{Kung2}, Kung showed the following:

\begin{theorem}\label{kungc}
	 Let $\ell \ge 2$ be an integer. If $M$ is a simple rank-$r$ $\bC$-representable matroid with no $U_{2,\ell+2}$-minor, then 
	 \[|M| \le \big(\ell^{2^{\ell-1}-1}-\ell^{2^{\ell-1}-2}\big)\binom{r+1}{2} - r.\]
\end{theorem}

That is, the upper bound drops from exponential to quadratic in $r$. In turn, no bound better than quadratic in $r$ is possible; there is a natural example showing this. Let $t \ge 2$ be an integer, and let $\zeta \in \bC$ be a $(t-1)$-st primitive root of unity. Let $D=D_{r,t-1}$ denote the $\bC$-matrix having as columns the standard basis vectors $b_1, \dotsc, b_r$, as well as all vectors of the form $b_i - \zeta^k b_j$, where $1 \le i < j \le r$ and $0 \le k < t-1$. The matroid represented by $D$ is called a \emph{rank-$r$ cyclic Dowling geometry of order $t-1$}, and can be shown to have no $U_{2,t+2}$-minor. The matrix $D$ has $r + (t-1)\binom{r}{2}$ columns, so the best improvement on the upper bound in Theorem~\ref{kungc} that one could hope for is $(\ell-1)\binom{r}{2} + r$, in which the doubly-exponential term in $\ell$ is replaced with a linear one. Our main result shows that this bound is correct for all large $r$. 

\begin{theorem}\label{mainc}
	Let $\ell \ge 2$ be an integer. If $r$ is sufficiently large, and $M$ is a simple rank-$r$ $\bC$-representable matroid with no $U_{2,\ell+2}$-minor, then 
	\[|M| \le (\ell-1)\binom{r}{2} + r.\]
	Moreover, if equality holds, then $M$ is isomorphic to a rank-$r$ cyclic Dowling geometry of order $\ell-1$. 	
\end{theorem}

If $M$ is instead required to be $\bR$-representable, most cyclic Dowling geometries also cease to be examples; the only exceptions are those of order $1$ or $2$, where all entries of $D$ are real. A rank-$r$ cyclic Dowling geometry of order $2$ has no $U_{2,\ell+2}$-minor for any $\ell \ge 3$, and has size $2\binom{r}{2} + r = r^2$; our second main result is that, up to a linear error term, this is the correct upper bound for $\bR$-representable matroids. 

\begin{theorem}\label{mainr}
	For each integer $\ell \ge 3$, there is a constant $c$ so that if $M$ is a simple rank-$r$ $\bR$-representable matroid with no $U_{2,\ell+2}$-minor, then $|M| \le r^2 + c\cdot r$. 
\end{theorem}

We prove both of the above theorems using a much stronger result, Theorem~\ref{main corollary}, that applies to all matroids, not just representable ones. Our theorem yields asymptotically best-possible upper bounds in a variety of settings where the extremal matroids are close to Dowling geometries, and further gives exact upper bounds in cases where Dowling geometries are the tight examples. One such family of examples arises by imposing representability over fields of different characteristic. Let $\bF^\times$ denote the multiplicative group of a field $\bF$. 

\begin{theorem}\label{mainf1f2}
	Let $\bF_1$ and $\bF_2$ be fields of different characteristic, for which $\bF_1$ is finite and $\bF_1^\times$ is a subgroup of $\bF_2^{\times}$. If $r$ is sufficiently large, and $M$ is a simple rank-$r$ matroid representable over $\bF_1$ and $\bF_2$, then $|M| \le (|\bF_1|-1)\binom{r}{2}+r$.
	Moreover, if equality holds, then $M$ is isomorphic to a rank-$r$ cyclic Dowling geometry of order $|\bF_1|-1$. 
\end{theorem}

These hypotheses apply in particular if $\bF_2 = \bC$, or if $\bF_2$ is finite and $|\bF_1|-1$ divides $|\bF_2|-1$. This problem was previously open for all $|\bF_1| > 3$; the best general result was due to Kung \cite{Kung2}, who proved that $|M|$ satisfies the bound given in Theorem \ref{kungc}, with $\ell=|\bF_1|$.

\subsection*{Extremal Functions}

In \cite{Kung1}, Kung introduced a piece of notation to facilitate discussions of Problem~\ref{Kung}: the \emph{extremal function}. For a nonempty class $\cM$ of matroids, the \emph{extremal function} of $\cM$ is the function $h_{\cM}(r)\colon \bZ^{\ge 0}\to (\bZ^{\ge 0}\cup\{\infty\})$ whose value at an integer $r\ge 0$ is the maximum number of elements of a simple matroid in $\cM$ of rank at most $r$. Much later, it was proved in \cite{GRT} that, if $\cM$ is closed under both minors and isomorphism (for convenience we use the term \emph{minor-closed} to describe this property) and does not contain all rank-$2$ uniform matroids, then the extremal function of $\cM$ has one of three behaviours: linear, quadratic, and exponential.

\begin{theorem}[Growth Rate Theorem] \label{GRT}
	If $\cM$ is a minor-closed class of matroids not containing all rank-$2$ uniform matroids, then there is a constant $c = c_{\cM}$ so that either
	\begin{enumerate}
		\item\label{lgrt} $h_{\cM}(r)\le c\cdot r$ for all $r\ge 0$, or
		
		\item\label{qgrt} ${r+1\choose 2}\le h_{\cM}(r)\le c \cdot r^2$ for all $r\ge 0$ and $\cM$ contains all cyclic Dowling geometries of order $1$, or
		
		\item\label{egrt} there is a prime power $q$ such that $\frac{q^r-1}{q-1}\le h_{\cM}(r)\le c \cdot q^r$ for all $r\ge 0$, and $\cM$ contains all $\GF(q)$-representable matroids.
	\end{enumerate}
\end{theorem}

Since rank-$r$ order-$1$ cyclic Dowling geometries have $\binom{r+1}{2}$ elements and rank-$r$ projective geometries over $\GF(q)$ have $\frac{q^r-1}{q-1}$ elements, these members of $\cM$ give the lower bounds in (\ref{qgrt}) and (\ref{egrt}).

While this theorem describes the qualitative behaviour of extremal functions, it does not determine the functions themselves nor the densest matroids in the classes, and the lower bounds it gives are correct in general only up to a large constant factor. In other words, there is ample room to strengthen and refine the statement. For instance, \cite{Densest Exp} does this for classes satisfying (\ref{egrt}), showing that there are in fact integers $k,d \ge 0$ such that $h_{\cM}(r) = \frac{q^{r+k}-1}{q-1} - qd$, for all large $r$. 

We focus on classes satisfying (2), which are called \emph{quadratically dense} classes.
Our main result, Theorem~\ref{main corollary}, substantially refines outcome (\ref{qgrt}), and has Theorems~\ref{mainc}, \ref{mainr} and \ref{mainf1f2} as corollaries. 
The statement concerns $\Gamma$-frame and $\Gamma$-lift matroids, which are two classes of matroids that arise by labelling the edges of a graph with elements of a group $\Gamma$; we direct the reader to Section \ref{prelims} for the detailed definitions of these classes.

\begin{theorem} \label{main corollary}
If $\cM$ is a minor-closed class of matroids not containing all simple rank-2 matroids, then there is a finite group $\Gamma$ and a constant $c=c_{\cM}$ so that either 
\begin{itemize}
\item $h_{\cM}(r)\le c \cdot r$ for all $r\ge 0$, or

\item $\Gamma$ is nontrivial and $\cM$ contains all $\Gamma$-lift matroids, or

\item $|\Gamma|{r\choose 2}+r\le h_{\cM}(r)\le |\Gamma|{r\choose 2}+c \cdot r$ for all $r\ge 0$, and $\cM$ contains all $\Gamma$-frame matroids.
\end{itemize}
\end{theorem}

The proof of this result occupies Sections 4-8;  Sections 5-7, which form the bulk of the proof, are taken essentially verbatim from \cite{Walsh}.
In Sections \ref{exact chap} and \ref{proofs} we apply Theorem \ref{main corollary} to prove Theorems \ref{mainc} and \ref{mainr}.
\end{section}


\begin{section}{Corollaries} \label{results}
The classes of $\Gamma$-lift matroids and $\Gamma$-frame matroids have algebraic properties that allow us to use Theorem \ref{main corollary} to derive results for classes of representable matroids.

\begin{subsection}{Exact Results}
Theorems \ref{mainc} and \ref{mainf1f2} are both special cases of the following result, which we prove in Chapter \ref{proofs}.

\begin{theorem} \label{rep exact}
Let $\cF$ be a family of fields having no common finite subfield, and let $t\ge 1$ be an integer.
If $r$ is sufficiently large, and $M$ is a simple rank-$r$ matroid representable over all fields in $\cF$ and with no $U_{2,t+3}$-minor, then $|M| \le t{r\choose 2}+r$.
Moreover, if equality holds, then $M$ is isomorphic to a rank-$r$ cyclic Dowling geometry of order $t$. 
\end{theorem}

Whenever the size of the largest common subgroup (up to isomorphism) of the multiplicative groups of the fields in $\cF$ has size $t$, this bound is tight for cyclic Dowling geometries of order $t$.
There are two notable cases for which this occurs.
The first is Theorem \ref{mainc}, which was conjectured independently by Nelson \cite{Nelson} and Kapadia \cite{Kapadia}.
The second was conjectured by Geelen, Gerards, and Whittle in \cite{Structure Theory}, and is an important special case of Theorem \ref{mainf1f2}.

\begin{theorem} \label{ff1}
Let $\bF_1$ and $\bF_2$ be finite fields of different characteristic such that $|\bF_1|-1$ divides $|\bF_2|-1$.
If $r$ is sufficiently large, and $M$ is a simple rank-$r$ matroid representable over $\bF_1$ and $\bF_2$, then  $|M| \le (|\bF_1|-1){r\choose 2}+r$.
Moreover, if equality holds, then $M$ is isomorphic to a rank-$r$ cyclic Dowling geometry of order $|\bF_1|-1$.
\end{theorem}

When $\bF_1=\GF(2)$, this is the class of regular matroids, and the result was proved by Heller \cite{Heller}. 
When $\bF_1=\GF(3)$, the result follows from results of Kung \cite{dyadic1}, Kung and Oxley \cite{dyadic2}, and Whittle \cite{Whittle}.
This problem was previously open for all cases with $\bF_1\notin\{\GF(2),\GF(3)\}$.
\end{subsection}

\begin{subsection}{Approximate Results}
We will show that Theorem \ref{main corollary} implies the following general result,
which determines the correct leading coefficient of the extremal function of any quadratically dense class of matroids defined by representability over a family of fields.

\begin{theorem} \label{rep approx}
Let $\ell\ge 2$ be an integer, and let $\cF$ be a family of fields having no common finite subfield. 
Let $\alpha$ be the size of the largest common subgroup (up to isomorphism) of size less than $\ell$ of the multiplicative groups of the fields in $\cF$.
Then there is a constant $c$ so that if $M$ is a simple rank-$r$ matroid representable over all fields in $\cF$ and with no $U_{2,\ell+2}$-minor, then $|M| \le  \alpha{r\choose 2}+c\cdot r$.
\end{theorem}

Since a rank-$r$ cyclic Dowling geometry of order $\alpha$ has size $\alpha{r\choose 2}+r$, this bound is tight up to the constant $c$.
Note that Theorem \ref{mainr} is a special case of Theorem \ref{rep approx}; we highlight one other interesting special case, which was conjectured in \cite{Structure Theory}.

\begin{theorem} \label{ff2}
Let $\bF_1$ and $\bF_2$ be finite fields of different characteristic, and let $\alpha=\gcd(|\bF_1|-1,|\bF_2|-1)$. Then there is a constant $c$ so that if $M$ is a simple rank-$r$ matroid representable over $\bF_1$ and $\bF_2$, then $|M| \le \alpha{r\choose 2} + c\cdot r$. 
\end{theorem}

The previous best upper bound on the leading coefficient was $q^{2^{q-1}-1}-q^{2^{q-1}-2}$, where $q=\min(|\bF_1|,|\bF_2|)$, which was proved by Kung \cite{Kung2}.

We also find the extremal functions for a large family of minor-closed classes of  $\GF(p^k)$-representable matroids, up to a linear error term.
For each integer $m\ge 1$, we write $\bZ_m$ for the cyclic group on $m$ elements.

\begin{theorem} \label{GF(p^k)}
Let $p$ be a prime, let $k\ge 1$ be an integer, and let $N$ be a $\bZ_p$-lift matroid.
Then there is a constant $c$ so that the class $\cM$ of matroids representable over $\GF(p^k)$ and with no $N$-minor satisfies $h_{\cM}(r)\le (p^k-1){r\choose 2} + c\cdot r$, for all $r\ge 0$.
If, in addition, $N$ is not a $(\bZ_{p^k-1})$-frame matroid, then $\cM$ also satisfies $h_{\cM}(r)\ge (p^k-1){r\choose 2}+r$ for all $r\ge 0$.
\end{theorem}

Our final approximate result, which we prove at the end of Section \ref{main proof chap}, can be used to show that certain integer programs can be solved efficiently; see \cite{Paat} and \cite{Glanzer} for a detailed discussion.
Given a positive integer $\Delta$, we say that an integer matrix $A$ is \emph{$\Delta$-modular} if the determinant of each $\rank(A)\times \rank(A)$ submatrix has absolute value at most $\Delta$. 

\begin{theorem} \label{IP}
For each positive integer $\Delta$, there is a constant $c$ so that if $A$ is a rank-$m$ $\Delta$-modular matrix, then the number of distinct rows of $A$ is at most $m^2+c\cdot m$.
\end{theorem}

The previous best upper bound was $\Delta^{2+2\log_2\log_2(\Delta)}\cdot m^2+1$, due to Glanzer et al. \cite{Glanzer}.

\end{subsection}


\begin{subsection}{Algebraic Matroids}
We will extend all of these results to classes of algebraic matroids, with the exceptions of Theorems \ref{rep approx}, \ref{ff2}, and \ref{IP}.
Given a field $\bF$, an extension field $\bK$ of $\bF$, and a finite set $E\subseteq \bK$, the collection of subsets of $E$ that are algebraically independent over $\bF$ is the collection of independent sets (sets which have equal rank and size) of a matroid on $E$.
A matroid which arises in this way is \emph{algebraic} over $\bF$.

Piff \cite{Piff} proved that if a matroid is representable over a field $\bF$, then it is algebraic over $\bF$.
Conversely, Ingleton \cite{Ingleton} proved that if $\bF$ is a field of characteristic zero, then a matroid is algebraic over $\bF$ if and only if it is $\bF$-representable.
By Piff's theorem, the following result implies Theorem \ref{rep exact}.

\begin{theorem} \label{alg exact}
Let $\cF$ be a family of fields having no common finite subfield, and let $t\ge 1$ be an integer.
If $r$ is sufficiently large, and $M$ is a simple rank-$r$ matroid algebraic over all fields in $\cF$ and with no $U_{2,t+3}$-minor, then $|M| \le t{r\choose 2}+r$.
Moreover, if equality holds, then $M$ is isomorphic to a rank-$r$ cyclic Dowling geometry of order $t$.
\end{theorem}

This implies that Theorem \ref{ff1} extends to algebraic matroids, although we must explicitly exclude the relevant rank-2 uniform matroid, since every uniform matroid is algebraic over every field.

\begin{theorem} \label{alg ff1}
Let $\bF_1$ and $\bF_2$ be finite fields of different characteristic such that $|\bF_1|-1$ divides $|\bF_2|-1$.
If $r$ is sufficiently large, and $M$ is a simple rank-$r$ matroid algebraic over $\bF_1$ and $\bF_2$ and with no $U_{2,|\bF_1|+2}$-minor, then $|M|\le (|\bF_1|-1){r\choose 2}+r$.
Moreover, if equality holds, then $M$ is isomorphic to a rank-$r$ cyclic Dowling geometry of order $|\bF_1|-1$.
\end{theorem}

It is unclear whether or not this result is equivalent to Theorem \ref{ff1}; it may be strictly stronger for some choices of $\bF_1$ and $\bF_2$.
We also extend Theorem \ref{GF(p^k)} to algebraic matroids.

\begin{theorem} \label{alg GF(p^k)}
Let $p$ be a prime, let $k\ge 1$ be an integer, and let $N$ be a $\bZ_p$-lift matroid.
Then there is a constant $c$ so that the class $\cM$ of matroids algebraic over $\GF(p^k)$ and with no $N$-minor and no $U_{2,p^k+2}$-minor satisfies $h_{\cM}(r)\le (p^k-1){r\choose 2} + c\cdot r$, for all $r\ge 0$.
If, in addition, $N$ is not a $(\bZ_{p^k-1})$-frame matroid, then $\cM$ also satisfies $h_{\cM}(r)\ge (p^k-1){r\choose 2}+r$ for all $r\ge 0$.
\end{theorem}

Again, it is unclear whether this theorem is equivalent to Theorem \ref{GF(p^k)}; it may be strictly stronger for some choices of $p$, $k$, and $N$.

\end{subsection}

\end{section}


\begin{section}{Preliminaries} \label{prelims}
We use the notation of Oxley \cite{Oxley}.
A rank-1 flat is a \emph{point} and a rank-2 flat is a \emph{line}. 
A \emph{long line} is a line that contains at least three points.
We write $|M|$ for $|E(M)|$ and $\elem(M)$ for $|\si(M)|$, the number of points of a matroid $M$.
For an integer $\ell\ge 2$, we write $\cU(\ell)$ for the class of matroids with no $U_{2,\ell+2}$-minor.
In this section we will discuss the basics of distance, biased graphs, and connectivity, and then state a structural result (Theorem \ref{spanning clique}) which is vital for the proof of Theorem \ref{main corollary}.

We first highlight one component of the Growth Rate Theorem which will be particularly useful \cite{PG}. 

\begin{theorem}[Geelen, Whittle] \label{clique minor}
There is a function $\alpha_{\ref{clique minor}}\colon \bZ^2\to \bZ$ so that for all integers $\ell,t\ge 2$, if $M\in \cU(\ell)$ satisfies $\elem(M)>\alpha_{\ref{clique minor}}(\ell,t) \cdot r(M)$, then $M$ has an $M(K_{t+1})$-minor.
\end{theorem}

\begin{subsection}{Distance}
An \emph{extension} of a matroid $M$ is a matroid $M^+$ with ground set $E(M)\cup\{e\}$ for some $e\notin E(M)$ such that $M=M^+\del e$.
If $e$ is in a nontrivial parallel class of $M^+$, then $M^+$ is a \emph{parallel extension} of $M$.
We say that $M^+$ is a \emph{trivial extension} of $M$ if $e$ is a loop or coloop of $M^+$, or $e$ is parallel to an element of $M$; otherwise $M^+$ is a \emph{nontrivial extension} of $M$.
A \emph{projection} of $M$ is a matroid of the form $M^+/e$, where $M^+$ is an extension of $M$ by $e$.
A \emph{lift} of $M$ is a matroid $N$ so that $M$ is a projection of $N$.

For matroids $M$ and $N$ with the same ground set, the \emph{distance} between $M$ and $N$, denoted $\dist(M,N)$, is the smallest integer $k$ so that $N$ can be obtained from $M$ by a sequence of $k$ operations, each of which is a projection or lift.
Note that $\dist(M,N)=\dist(N,M)$, and that $\dist(M,N)\le r(M)+r(N)$, since the matroid on ground set $E$ consisting of only loops can be obtained from $M$ by $r(M)$ projections.
\end{subsection}


\begin{subsection}{Frame Matroids}
A matroid $M$ is \emph{framed} by $B$, and $B$ is a \emph{frame} for $M$, if $B$ is a basis of $M$ and each element of $M$ is spanned by a subset of $B$ with at most two elements.
A \emph{frame matroid} is a matroid of the form $M\del B$ where $M$ is a matroid framed by $B$.
We write $\cF$ for the class of frame matroids.

An important special case of a frame matroid is a \emph{$B$-clique}, which is a matroid framed by $B$ so that each pair of elements of $B$ is contained in a common long line.
For example, $M(K_n)$ is a $B$-clique for any spanning star $B$ of $K_n$. 
Since each graphic matroid is a restriction of $M(K_n)$ for some $n\ge 3$, each graphic matroid is a frame matroid.

The class of frame matroids is minor-closed; if $M$ is framed by $B$ and $e\in E(M)$, then $M/e$ is framed by any spanning subset of $B-\{e\}$.
This implies that if $M$ is a $B$-clique and $N$ is a contract-minor of $M$, then there is some $B'\subseteq B$ so that $N$ is a $B'$-clique.

The following lemma states some basic properties of frame matroids.

\begin{lemma} \label{frame obvious}
Let $M$ be a matroid framed by $B$. Let $b\in B$, and let $f$ be an element of $M$ which is not parallel to any element of $B$. Then
\begin{enumerate}[(i)]
\item each line of $M$ of length at least four spans two elements of $B$,
\item $f$ is on at most one line of length at least four, and 
\item each long line of $M$ through $b$ contains some $b'\in B-\{b\}$.  
\end{enumerate}
\end{lemma}

We will be interested in the interplay between frame matroids and spikes.
A \emph{spike} is a simple matroid $S$ with an element $t$, called a \emph{tip}, so that $\si(S/t)$ is a circuit, and each parallel class of $S/t$ has size two.
Note that Lemma \ref{frame obvious} (iii) implies that $b$ is not the tip of a spike restriction of $M$, since $B-\{b\}$ is independent in $M/b$.
\end{subsection}


\begin{subsection}{Biased Graphs}
Zaslavsky \cite{Z3} showed that frame matroids can be encoded by graphs.
A \emph{theta graph} consists of two distinct vertices $x$ and $y$, and three internally disjoint paths from $x$ to $y$.
A set $\cB$ of cycles of a graph $G$ satisfies the \emph{theta property} if no theta subgraph of $G$ contains exactly two cycles in $\cB$. 
A \emph{biased graph} is a pair $(G,\cB)$ where $G$ is a graph, and $\cB$ is a collection of cycles of $G$ which satisfies the theta property. 
The cycles in $\cB$ are \emph{balanced}, and the cycles not in $\cB$ are \emph{unbalanced}.

Biased graphs were first described by Zaslavsky in \cite{Z1}, and in \cite{Z2} he defined two matroids associated with a given biased graph $(G,\cB)$.
The \emph{frame matroid} $M(G,\cB)$ of $(G,\cB)$ is the matroid with ground set $E(G)$ so that $C\subseteq E(G)$ is a circuit of $M(G,\cB)$ if and only if the edges of $C$ form 
\begin{itemize}
\item a balanced cycle,
\item two vertex-disjoint unbalanced cycles with a path between them,
\item two unbalanced cycles which share a single vertex, or
\item a theta graph with all cycles unbalanced.
\end{itemize}
If $\cB$ is the set of all cycles of $G$, then $M(G,\cB)\cong M(G)$, the cycle matroid of $G$.
The matroid $M(G,\cB)$ is in fact a frame matroid; if we add an unbalanced loop to each vertex of $G$, then this set of elements frames the resulting matroid.

The \emph{lift matroid} $L(G,\cB)$ of $(G,\cB)$ is the matroid with ground set $E(G)$ so that $C\subseteq E(G)$ is a circuit of $L(G,\cB)$ if and only if the edges of $C$ form 
\begin{itemize}
\item a balanced cycle,
\item two vertex-disjoint unbalanced cycles,
\item two unbalanced cycles which share a single vertex, or
\item a theta graph with all cycles unbalanced.
\end{itemize}
Again, if $\cB$ is the set of all cycles of $G$, then $L(G,\cB)\cong M(G)$.
The lift matroid $L(G,\cB)$ of a biased graph $(G,\cB)$ is always a lift of the cycle matroid of $G$; we define the \emph{extended lift matroid} $L^+(G,\cB)$ to be the matroid which is an extension of $L(G,\cB)$ by an element $e_0$ so that $L^+(G,\cB)/e_0\cong M(G)$.
This terminology differs from that of Zaslavsky \cite{Z2}, who calls this the \emph{complete lift matroid}.

A natural family of biased graphs arises from graphs whose edges are labeled by elements of a group.
We define group-labeled graphs using the notation of \cite{DeVos}.
A \emph{group-labeling} of a graph $G$ consists of an orientation of the edges of $G$ and a function $\phi\colon E(G)\to \Gamma$ for some (multiplicative) group $\Gamma$.
For each walk $W$ on $G$ with edge sequence $e_1,e_2,\dots,e_k$, define $\epsilon_i(W)$ by
$$\epsilon_i(W)=\begin{cases}
1 & \text{ if $e_i$ is traversed forward in $W$},  \\
-1 & \text{ if $e_i$ is traversed backward in $W$}, 
\end{cases}$$
and define $\phi(W)=\prod_{i=1}^{k}\phi(e_i)^{\epsilon_i(W)}$.
We define $\cB_{\phi}$ to be the set of all cycles $C$ for which there is some simple closed walk $W$ around $C$ so that $\phi(W)=1$.
This set $\cB_{\phi}$ of cycles satisfies the theta property, so $(G,\cB_{\phi})$ is a biased graph.
The frame matroid of $(G,\cB_{\phi})$ is a \emph{$\Gamma$-frame matroid}, and the lift matroid of $(G,\cB_{\phi})$ is a \emph{$\Gamma$-lift matroid}.
\end{subsection}


\begin{subsection}{Dowling Geometries}
We have already seen cyclic Dowling geometries; we now define Dowling geometries over any finite group $\Gamma$, as a special class of $\Gamma$-frame matroids.
Dowling geometries were first introduced by Dowling in \cite{Dowling}, although the definition presented here is due to Zaslavsky \cite{Z2}.

Let $\Gamma$ be a finite group of size at least two with identity $1$, and let $G$ be a graph with $k\ge 3$ vertices, a single loop at each vertex with label in $\Gamma-\{1\}$, and exactly $|\Gamma|$ parallel edges between each pair of vertices, so that between each pair of vertices each edge has the same orientation and each element of $\Gamma$ appears as a label.
The frame matroid obtained from this group-labeled graph is the rank-$k$ \emph{Dowling geometry} over the group $\Gamma$, and is denoted $\DG(k,\Gamma)$.
If $|\Gamma|=1$, then we define $\DG(k,\Gamma)$ to be the frame matroid constructed from the graph $K_{k}$ with a loop at each vertex so that all nonloop cycles are unbalanced.
It is not hard to show that $\DG(k,\{1\})\cong M(K_{k+1})$, which implies that each graphic matroid is the frame matroid of a $\{1\}$-labeled graph.
Note that  $|\DG(k,\Gamma)|=|\Gamma|{k\choose 2}+k$, since $G$ has $k$ loops and $|\Gamma|$ edges between each pair of vertices.

The class of Dowling geometries with group $\Gamma$ has the useful property that $\si(\DG(k,\Gamma)/e)\cong \DG(k-1,\Gamma)$ for each $e\in E(\DG(k,\Gamma))$.
It is also not hard to see that each simple rank-$k$ $\Gamma$-frame matroid is a restriction of $\DG(k,\Gamma)$, which is analogous to the fact that each simple rank-$k$ $\GF(q)$-representable matroid is a restriction of $\PG(k-1,q)$.

In \cite{Dowling} Dowling proved that $\DG(k,\Gamma)$ and $\DG(k,\Gamma')$ are isomorphic matroids if and only if $\Gamma$ and $\Gamma'$ are isomorphic groups.
He also proved the following beautiful theorem, which we use to prove applications of our main result.

\begin{theorem} \label{Dowling fields}
The Dowling geometry $\DG(k,\Gamma)$ is representable over a field $\bF$ if and only if $\Gamma$ is isomorphic to a subgroup of the multiplicative group of $\bF$.
\end{theorem}

Dowling geometries also have the attractive property that if $\Gamma'$ is a subgroup of $\Gamma$, then $\DG(k,\Gamma')$ is a restriction of $\DG(k,\Gamma)$.
In particular, since $\DG(k,\{1\})\cong M(K_{k+1})$ and every group has a subgroup of size one, every Dowling geometry has a spanning clique restriction.
As we shall see in Section \ref{frame chap}, Dowling geometries are occasionally easier to work with if we delete the frame $B$.
We define $\DG^-(k,\Gamma)= \DG(k,\Gamma)\del B$, where $B$ is a frame of $\DG(k,\Gamma)$. 
\end{subsection}


\begin{subsection}{Connectivity}
Let $M=(E,r)$ be a matroid.
The \emph{connectivity function} of $M$ is the function $\lambda_M\colon 2^{E}\to \bZ$ defined by $\lambda_M(X)=r(X)+r(E-X)-r(M)$ for all $X\subseteq E$.

For disjoint sets $A,B\subseteq E$ we write $\kappa_M(A,B)=\min(\lambda(Z)\colon A\subseteq Z \text{ and } B\subseteq E-Z)$.
If $B\subseteq C\subseteq E-A$, then it is clearly the case that $\kappa_M(A,B)\le \kappa_M(A,C)$.
It is also not hard to check that if $N$ is a minor of $M$ for which $A\cup B\subseteq E(N)$, then $\kappa_N(A,B)\le \kappa_M(A,B)$.
We will use the following strengthening of Tutte's linking theorem \cite{TLT}, which was proved by Geelen, Gerards, and Whittle in \cite{STLT}.

\begin{theorem} \label{STLT}
Let $X$ and $Y$ be disjoint sets of elements of a matroid $M$. Then $M$ has a minor $N$ with ground set $X\cup Y$ such that $\kappa_N(X,Y)=\kappa_M(X,Y)$, while $N|X=M|X$ and $N|Y=M|Y$.
\end{theorem}

A \emph{vertical $k$-separation} of $M$ is a partition $(X,Y)$ of $E$ so that $r(X)+r(Y)-r(M)<k$ and $\min(r(X),r(Y))\ge k$.
We say that $M$ is \emph{vertically $k$-connected} if $M$ has no vertical $j$-separation with $j<k$.
Equivalently, $M$ is vertically $k$-connected if there is no partition $(X,Y)$ of $E$ so that $r(X)+r(Y)-r(M)<k-1$ and $\max(r(X),r(Y))<r(M)$.

We also need some notation to capture how much subsets of $E$ interact in $M$.
For sets $A,B\subseteq E$, the \emph{local connectivity} between $A$ and $B$ is $\sqcap_M(A,B)=r_M(A)+r_M(B)-r_M(A\cup B)$.
Note that if $A$ and $B$ are disjoint, then $\sqcap_M(A,B)=r_M(A)-r_{M/B}(A)=r_M(B)-r_{M/A}(B)$.
We say that sets $A$ and $B$ are \emph{skew} in $M$ if $\sqcap_M(A,B)=0$; roughly speaking, this means that $A$ and $B$ do not interact at all in $M$.
More generally, sets $A_1,\dots,A_k\subseteq E$ are \emph{mutually skew} in $M$ if $r_M(A_1\cup\dots\cup A_k)=\sum_i r_M(A_i)$.
\end{subsection}


\begin{subsection}{Matroids with a Spanning Clique}
The following deep result from \cite{Spanning Clique} is vital for the proof of Theorem \ref{main corollary}.
It shows that every matroid with a spanning clique restriction in a quadratically dense class is a bounded distance from a frame matroid. 

\begin{theorem}[Geelen, Nelson] \label{spanning clique}
There is a function $h_{\ref{spanning clique}}\colon \bZ^2\to \bZ$ so that for all integers $\ell,n\ge 2$, if $M\in \cU(\ell)$ has a spanning $B$-clique restriction and no rank-$n$ projective-geometry minor,  then there is a set $\hat B\subseteq B$ and a $\hat B$-clique $N$ such that $\dist(M,N)\le h_{\ref{spanning clique}}(\ell,n)$.
Moreover, there are disjoint sets $C_1,C_2\subseteq E(M)$ with $r_M(C_1\cup C_2)\le h_{\ref{spanning clique}}(\ell,n)$ such that 
\begin{itemize}
\item $\si(N)$ is isomorphic to a restriction of $M/C_1$, and
\item for all $X\subseteq E(M)-(C_1\cup C_2)$, if $(M/(C_1\cup C_2))|X$ is simple, then $N|X=(M/C_1)|X$. 
\end{itemize}
\end{theorem}

The second condition implies that $\elem(N)\ge \elem(M/(C_1\cup C_2))$, by taking $(M/(C_1\cup C_2))|X$ to be a simplification of $M/(C_1\cup C_2)$.

\end{subsection}

\end{section}


\begin{section}{$\Gamma$-Lift Matroids} \label{lifts}
In this section we define an analogue of Dowling geometries for $\Gamma$-lift matroids, and show that they share several natural properties with Dowling geometries.
We then prove that these matroids arise as minors of large-rank  lift matroids of `complete' biased graphs, which we call doubled cliques.

Let $\Gamma$ be a finite group of size at least two, and let $G$ be a graph with $k\ge 3$ vertices, and exactly $|\Gamma|$ parallel edges between each pair of vertices, so that between each pair of vertices, each edge has the same orientation, and each element of $\Gamma$ appears as a label.
We say that the lift matroid obtained from this group-labeled graph is the rank-$k$ \emph{lift geometry} over $\Gamma$, and write $\LG(k,\Gamma)$.
We write $\LG^+(k,\Gamma)$ for the corresponding extended lift matroid, and we use the convention that $e_0$ is the element so that $\LG^+(k,\Gamma)\del e_0=\LG(k,\Gamma)$.
We define $\LG(k,\{1\})\cong M(K_{k+1})$. 

If $\Gamma'$ is a nontrivial subgroup of $\Gamma$, then clearly $\LG(k,\Gamma')$ is a restriction of $\LG(k,\Gamma)$, by only considering edges labeled by elements of $\Gamma'$.
Also, if $N$ is a restriction of $\LG^+(k,\Gamma)$ so that $e_0\in E(N)$, each line of $N$ through $e_0$ has length $|\Gamma|+1$, and $\si(N/e_0)\cong M(K_m)$, then $N\cong \LG^+(m,\Gamma)$.
This implies that for each integer $k\ge 4$, each element $e$ of $\LG^+(k,\Gamma)$ other than $e_0$ satisfies $\si(\LG^+(k,\Gamma)/e)\cong \LG^+(k-1,\Gamma)$.

We will need the following lemma from \cite{Dowling}.


\begin{lemma}[Dowling] \label{groups}
If $\Gamma$ and $\Gamma'$ are (additive) groups, and there exist bijections $\phi_i\colon \Gamma\to\Gamma'$ for $i=1,2,3$ such that, for all $\alpha,\beta\in\Gamma$,
\begin{align*}
\phi_3(\alpha+\beta)=\phi_1(\alpha)+\phi_2(\beta),
\end{align*}
then $\Gamma\cong\Gamma'$.
\end{lemma}


We first prove that the lift geometry encodes the structure of the group, following the proof of Theorem 8 from \cite{Dowling}.

\begin{theorem}
If $\LG^+(k,\Gamma)$ and $\LG^+(k,\Gamma')$ are isomorphic matroids, then $\Gamma$ and $\Gamma'$ are isomorphic groups.
\end{theorem}
\begin{proof}
We will write $\Gamma$ and $\Gamma'$ as additive groups with identity $0$.
Clearly $|\Gamma|=|\Gamma'|$.
The theorem holds if $|\Gamma|<4$, so we may assume that $|\Gamma|\ge 4$.
Let $\sigma$ be an isomorphism from $\LG^+(k,\Gamma)$ to $\LG^+(k,\Gamma')$.
Then $\sigma(e_0)=e_0'$, since these are the only elements on more than one line of length at least five.

Let $N$ be an $\LG^+(3,\Gamma)$-restriction of $\LG^+(k,\Gamma)$.
Then $N$ consists of three lines through $e_0$, and $\sigma(E(N))$ thus consists of three lines through $e_0'$.
Since this set of elements has rank three in $\LG^+(k,\Gamma')$, this implies that the image of $E(N)$ under $\sigma$ is an $\LG^+(3,\Gamma')$-restriction $N'$ of $\LG^+(k,\Gamma)$.
Let $G$ be a $\Gamma$-labeled graph on vertex set $\{1,2,3\}$ whose associated extended lift matroid is $N$, and whose edges are oriented from $1$ to $2$, from $2$ to $3$, and from $3$ to $1$.
Then a triangle of $G$ is a circuit of $N$ if and only if the sum of their labels in $\Gamma$ is $0$.
Let $G'$ be a similarly-defined group-labeled graph for $N'$.

Let $L_1,L_2,L_3$ denote the points of $N/e_0$, and for each $i\in \{1,2,3\}$, let $L_i'=\sigma(L_i)$.
For each $i\in \{1,2,3\}$, we may regard the restriction of $\sigma$ to $L_i$ as a bijection $\sigma_i$ from $\Gamma$ to $\Gamma'$, since  $L_i$ is fully labeled by $\Gamma$ and $L_i'$ is fully labeled by $\Gamma'$.
Since $\sigma$ is an isomorphism from $N$ to $N'$, for all $\alpha_1,\alpha_2,\alpha_3\in \Gamma$ we have 
$$\alpha_1+\alpha_2+\alpha_3=0 \iff \sigma_1(\alpha_1)+\sigma_2(\alpha_2)+\sigma_3(\alpha_3)=0,$$
which implies that $$\sigma_3(-(\alpha_1+\alpha_2))=-(\sigma_1(\alpha_1)+\sigma_2(\alpha_2)).$$
Let $\tau_3\colon \Gamma\to \Gamma'$ be the bijection $\tau_3(\alpha)=-\sigma_3(-\alpha)$.
Then for all $\alpha_1,\alpha_2\in\Gamma$,
\begin{align*}
\tau_3(\alpha_2+\alpha_1)&=-\sigma_3(-(\alpha_2+\alpha_1))\\
&=(\sigma_1(\alpha_1)+\sigma_2(\alpha_2))\\
&=\sigma_2(\alpha_2)+\sigma_1(\alpha_1).
\end{align*}
Now, Lemma \ref{groups} applied with $(\sigma_1,\sigma_2,\sigma_3)=(\sigma_2,\sigma_1,\tau_3)$ shows that $\Gamma\cong\Gamma'$, as desired.
\end{proof}


\begin{subsection}{Representability of Lift Geometries}
We now determine when $\LG^+(k,\Gamma)$ is representable over a field $\bF$.
Zaslavsky proved in \cite{Z4} that if $\Gamma$ is isomorphic to a subgroup of the additive group $\bF^+$ of $\bF$, then $\LG^+(k,\Gamma)$ is representable over $\bF$.
He showed that the matrix with columns consisting of the unit vector $e_i$ for each $i\in [k]$, and $-e_i+e_j+\alpha e_1$ for each $2\le i<j\le k$ and $\alpha\in \bF^+$, is a representation of $\LG^+(k,\bF^+)$ over $\bF$.
This implies that every rank-$k$ $\bF$-representable extended lift matroid of a biased graph is a restriction of $\LG^+(k,\bF^+)$.

We first prove a theorem which, combined with Zaslavsky's result, characterizes the fields over which $\LG^+(k,\Gamma)$ is representable.
This is the analogue of Theorem \ref{Dowling fields} for lift geometries.
We follow the proof of Theorem 9 in \cite{Dowling}.

\begin{theorem} \label{LG rep}
Let $\Gamma$ be a finite group, and let $\bF$ be a field with additive group $\bF^+$.
Then $\LG^+(k,\Gamma)$ is representable over $\bF$ if and only if $\Gamma$ is isomorphic to a subgroup of $\bF^+$.
\end{theorem}
\begin{proof}
We will prove the forwards direction, while the backwards direction is due to Zaslavsky \cite{Z4}.
We write $\Gamma$ as an additive group with identity $0$.
Since each element $e$ of $\LG^+(k,\Gamma)$ other than $e_0$ satisfies $\si(\LG^+(k,\Gamma)/e)\cong \LG^+(k-1,\Gamma)$, and the class of $\bF$-representable matroids is minor-closed, we may assume that $k=3$.
Let $Q$ be a $\Gamma$-labeled graph on vertex set $\{1,2,3\}$ whose associated extended lift matroid is $\LG^+(3,\Gamma)$, and whose edges are oriented from $1$ to $2$, from $2$ to $3$, and from $3$ to $1$.
Then a triangle of $Q$ is a circuit of $\LG^+(3,\Gamma)$ if and only if the sum of the edge-labels is $0$.

Let $L_1,L_2,L_3$ denote the three points of $\LG^+(3,\Gamma)/e_0$.
By row and column scaling, we may consider a representation of $\LG^+(3,\Gamma)$ over $\bF$ so that $e_0$ is represented by $[1,0,0]^T$, each element of $L_1-\{e_0\}$ is represented by a vector of the form $[a_1,1,0]^T$, each element of $L_2-\{e_0\}$ is represented by a vector of the form $[a_2,0,1]^T$, and each element of $L_3-\{e_0\}$ is represented by a vector of the form $[a_3,-1,-1]^T$, where $a_1,a_2,a_3\in \bF$.
This representation naturally defines a function $\sigma\colon E(\LG^+(3,\Gamma))\to \bF\times \bF\times \bF$.
Let 
\begin{align*}
H_1&=\{a_1\colon (a_1,1,0)\in \sigma(L_1)\},\\
H_2&=\{a_2\colon (a_2,0,1)\in \sigma(L_2)\},\\
H_3&=\{a_3\colon (a_3,-1,-1)\in \sigma(L_3)\}.
\end{align*}
For each $i\in \{1,2,3\}$, we define a bijection $\phi_i\colon \Gamma\to H_i$ as follows: given $\alpha\in\Gamma$, let $\alpha_i$ denote the element on $L_i$ labeled by $\alpha$ in $Q$, and define $\phi_i(\alpha)$ to be the first coordinate of the vector $\sigma(\alpha_i)$.

Since these vectors form a representation of $\LG^+(3,\Gamma)$, this implies that all elements $\alpha_1,\alpha_2,\alpha_3\in \Gamma$ satisfy $\phi_1(\alpha_1)+\phi_2(\alpha_2)+\phi_3(\alpha_3)=0$ (over $\bF$) if and only if $\alpha_1+\alpha_2+\alpha_3=0$ (over $\Gamma$).
Thus, all $\alpha_1,\alpha_2\in\Gamma$ satisfy 
\begin{align}
\phi_3(-(\alpha_1+\alpha_2))=-(\phi_1(\alpha_1)+\phi_2(\alpha_2)).
\end{align}
Since interchanging $\alpha_1$ and $\alpha_2$ does not affect the right-hand side since $\bF^+$ is Abelian, (1) implies that $\Gamma$ is Abelian.

\begin{claim}
For each $i\in \{1,2,3\}$, there is a subgroup $G_i$ of $\bF^+$ and an element $c_i\in H_i$ so that $H_i=G_i+c_i$ is a coset of $G_i$.
\end{claim}
\begin{proof}
We show this for $i=1$; the other cases are identical.
Since $\LG^+(3,\Gamma)$ has no $U_{2,|\Gamma|+2}$-minor, for any permutation $(i,j,k)$ of $(1,2,3)$, 
each element $\alpha_i\in L_i$ if and only if there exist $\alpha_j\in L_j$ and $\alpha_k\in L_k$ so that $\{\alpha_i,\alpha_j,\alpha_k\}$ is a line of $\LG^+(3,\Gamma)$.
This implies that $a_i\in H_i$ if and only if there exist $a_j\in H_j$ and $a_k\in H_k$ so that $a_i=-a_j-a_k$.

Let $a_1,b_1,c_1$ be any three (not necessarily distinct) elements of $H_1$.
Choose any $a_2\in H_2$.
Then
\begin{align}
&a_2\in H_2 \textrm{ and } a_1\in H_1 \implies -a_2-a_1\in H_3,\\
& -a_2-a_1\in H_3 \textrm{ and } c_1\in H_1 \implies a_2+a_1-c_1\in H_2,\\
& a_2+a_1-c_1\in H_2 \textrm{ and } b_1\in H_1 \implies c_1-a_2-a_1-b_1\in H_3,\\
& c_1-a_2-a_1-b_1\in H_3 \textrm{ and } a_2\in H_2 \implies a_1+b_1-c_1\in H_1.
\end{align}
Thus, 
\begin{align}
a_1,b_1,c_1\in H_1 \implies a_1+b_1-c_1\in H_1.
\end{align}
Fix $c_1\in H_1$, and define $$G_1=\{a_1-c_1\colon a_1\in H_1\}.$$
We will show that $G_1$ is a group.
Clearly $0\in G_1$, since $c_1\in H_1$.
If $a_1-c_1\in G_1$, then (6) with $(a_1,b_1,c_1)=(c_1,c_1,a_1)$ implies that $c_1+c_1-a_1\in H_1$, so $c_1-a_1\in G_1$, and so $G_1$ is closed under inverses.
Finally, if $a_1-c_1\in G_1$ and $b_1-c_1\in G_1$, then (6) implies that $a_1+b_1-c_1\in H_1$, so $(a_1-c_1)+(b_1-c_1)\in G_1$.
Thus, $G_1$ is a subgroup of $\bF^+$, and $H_1=G_1+c_1$ is a coset of $G_1$.
\end{proof}

These groups are in fact equal.
We define $\lambda_0=c_1+c_2+c_3$, for convenience.

\begin{claim}
$G_1=G_2=G_3$, and $\lambda_0\in G_1$.
\end{claim}
\begin{proof}
We will show that $G_1=G_3$; the other cases are identical.
Let $d_1\in G_1$, so $d_1+c_1\in H_1$. 
Then since $c_2\in H_2$, we have $-d_1-c_1-c_2\in H_3$, so $-d_1-(c_1+c_2+c_3)=-d_1-\lambda_0\in G_3$.
When $d_1=0$, this implies that $-\lambda_0\in G_3$, so $\lambda_0\in G_3$.
Then $-d_1-\lambda_0\in G_3$ implies that $-d_1\in G_3$, so $d_1\in G_3$.
Thus, $G_1\subseteq G_3$, and the same argument implies that $G_3\subseteq G_1$.
\end{proof}

Let $G=G_1=G_2=G_3$.
We use $G$ to construct a different representation of $\LG^+(k,\Gamma)$. 
Let $\mu_1,\mu_2,\mu_3\in G$ be elements so that $\mu_1+\mu_2+\mu_3=\lambda_0$ (we may take $\mu_1=\mu_2=0$ and $\mu_3=\lambda_0$).
For each $i\in \{1,2,3\}$, define a bijection $f_i\colon H_i\to G$ by $f_i(a_i)=\mu_i+a_i-c_i$.
Then
\begin{align*}
f_1(a_1)+f_2(a_2)+f_3(a_3)=a_1+a_2+a_3,
\end{align*}
for all $a_1\in H_1$, $a_2\in H_2$, and $a_3\in H_3$.
This implies that for all $\alpha_1,\alpha_2,\alpha_3\in\Gamma$,
\begin{align*}
\alpha_1+\alpha_2+\alpha_3=0 &\iff \phi_1(\alpha_1)+\phi_2(\alpha_2)+\phi_3(\alpha_3)=0\\
&\iff f_1(\phi_1(\alpha_1))+f_2(\phi_2(\alpha_2))+f_3(\phi_3(\alpha_3))=0,
\end{align*}
where $\phi_1,\phi_2,\phi_3$ are the bijections satisfying (1).
Thus, we may assume that $H_1=H_2=H_3=G$.

Finally, we show that $G\cong \Gamma$.
Define a bijection $\tau\colon \Gamma\to G$ by $\tau(\alpha)=-\phi_3(-\alpha)$, for all $\alpha\in\Gamma$.
Then (1) implies that $\tau(\alpha_1+\alpha_2)=\phi_1(\alpha_1)+\phi_2(\alpha_2)$, for all $\alpha_1,\alpha_2\in\Gamma$.
By Lemma \ref{groups} with $(\phi_1,\phi_2,\phi_3)=(\phi_1,\phi_2,\tau)$, the groups $\Gamma$ and $G$ are isomorphic, and thus $\Gamma$ is isomorphic to a subgroup of $\bF^+$.
\end{proof}


Since fields of characteristic zero have no finite additive subgroup, and each nontrivial, finite additive subgroup of a field of characteristic $p$ is isomorphic to $\bZ_p^m$ for some $m\ge 1$, Theorem \ref{LG rep} has the following corollary.

\begin{corollary}
Let $k\ge 3$ be an integer, let $\Gamma$ be a nontrivial finite group, and let $\bF$ be a field. Then
\begin{enumerate}[(i)]
\item if $\bF$ has characteristic zero, then $\LG^+(k,\Gamma)$ is not $\bF$-representable,

\item if $\bF=\GF(p^m)$, then $\LG^+(k,\Gamma)$ is $\bF$-representable if and only if $\Gamma=\bZ_p^j$ for some $1\le j\le m$, and

\item if $\bF$ has characteristic $p$ and $\LG^+(k,\Gamma)$ is representable over $\bF$, then $\Gamma=\bZ_p^j$ for some $j\ge 1$.
\end{enumerate}
\end{corollary}

In particular, this implies that $\PG(k,p^j)$ has an $\LG^+(k,\bZ_p^j)$-restriction.

\end{subsection}

\begin{subsection}{Finding a Lift-Geometry Minor}
We need the following theorem about biased graphs, which we will use to identify $\LG^+(k,\Gamma)$ and $\DG(k,\Gamma)$.
The proof is based on Section 7 of \cite{Kahn/Kung}.
For each integer $k\ge 3$ and each finite group $\Gamma$, we write $(G,\cB)_{(k,\Gamma)}$ to denote the biased graph obtained from the $\Gamma$-labeled graph on $k$ vertices for which, between each pair of vertices there are $|\Gamma|$ edges, all oriented in the same way, and all with distinct labels.
This is the $\Gamma$-labeled graph used to define both $\DG(k,\Gamma)$ and $\LG(k,\Gamma)$.


\begin{theorem} \label{biased graphs}
Let $t\ge 1$ and $m\ge 4$ be integers, and let $(G,\cB)$ be a biased graph on $m$ vertices with $t$ parallel edges between each pair of vertices. If each cycle in $\cB$ has size at least three, and 
\begin{enumerate}[$(*)$]
\item for all distinct $b_1,b_2,b_3\in V(G)$, each edge $a$ with ends $b_1$ and $b_2$, and each edge $b$ with ends $b_1$ and $b_3$, there is an edge $c$ with ends $b_2$ and $b_3$ so that $\{a,b,c\}\in\cB$,
\end{enumerate}
then there is a group $\Gamma$ of size $t$ so that $(G,\cB)=(G,\cB)_{(m,\Gamma)}$.
\end{theorem}
\begin{proof}
Let $b_1,b_2,\dots,b_m$ be an ordering of $V(G)$.
For all $1\le i<j\le m$, let $L_{ij}$ denote the set of edges of $G$ with ends $b_i$ and $b_j$.
For all $1\le i<j<k\le m$, let $t_{ijk}=L_{ij}\cup L_{ik}\cup L_{jk}$. 
We will refer to $t_{ijk}$ as a \emph{facet} of $G$. 

Note that the element $c$ from $(*)$ is necessarily unique; otherwise $c$ is in a size-2 cycle in $\cB$, by the theta-property.
We will refer to 3-element cycles in $\cB$ as \emph{triangles}.

\setcounter{claim}{0}
\begin{claim}\label{useful}
For all $1\le i<j<k<h\le m$, if $\{a,b,c\}$, $\{c,d,e\}$, and $\{a,f,e\}$ are triangles on three distinct facets out of $t_{ijk}$, $t_{ijh}$, $t_{ikh}$, and $t_{jkh}$, then $\{b,d,f\}$ is a triangle on the fourth facet. 
\end{claim}
\begin{proof}
Without loss of generality, we may assume that $(i,j,k,h)=(1,2,3,4)$, and $\{a,b,c\}$, $\{c,d,e\}$ and $\{a,f,e\}$ are  triangles of $t_{123}$, $t_{134}$, and $t_{124}$, respectively. 
This implies that $a\in L_{12}$, $b\in L_{23}$, $c\in L_{13}$, $d\in L_{34}$, $e\in L_{14}$, and $f\in L_{24}$.
By the theta-property, the cycles $\{a,b,c\}$ and $\{a,f,e\}$ imply that $\{b,c,e,f\}\in\cB$.
Again by the theta-property, the cycles $\{b,c,e,f\}$ and $\{c,d,e\}$ imply that $\{b,d,f\}\in \cB$.
\end{proof}

Let $\Gamma$ be a set of size $t$, and let $\ep\in \Gamma$. 
We will label $E(G)$ by elements of $\Gamma$ according to three rules, which force $\epsilon\in\Gamma$ to be the group identity element.
Rules (A) and (B) look very similar, but we need both in order to show that $\ep$ commutes with all elements of $\Gamma$ once we define a group operation on $\Gamma$.

\begin{claim}
There exists a function $f\colon E(G)\to \Gamma$ so that for all $1\le i<j<k\le m$,
\begin{enumerate}[(A)]
\item if $\{a,b,c\}$ is a triangle in $t_{ijk}$ such that $b\in L_{jk}$ and $f(b)=\epsilon$, then $f(a)=f(c)$, and 

\item if $\{a,b,c\}$ is a triangle in $t_{ijk}$ such that $b\in L_{ij}$ and $f(b)=\epsilon$, then $f(a)=f(c)$, and  

\item $f$ restricted to $L_{i'j'}$ is a bijection for all $1\le i'<j'\le m$.
\end{enumerate}
\end{claim}
\begin{proof}
For each $2\le j\le m$, arbitrarily choose an element $a\in L_{1j}$ and set $f(a)=\epsilon$.
Arbitrarily assign labels to elements of $L_{12}$ so that $f$ restricted to $L_{12}$ is a bijection.
We use the following three steps to define $f$, relying on the fact that $\cB$ satisfies $(*)$.

\begin{enumerate}[(1)]
\item 
For each $2\le j<k\le m$, let $a\in L_{jk}$ be the element in a triangle with the elements of $L_{1j}$ and $L_{1k}$ labeled $\epsilon$, and set $f(a)=\epsilon$.
Note that this element $a\in L_{jk}$ exists since $\cB$ satisfies $(*)$.

\item 
For each $3\le k\le m$, let $a\in L_{1k}$ be the element in a triangle with the element of $L_{12}$ labeled $\alpha$ and the element of $L_{2k}$ labeled $\epsilon$, and set $f(a)=\alpha$.
This shows that $t_{12k}$ satisfies (A) for all $3\le k\le m$.

\item 
For each $2\le j<k\le m$, let $a\in L_{jk}$ be the element in a triangle with the element of $L_{1j}$ labeled $\epsilon$ and the element of $L_{1k}$ labeled $\alpha$, and set $f(a)=\alpha$.
This show that $t_{1jk}$ satisfies (B) for all $2\le j<k\le m$.
\end{enumerate}

We now have a function $f\colon E(G)\to \Gamma$, and we will show that it satisfies (A), (B), and (C).
If $f$ restricted to $L_{1k}$ is not a bijection for some $3\le k\le m$, then there is some $a\in L_{1k}$ and $e\in L_{2k}$ which are in triangles with two distinct elements of $L_{12}$, which contradicts the uniqueness of $c$ in $(*)$.
Similarly, if $f$ restricted to $L_{jk}$ is not a bijection for some $2\le j<k\le m$, then there is some $a\in L_{jk}$ and $e\in L_{1j}$ which are in triangles with two distinct elements of $L_{1k}$, which contradicts the uniqueness of $c$ in $(*)$.
Thus, $f$ satisfies (C).

We now prove a sequence of claims to show that $f$ satisfies (A) and (B).
For each $\alpha\in \Gamma$ and $1\le i<j\le m$, let $\alpha_{ij}$ denote the element in $L_{ij}$ such that $f(\alpha_{ij})=\alpha$.

\begin{enumerate}[(i)]
\item By (1) and  \ref{useful} with $(a,b,c,d,e,f)$ all labeled $\epsilon$, we find that $\{\epsilon_{ij},\ep_{ik},\ep_{jk}\}\in\cB$  for all $1\le i<j<k\le m$.

\item By (2), $f$ satisfies (A) in $t_{12k}$ for all $3\le k\le m$.

\item By (i), (ii) and  \ref{useful} with $(a,b,c,d,e,f)=(\alpha_{12},\alpha_{1j},\ep_{2j},\ep_{jk},\ep_{2k},\alpha_{1k})$ we find that $\{\alpha_{1j},\ep_{jk},\alpha_{1k}\}\in\cB$, so $f$ satisfies (A) in $t_{1jk}$ for all $2\le j<k\le m$.

\item By (3), $f$ satisfies (B) in $t_{1jk}$ for all $2\le j<k\le m$. 

\item By (i), (iv) and  \ref{useful} with $(a,b,c,d,e,f)=(\alpha_{1k},\alpha_{ik},\ep_{1i},\ep_{ij},\ep_{1j},\alpha_{jk})$ we find that $\{\ep_{ij},\alpha_{jk},\alpha_{ik}\}\in\cB$, so (B) holds for all $2\le i<j<k\le m$.

\item By (iv), (iii) and  \ref{useful} with $(a,b,c,d,e,f)=(\ep_{1i},\alpha_{ij},\alpha_{1j},\ep_{jk},\alpha_{1k},\alpha_{ik})$, we find that $\{\alpha_{ij},\ep_{jk},\alpha_{ik}\}\in\cB$, so (A) holds for all $2\le i<j<k\le m$. 
\end{enumerate}
Thus, $f$ satisfies (A) by (iii) and (vi) and $f$ satisfies (B) by (iv) and (v).
\end{proof}

Now for each facet $t_{ijk}$, we define a binary operation $\circ_{ijk}:\Gamma\times\Gamma\to\Gamma$ by $\circ_{ijk}(\alpha,\beta)=\gamma$ if $\{\alpha_{ij},\beta_{jk},\gamma_{ik}\}\in\cB$. 
These operations are well-defined since $\cB$ satisfies $(*)$.

\begin{claim}
$\circ_{ijk}=\circ_{i'j'k'}$ for all $1\le i<j<k\le m$ and $1\le i'<j'<k'\le m$.
\end{claim}
\begin{proof}
We will show that any two facets with two common indices have the same binary operation.
Without loss of generality, it suffices to show that $\circ_{123}=\circ_{124}$ and $\circ_{134}=\circ_{234}$ and $\circ_{124}=\circ_{134}$.
We first show that $\circ_{123}=\circ_{124}$. 
Let $\alpha,\beta\in \Gamma$.
By \ref{useful} with $(a,b,c,d,e,f)=(\beta_{23},(\alpha\circ_{123}\beta)_{13},\alpha_{12},(\alpha\circ_{124}\beta)_{14},\beta_{24},\ep_{34})$, we find that $\{\ep_{34},(\alpha\circ_{123}\beta)_{13},(\alpha\circ_{124}\beta)_{14}\}\in\cB$.
This shows that $\alpha\circ_{123}\beta=\alpha\circ_{124}\beta$ since rule (A) holds for $t_{134}$. 

Similarly,  \ref{useful} with $(a,b,c,d,e,f)=(\alpha_{13},\ep_{12},\alpha_{23},(\alpha\circ_{234}\beta)_{24},\beta_{34},(\alpha\circ_{134}\beta)_{14})$ shows that $\circ_{134}=\circ_{234}$.
Lastly,  \ref{useful} with $(a,b,c,d,e,f)=(\alpha_{12},\alpha_{13},\ep_{23},\beta_{34},\beta_{24},(\alpha\circ_{124}\beta)_{14})$ shows that $\circ_{124}=\circ_{134}$.
Thus, any two facets with two common indices have the same binary operation, and this implies that all facets have the same binary operation on $\Gamma$.
\end{proof}

Thus, there is a single binary operation $\circ$ on $\Gamma$ defined by $\circ=\circ_{123}$.

\begin{claim}
$(\Gamma,\circ)$ is a group.
\end{claim}
\begin{proof}
Clearly $\alpha\circ\ep=\ep\circ\alpha=\alpha$ for all $\alpha\in \Gamma$, since (A) and (B) hold for $t_{123}$, so $\ep$ is the identity. 
Also, the inverse of $\alpha$ is the unique element $\beta$ such that $\{\alpha_{12},\beta_{23},\ep_{13}\}\in\cB$, and this element exists by $(*)$. 
For all $\alpha,\beta,\gamma\in\Gamma$, by  \ref{useful} with 
$$(a,b,c,d,e,f)=(\beta_{23},\alpha_{12},(\alpha\circ \beta)_{13},((\alpha\circ\beta)\circ \gamma)_{14}, \gamma_{34}, (\beta\circ\gamma)_{24})$$
 we find that $\{\alpha_{12},(\beta\circ\gamma)_{24},((\alpha\circ\beta)\circ \gamma)_{14}\}\in\cB$.
 Thus, $\alpha\circ(\beta\circ\gamma)=(\alpha\circ\beta)\circ\gamma$, so $\circ$ is a group operation on $\Gamma$. 
 Note that we require $m\ge 4$ in order to prove associativity.
\end{proof}

Now, label each $e\in E(G)$ by $f(e)$, and assume that if $e\in L_{ij}$, then $e$ is oriented from $b_i$ to $b_j$. 
Then the biased graph obtained from this $\Gamma$-labeled graph is $(G,\cB)_{(m,\Gamma)}$; we will write $(G,\cB)_{(m,\Gamma)}=(G,\cB')$ for convenience.
It is easy to see that $\cB'$ satisfies $(*)$, by taking $c$ to be an element so that the product of the group labels is $\ep$.

We now complete the proof by showing that $\cB=\cB'$.
We first show that $\{a,b,c\}\in \cB$ if and only if $\{a,b,c\}\in \cB'$. Let $a\in L_{ij}$, $b\in L_{jk}$, and $c\in L_{ik}$. If $\{a,b,c\}\in\cB$, then $f(a)\circ f(b)=f(c)$, and due to the orientation of edges of $G$, the value of the cycle is $f(a)\circ f(b)\circ(f(c))^{-1}=\epsilon$, so $\{a,b,c\}\in \cB'$. 
If $\{a,b,c\}\in\cB'$, then $f(a)\circ f(b)\circ(f(c))^{-1}=\epsilon$, and so $f(a)\circ f(b)=f(c)$ and thus $\{a,b,c\}\in\cB$.
Therefore, $\cB$ and $\cB'$ contain the same 3-cycles of $G$.

We now show that $\cB\subseteq \cB'$. 
Let $C\in \cB$ have minimum size so that $C\notin \cB'$.
Then $|C|\ge 4$ since each cycle of $\cB$ has size at least three,  and $\cB$ and $\cB'$ contain the same 3-cycles.
Without loss of generality assume that $C=\{a_1,a_2,\dots,a_{k-1},a_k\}$ such that $a_j\in L_{j,j+1}$ for each $j\in [1,k-1]$, and $a_k\in L_{1k}$. 

By $(*)$ for $\cB$ there is some $d\in L_{13}$ so that $\{a_1,a_2,d\}\in \cB$.
Since $\cB$ satisfies the theta property, the cycle $\{d,a_3,a_4,\dots,a_k\}\in\cB$. 
This cycle has size less than $k$ and is thus in $\cB'$, by the minimality of $k$.
Then since $\{a_1,a_2,d\}\in\cB'$ and $\cB'$ satisfies the theta property, $C\in \cB'$.
The same argument can be applied to show that $\cB'\subseteq \cB$.
Thus, $(G,\cB)=(G,\cB)_{m,\Gamma}$.
\end{proof}


We have the following corollary for lift geometries.

\begin{theorem} \label{identify LG}
Let $M$ be a simple matroid with an element $x$ so that $\si(M/x)\cong M(K_n)$ for some $n\ge 4$, and each parallel class of $M/x$ has size $t\ge 2$.
If there is no element $e\in E(M)$ for which $M/e$ has a $U_{2,t+2}$-restriction, then there is a group $\Gamma$ so that $M\cong \LG^+(n,\Gamma)$.
\end{theorem}
\begin{proof}
Let $(G,\cB)$ be a biased graph on $n$ vertices so that the extended lift matroid of $(G,\cB)$ is $M$.
Then $\cB$ satisfies the following property:
\begin{enumerate}[$(*)$]
\item for all distinct $b_1,b_2,b_3\in V(G)$, each edge $a$ with ends $b_1$ and $b_2$, and each edge $b$ with ends $b_1$ and $b_3$, there is an edge $c$ with ends $b_2$ and $b_3$ so that $\{a,b,c\}\in\cB$.
\end{enumerate}
Otherwise, $M/a$ has a $U_{2,t+2}$-restriction consisting of $x,b$ and all edges of $G$ with ends $b_2$ and $b_3$. 
By Theorem \ref{biased graphs}, there is a group $\Gamma$ so that the extended lift matroid of $(G,\cB)$ is isomorphic to $\LG^+(n,\Gamma)$.
\end{proof}

For each integer $n\ge 3$, we say that a simple matroid $M$ is a rank-$n$ \emph{doubled-clique} if there is an element $e\in E(M)$ so that $\si(M/e)\cong M(K_n)$, and each parallel class of $M/e$ has size two.
Doubled cliques play the role for lift matroids that $B$-cliques play for frame matroids; they are the `complete' extended lift matroids.
We show that every doubled clique in $\cU(\ell)$ with sufficiently large rank has an $\LG^+(k,\Gamma)$-minor for some group $\Gamma$ of size at most $\ell$.

\begin{proposition} \label{LG minor}
For all integers $\ell\ge 2$ and $m\ge 3$, if $M\in\cU(\ell)$ is a rank-$m^{2^{\ell}}$ doubled clique, then there is a group $\Gamma$ so that $2\le |\Gamma|\le \ell$, and $M$ has an $\LG^+(m,\Gamma)$-minor.
\end{proposition}
\begin{proof}
Define a function $g\colon \{2,3,\dots,\ell\}\to \bZ$ by $g(\ell)=m$, and $g(s)=m{g(s+1)\choose 2}$ for $s<\ell$. 
It is an easy induction proof to show that $g(s)\le m^{(2^{\ell-s+1}-1)}$ for all $s\ge 2$.
Let $k=g(2)$, and note that $k\le m^{2^{\ell}}$. 
For all integers $d\ge 2$ and $r\ge 3$, a \emph{rank-$r$ $d$-doubled clique} is a simple matroid $M$ 
with an element $e$ so that $\si(M/e)\cong M(K_r)$, and each parallel class of $M/e$ has size $d$.
Note that a 2-doubled clique is a doubled clique. 

Let $s\in \{2,3,\dots,\ell\}$ be maximal so that there exists a rank-$g(s)$ $s$-doubled clique $M\in \cU(\ell)$ with no $\LG^+(m,\Gamma)$-minor for some nontrivial group $\Gamma$.
If $s=\ell$, then since $M$ has no $U_{2,\ell+2}$-minor, we have $M\cong \LG^+(m,\Gamma)$, by Theorem \ref{identify LG}.
Thus,  $s<\ell$.

Let $(M/e)|X$ be a simplification of $M/e$, and note that $(M/e)|X\cong M(K_{g(s)})$. 
Let $\cX$ be a collection of ${g(s+1)\choose 2}$ size-$m$ subsets of $X$ which correspond to pairwise vertex-disjoint cliques of $(M/e)|X$. 
For each $F\in\cX$, there are distinct $x_F,y_F\in F$ so that $\cl_{M/x_F}(\{e,y_F\})$ has a $U_{2,s+2}$-restriction; otherwise, the restriction of $M$ consisting of all lines through $t$ and an element of $F$ is $\LG^+(m,\Gamma)$, by Theorem \ref{identify LG}.

Let $C=\{x_F\colon F\in\cX\}$, and let $X'=\{y_F\colon F\in\cX\}$.
Then $\si((M/e/C)|X)$ is isomorphic to a clique, and $X'$ corresponds to a size-${g(s+1)\choose 2}$ matching of $(M/e/C)|X$. 
Thus, $(M/e/C)|X$ has a minor with ground set $X'$ which is isomorphic to $M(K_{g(s+1)})$.
So there is some $C'\subseteq X-X'$ so that $(M/e/C/C')|X'\cong M(K_{g(s+1)})$ and $e\notin\cl_M(C\cup C')$.

For each $y\in X'$, the matroid $(M/C)|\cl_{M/C}(\{e,y\})$ has a $U_{2,s+2}$-restriction, by the definitions  of $C$ and $X'$.
In $M/C/C'$, there is no line through $e$ which contains two elements of $X'$, since $(M/e/C/C')|X'$ is simple.
Thus, for each $y\in X'$, the set $\{e,y\}$ spans a distinct $U_{2,s+2}$-restriction of $M/C/C'$.
Since $(M/e/C/C')|X'\cong M(K_{g(s+1)})$, this implies that $(M/C/C')|(\{e\}\cup X')$ has an  $(s+1)$-doubled-clique restriction of rank $g(s+1)$. 
By the maximality of $s$, this matroid has an $\LG^+(m,\Gamma)$-minor with $|\Gamma|\ge 2$.  
\end{proof}
\end{subsection}
\end{section}


\begin{section}{Finding a Dowling Geometry} \label{frame chap}
In this section we prove some fundamental properties of Dowling geometries.


\begin{subsection}{Finding a Dowling-Geometry Minor} \label{Dowling props}
Recall that a \emph{B-clique} is a matroid $M$ framed by $B$ so that each pair of elements of $B$ is contained in a long line of $M$.
The following theorem of Kahn and Kung \cite{Kahn/Kung} characterizes when a simple $B$-clique is a Dowling geometry; Theorem \ref{identify LG} is the analogue of this result for lift geometries.
We prove it using Theorem \ref{biased graphs}, the proof of which is based on \cite{Kahn/Kung}.


\begin{theorem}[Kahn, Kung] \label{dowling}
Let $t\ge 1$ be an integer, and let $M$ be a simple $B$-clique of rank at least four so that $|\cl_M(\{b_1,b_2\})|=t+2$ for all distinct $b_1,b_2\in B$. 
If there is no element $e\in E(M)$ so that $M/e$ has a $U_{2,t+3}$-restriction,
then there exists a group $\Gamma$ with $|\Gamma|=t$ so that $M\cong \DG(r(M),\Gamma)$.
\end{theorem}
\begin{proof}
Let $(G,\cB)$ be a biased graph with vertex set $B$ so that the frame matroid of $(G,\cB)$ is $M$; note that each element of $B$ is an unbalanced loop of $(G,\cB)$.
Then 
\begin{enumerate}[$(*)$]
\item for all distinct $b_1,b_2,b_3\in B$, each edge $a$ with ends $b_1$ and $b_2$, and each edge $b$ with ends $b_1$ and $b_3$, there is an edge $c$ with ends $b_2$ and $b_3$ so that $\{a,b,c\}\in\cB$,
\end{enumerate}
or else $M/a$ has a $U_{2,t+3}$-restriction consisting of $b,b_2,b_3$ and all edges of $G$ with ends $b_2$ and $b_3$. 
By Theorem \ref{biased graphs}, there is a group $\Gamma$ so that the frame matroid of $(G,\cB)$ is isomorphic to $\DG(r(M),\Gamma)$.
\end{proof}


We will need the following straightforward lemma about $B$-cliques, which is essentially equivalent to the statement that a complete graph with a matching $J$ of size ${k\choose 2}$ has a $K_k$-minor with edge-set $J$.

\begin{lemma} \label{skew lines}
Let $k,t\ge 3$ be integers,  let $M$ be a $B$-clique, and let $\cX$ be a collection of ${k\choose 2}$ pairwise disjoint 2-element subsets of $B$.
Then there is some $B'\subseteq B$ with $|B'|=k$, and a $B'$-clique minor $N$ of $M$ such that $r_N(X)=r_M(X)=2$ for all $X\in\cX$, and for all distinct $b_1,b_2\in B'$ there is some $X\in \cX$ so that $\cl_M(X)\subseteq \cl_N(\{b_1,b_2\})$.
\end{lemma}
\begin{proof}
We write $\cup\cX$ for $\cup_{X\in \cX}X$.
We may assume that $B=\cup \cX$.
Let $(B_1,B_2,\dots, B_k)$ be a partition of $\cup \cX$ so that $|B_i|=k-1$ for each $i\in [k]$, and for each $1\le i<j\le k$ there is some $X\in \cX$ such that $|X\cap B_i|=1$ and $|X\cap B_j|=1$. 
This amounts to arbitrarily putting $\cX$ into bijection with $E(K_k)$.
For each $i\in [k]$, fix some $b_i\in B_i$ and let $C_i\subseteq \cl_M(B_i)-B$ be a set of size $k-2$ so that for each $b\in B_i-\{b_i\}$ there is some $e\in (C_i\cap \cl_M(\{b_i,b\}))-B$.
Such a set $C_i$ exists since $M$ is a $B$-clique.
Then $B_i\cap\cl_M(C_i)=\varnothing$ or else $C_i$ spans $B_i$.
Since $|B_i|=|C_i|+1$ we have $r_{M/C_i}(B_i)=1$.

Let $N=M/(\cup_i C_i)$, and let $B'$ be a transversal of $(B_1,B_2,\dots,B_k)$ which is independent in $N$.
Then $B'$ spans $N$ since $B=\cup \cX$, so $N$ is a $B'$-clique.
Let $b_1,b_2\in B'$ with $b_1\ne b_2$, and assume without loss of generality that $b_1\in B_1$ and $b_2\in B_2$. 
There is some $X\in\cX$ such that $|X\cap B_1|=1$ and $|X\cap B_2|=1$. 
Let $b_1'\in X\cap B_1$ and $b_2'\in X\cap B_2$, and note that $b_j$ and $b_j'$ are parallel in $N$ for each $j\in \{1,2\}$ since $B_i\cap \cl_M(C_i)=\varnothing$ and $r_{M/C_i}(B_i)=1$ for each $i\in [k]$.
Thus, $r_N(X)=r_M(X)=2$, and $\cl_M(X)\subseteq \cl_N(X)=\cl_N(\{b_1,b_2\})$.
\end{proof}


We now combine Theorem \ref{dowling} and Lemma \ref{skew lines} to show that every $B$-clique in $\cU(\ell)$ with very large rank contains a big Dowling-geometry minor with group size as large as possible. 
The idea is that we contract elements and increase density until we have a $B$-clique minor such that each pair of elements of $B$ spans a $U_{2,\ell+1}$-restriction, and this must be a Dowling geometry by Theorem \ref{dowling}.
This easily implies a result for $B$-cliques which is analogous to Proposition \ref{LG minor}.

\begin{lemma} \label{DGminor}
Let $\ell-1\ge t\ge 1$ and $k\ge 3$ be integers.
If $M\in \cU(\ell)$ is a $B$-clique and $\cX$ is a collection of pairwise disjoint 2-element subsets of $B$ such that $|\cX|\ge k^{2^{\ell}}$, and each $X\in\cX$ satisfies $\elem(M|\cl_M(X))\ge t+2$, then there is some $B'\subseteq B$ and a $B'$-clique minor of $M$ which is isomorphic to $\DG(k,\Gamma)$ with $|\Gamma|\ge t$. 
\end{lemma}
\begin{proof}
Define a function $f_1\colon \bZ^3\to \bZ$ by $f_{1}(\ell,\ell,k)=1$, and $f_{1}(\ell,t,k)={(k+1)\cdot f_1(\ell,t+1,k)\choose 2}$ for $1\le t<\ell$. 
Define $f_{\ref{DGminor}}(\ell,k)=f_1(\ell,1,k)$, and note that $f_{\ref{DGminor}}(\ell,k)\ge f_1(\ell,t,k)$ for all $1\le t<\ell$.
One can show using induction that $f_{\ref{DGminor}}(\ell,k)\le k^{2^{\ell}-2}$. 
Fix $\ell\ge 2$, and let $t$ be maximal so that there exists a $B$-clique $M\in \cU(\ell)$ with a collection $\cX$ of pairwise disjoint 2-element subsets of $B$ such that $|\cX|\ge f_{\ref{DGminor}}(\ell,k)$ and each $X\in\cX$ satisfies $\elem(M|\cl_M(X))\ge t+2$, but there is no $B'\subseteq B$ and a $B'$-clique minor of $M$ which is isomorphic to $\DG(k,\Gamma)$ with $|\Gamma|\ge t$.

By Lemma \ref{skew lines}, there is some $B'\subseteq B$ with $|B'|=(k+1)\cdot f_{1}(\ell,t+1,k)$ and a $B'$-clique minor $N$ of $M$ such that $\elem(N|\cl_N(\{b_1,b_2\}))\ge t+2$ for all distinct $b_1,b_2\in B'$.
If $t=\ell-1$ then $N$ is a rank-$(k+1)$ Dowling geometry by Theorem \ref{dowling}, so $M$ is not a counterexample. Thus, $t<\ell-1$. 
Let $(B_1,B_2,\dots,B_h)$ be a partition of $B'$ such that $h=f_{1}(\ell,t+1,k)$ and $|B_i|=k+1$ for each $i\in [h]$. Since $M$ is a counterexample and $r_M(B_i) \ge k+1\ge 4$, for each $i\in [h]$ there is some $e_i\in \cl_N(B_i)$ such that $N/e_i$ has a $U_{2,t+3}$-restriction spanned by two elements of $B_i$, by Theorem \ref{dowling}. 
By the maximality of $t$, there is some $B''\subseteq B'$ and a $B''$-clique minor of $N/(\cup_i e_i)$ which is isomorphic to $\DG(k,\Gamma)$ with $|\Gamma|\ge t+1$, which contradicts that $M$ is a counterexample. 
\end{proof}

As an easy corollary, we can bound the size of any simple $B$-clique in $\cU(\ell)$ with no Dowling-geometry minor over a large group.  

\begin{corollary} \label{mader}
For all integers $\ell\ge t\ge 2$ and $k\ge 3$, if $M\in \cU(\ell)$ is a $B$-clique with no $\DG(k,\Gamma)$-minor with $|\Gamma|\ge t$, then $\elem(M)\le (t-1){r(M)\choose 2}+2\ell k^{2^{\ell}}\cdot r(M)$.
\end{corollary}
\begin{proof}
Let $m=k^{2^{\ell}}$.
Let $M\in\cU(\ell)$ be a $B$-clique, and let $\cX$ be a maximum-size collection of pairwise disjoint 2-element subsets of $B$ so that each $X\in\cX$ satisfies $\elem(M|\cl_M(X))\ge t+2$.
By Lemma \ref{DGminor} we have $|\cX|<m$, so by the maximality of $|\cX|$ there are at most $2(m-2)\cdot r(M)$ pairs $b,b'\in B$ such that $\elem(M|\cl_M(\{b,b'\}))\ge t+2$. 
Since $M$ is a $B$-clique and has no $U_{2,\ell+2}$-restriction, each pair $b,b'\in B$ satisfies $\elem(M|(\cl_M(\{b,b'\})-\{b,b'\}))\le \ell-1$.
Thus, $\elem(M)\le (t-1){r(M)\choose 2}+(2(m-2)\ell+1)r(M)$. 
\end{proof}

\end{subsection}


\begin{subsection}{Finding a Dowling-Geometry Restriction} \label{Dowling rest}
Recall that $\DG^-(k,\Gamma)= \DG(k,\Gamma)\del B$, where $B$ is a frame of $\DG(k,\Gamma)$. 
We will show that each matroid in $\cU(\ell)$ with no $\LG^+(n,\Gamma)$-minor with $|\Gamma|\ge 2$, and with a $\DG(r,\Gamma)$-minor, has a $\DG^-(k,\Gamma)$-restriction if $r$ is sufficiently large.
The key is the following lemma, which provides sufficient conditions for a $\DG^-(k,\Gamma)$-minor of a matroid $M$ to be a restriction of $M$.

Before stating the lemma we define a special type of circuit of $\DG^-(k,\Gamma)$.
Let $G$ be a Dowling geometry over a group $\Gamma$, and let $B$ be a frame for $G$.
We say that a circuit $C$ of $G\del B$ is \emph{balanced} if $\cl_G(C)\cap B=\varnothing$.
The balanced circuits are precisely the circuits of $G$ which correspond to balanced cycles of a $\Gamma$-labeled graph associated with $G$.


\begin{lemma} \label{jointless rest}
Let $k\ge 4$ be an integer, let $\Gamma$ be a finite group, and let $M$ be a matroid with a $\DG^-(k,\Gamma)$-minor $G$.
If each balanced circuit of $G$ of size at most four is also a circuit of $M$, then $G$ is a restriction of $M$.
\end{lemma}
\begin{proof}
Let $G_1$ be a Dowling geometry with frame $B=\{b_1,b_2,\dots,b_k\}$ so that $G_1\del B=G$.
Assume for a contradiction that $r(G)<r_M(E(G))$.
Let $H$ be a $\Gamma$-labeled graph so that $V(H)=B$ and the associated frame matroid is $G$, and each edge between vertices $b_i$ and $b_j$ with $i<j$ is oriented from $b_i$ to $b_j$.
For each $\alpha\in\Gamma$ and $1\le i<j\le k$, let $(i,j)_{\alpha}$ denote the element of $G$ spanned by $\{b_i,b_j\}$ and labeled by $\alpha$ in the $\Gamma$-labeling of $H$, and let $E_{\alpha}$ denote the set of elements of $G$ labeled by $\alpha$.
Let $\ep$ denote the identity element of $\Gamma$.

\begin{claim} \label{identity edges}
$M|E_{\ep}=G|E_{\ep}$.
\end{claim}
\begin{proof}
If not, then let $C'$ be a circuit of $G|E_{\ep}$ of minimum size so that $C'$ is independent in $M$.
Since each circuit of $G|E_{\ep}$ is a balanced circuit of $G$ we have $|C'|>4$.
Let $C_1$ and $C_2$ be distinct circuits of $G|E_{\ep}$ of size less than $|C'|$ so that $C_1,C_2,C'$ are the cycles of a theta subgraph of $H$; these circuits may be obtained by adding a chord to $C'$.
Then $C'\subseteq C_1\cup C_2$ and $r_{G}(C_1\cup C_2)= |C_1\cup C_2|-2$.
By the minimality of $|C'|$, both $C_1$ and $C_2$ are circuits of $M$.
Then $M|(C_1\cup C_2)$ has at least two circuits, so $r_M(C_1\cup C_2)\le |C_1\cup C_2|-2\le r_{G}(C_1\cup C_2)$. 
Since $r_G(C_1 \cup C_2) \ge r_M(C_1 \cup C_2)$ and $G$ is a minor of $M$, it follows that $M|(C_1\cup C_2)=G|(C_1\cup C_2)$.
This implies that $M|C'=G|C'$, a contradiction.
\end{proof}

If $|\Gamma|=1$, then \ref{identity edges} shows that $M|E(G)=G$, so we may assume that $|\Gamma|\ge 2$.
This implies that $r_M(E_{\ep})<r(G)<r_M(E(G))$.


\begin{claim} \label{same column}
Each $\gamma\in \Gamma-\{\ep\}$ satisfies 
$M|(E_{\ep}\cup E_{\gamma})=G|(E_{\ep}\cup E_{\gamma})$, and $E_{\gamma}\subseteq \cl_M(E_{\ep}\cup \{(1,2)_{\gamma}\})$. 
\end{claim} 
\begin{proof}
We show that $(i,j)_{\gamma}\in\cl_M(E_{\ep}\cup\{(1,2)_{\gamma}\})$ for all $1\le i<j\le k$, which implies that $$r_M(E_{\ep}\cup E_{\gamma})\le r_M(E_{\ep})+1=r_{G}(E_{\ep})+1\le r_{G}(E_{\ep}\cup E_{\gamma}),$$
where the first equality holds by \ref{identity edges}.
We first prove this in the case that $2<i<j\le k$.
The set $\{(1,2)_{\gamma},(2,j)_{\ep},(i,j)_{\gamma},(1,i)_{\ep}\}$ is a balanced circuit of $G$ of size at most four and is thus a circuit of $M$, which implies that $(i,j)_{\gamma}\in\cl_M(E_{\ep}\cup \{(1,2)_{\gamma}\})$.
Similarly, the balanced circuit $\{(1,2)_{\gamma},(2,j)_{\ep},(1,j)_{\gamma}\}$ shows that $(1,j)_{\gamma}\in\cl_M(E_{\ep}\cup \{(1,2)_{\gamma}\})$ for all $3\le j\le k$.
Using this, the balanced circuit $\{(1,2)_{\ep},(2,j)_{\gamma},(1,j)_{\gamma}\}$ shows that $(2,j)_{\gamma}\in\cl_M(E_{\ep}\cup \{(1,2)_{\gamma}\})$ for all $3\le j\le k$.
\end{proof}


Fix some $\gamma\in \Gamma-\{\ep\}$.
We will show that $E(G)\subseteq \cl_M(E_{\ep}\cup E_{\gamma})$; then 
$$r_M(E(G))\le r_M(E_{\ep}\cup E_{\gamma})= r_{G}(E_{\ep}\cup E_{\gamma})\le r(G),$$
where the first equality holds by \ref{same column}.
By \ref{same column} it suffices to show that each $\beta\in \Gamma-\{\ep,\gamma\}$ satisfies $(1,2)_{\beta}\in \cl_M(E_{\ep}\cup E_{\gamma})$.
Let $\beta\in \Gamma-\{\ep,\gamma\}$, and 
note that $C_1=\{(1,2)_{\beta},(2,4)_{\ep},(3,4)_{\beta},(1,3)_{\ep}\}$ and $C_2=\{(1,2)_{\beta},(2,4)_{\beta},(3,4)_{\beta},(1,3)_{\beta}\}$ are both balanced circuits of $G$.
Also, the set $\{(1,3)_{\beta},(1,3)_{\ep},(2,4)_{\ep},(2,4)_{\beta}\}$ is independent in $G$ since it contains no balanced cycle of $H$, and no handcuff or theta subgraph of $H$.
Thus, $$r_G(C_1\cup C_2)\ge 4=|C_1\cup C_2|-2\ge r_M(C_1\cup C_2),$$ 
where the last inequality holds because $M|(C_1\cup C_2)$ contains two distinct circuits.
Since $r_G(C_1 \cup C_2) \ge r_M(C_1 \cup C_2)$ and $G$ is a minor of $M$, it follows that $M|(C_1\cup C_2)=G|(C_1\cup C_2)$.
The set $C=\{(1,3)_{\beta},(1,3)_{\ep},(1,2)_{\beta},(2,4)_{\beta},(2,4)_{\ep}\}$ is a loose handcuff of $H$ which contains no balanced cycle of $H$, and is thus a circuit of $G$.
Since $C\subseteq C_1\cup C_2$ and $M|(C_1\cup C_2)=G|(C_1\cup C_2)$ it follows that $C$ is a circuit of $M$.
Then $(1,2)_{\beta}\in \cl_M(C-\{(1,2)_{\beta}\})$, since $M|C$ is a circuit.
Since $(C-\{(1,2)_{\beta}\})\subseteq E_{\ep}\cup E_{\gamma}$, this implies that $(1,2)_{\beta}\in \cl_M(E_{\ep}\cup E_{\gamma})$, as desired.
\end{proof}


We also need to be able to recognize when a matroid has an $\LG^+(n,\Gamma)$-minor for some nontrivial group $\Gamma$.

\begin{lemma} \label{find d clique}
Let $\ell\ge 2$ and $n\ge 3$ be integers, and let $M\in \cU(\ell)$ be simple with an element $e$ so that $M/e$ has a spanning $B$-clique restriction. 
If there is a collection $\cS$ of $n^{2^{2\ell+1}}$ pairwise disjoint 2-subsets of $B$ so that each $S\in \cS$ spans a nontrivial parallel class of $M/e$ which is disjoint from $S$, then $M$ has an $\LG^+(n,\Gamma)$-minor for some nontrivial group $\Gamma$ so that $e_0=e$.
\end{lemma}
\begin{proof}
Let $m=n^{2^{\ell}}$, and note that $n^{2^{2\ell+1}}=m^{2^{\ell+1}}$.
We will show that $M$ has a rank-$m$ doubled-clique minor, and then apply Proposition \ref{LG minor}.
Note that $m^{2^{\ell+1}}\ge {m^{2^{\ell}}\choose 2}$, and assume that $|\cS|={m^{2^{\ell}}\choose 2}$.
Let $X$ and $Y$ be disjoint transversals of the parallel classes of $M/e$ spanned by sets in $\cS$, so $|X|=|Y|={m^{2^{\ell}}\choose 2}$. 

\begin{claim} \label{find C}
There is some $C\subseteq E(M/e)$ and $X'\subseteq X$ so that $(M/e/C)|X'\cong M(K_m)$ and $e\notin\cl_{M}(C)$.
\end{claim}
\begin{proof}
By Lemma \ref{skew lines} with $\cX=\cS$ and $k=m^{2^{\ell}}$, there is some $B_1\subseteq B$ and a $B_1$-clique minor $M_1$ of $M/e$ so that $|B_1|=m^{2^{\ell}}$ and for all $b,b'\in B_1$ with $b\ne b'$ there is some $x\in X$ for which $\{b,b',x\}$ is a circuit of $M_1$.  
Then $M_1|(B_1\cup X)$ is a $B_1$-clique, so by Lemma \ref{DGminor} with $k=m$ there is some $B_2\subseteq B_1$ and a $B_2$-clique minor $M_2$ of $M_1|(B_1\cup X)$ so that $M_2\cong M(K_{m+1})$.

Then $B_2$ corresponds to a spanning star of $K_{m+1}$, so $M_2\del B_2\cong M(K_m)$.
Since $M_1|(B_1\cup X)$ is a $B_1$-clique, each element of $B_1-B_2$ is parallel in $M_2$ to an element of $B_2$.
Since $M_2$ is simple and $E(M_2)\subseteq B_1\cup X$, this implies that $E(M_2\del B_2)\subseteq X$.
Thus, the claim holds by taking $X'=E(M_2\del B_2)$, and $C\subseteq E(M)$ to be a set so that $M_2$ is a restriction of $M/e/C$.
\end{proof}

Let $Y'\subseteq Y$ be the set of elements on a line of $M$ through $e$ and an element of $X'$. 
In $M/C$ there is no line through $e$ which contains two elements of $X'$, since $(M/e/C)|X'$ is simple.
Therefore, $(M/C)|(X'\cup Y'\cup \{e\})$ is simple, and is thus a doubled clique by  \ref{find C}.
By Proposition \ref{LG minor} with $m=n$, this implies that $M/C$ has an $\LG^+(n,\Gamma)$-minor with $|\Gamma|\ge 2$.
\end{proof}


We also need to be able to recover the frame of a $\DG^-(k,\Gamma)$-restriction so that we can apply Lemma \ref{find d clique}.
The following lemma shows that we can recover $\DG(k,\Gamma)$ from $\DG^-(k+2,\Gamma)$ by contracting two elements.

\begin{lemma} \label{jointless}
Let $G_1$ be a matroid so that $G_1\cong \DG(r(G_1),\Gamma)$, and let $B$ be a frame of $G_1$.
If $G=G_1\del B$, then each pair of elements of $B$ spans a subset $C$ of $E(G)$ so that $|C|\le 2$, and $G/C$ has a $\DG(r(G/C),\Gamma)$-restriction. 
\end{lemma}
\begin{proof}
If $|\Gamma|=1$, then $G\cong \DG(r(G),\Gamma)$ and the lemma holds with $C=\varnothing$, so we may assume that $|\Gamma|\ge 2$.
Let $b_1$ and $b_2$ be distinct elements of $B$, and let $e_1$ and $e_2$ be distinct elements in $\cl_{G_1}(\{b_1,b_2\})\cap E(G)$.
We will show that each $b\in B-\{b_1,b_2\}$ is parallel to an element of $G$ in $G_1/\{e_1,e_2\}$.
Let $b\in B-\{b_1,b_2\}$, and let $x\in \cl_{G_1}(\{b,b_1\})\cap E(G)$.
Then $r_{G_1}(\{e_1,e_2,b,x\})\le 3$ since this set is spanned by $\{b_1,b_2,b\}$ in $G_1$.
Also, $\{b,x\}$ is disjoint from $\cl_{G_1}(\{e_1,e_2\})$, since $\{e_1,e_2\}$ and $\{b_1,b_2\}$ span the same flat of $G_1$ and $\{b,x\}$ is disjoint from $\cl_{G_1}(\{b_1,b_2\})$.
Thus, $\{b,x\}$ is a parallel pair of $G_1/\{e_1,e_2\}$, and so $\si(G/\{e_1,e_2\})$ is isomorphic to $\si(G_1/\{e_1,e_2\})$.
Since $\si(G_1/\{e_1,e_2\})\cong \DG(r(G)-2,\Gamma)$, this implies that $G/\{e_1,e_2\}$ has a $\DG(r(G)-2,\Gamma)$-restriction.
\end{proof}


We now prove a result which lets us move from a $\DG^-(r,\Gamma)$-minor of a matroid $M$ to a $\DG(m,\Gamma)^-$-restriction whenever $M\in\cU(\ell)$ has no $\LG^+(n,\Gamma)$-minor with $|\Gamma|\ge 2$, and $r$ is sufficiently large.

\begin{lemma} \label{Ramsey}
There is a function $f_{\ref{Ramsey}} \colon\bZ^3\to\bZ$ so that for all integers $\ell\ge 2$,  $m,n\ge 3$, and $d\ge 0$, and each finite group $\Gamma$, 
 if $M\in\cU(\ell)$ has no $\LG^+(n,\Gamma')$-minor with $|\Gamma'|\ge 2$, and has a $\DG^-(m+f_{\ref{Ramsey}}(\ell,n,d),\Gamma)$-minor $G$ for which $r_M(E(G))-r(G)\le d$, then $M|E(G)$ has a $\DG^-(m,\Gamma)$-restriction.
\end{lemma}
\begin{proof}
Let $m_1=n^{2^{2 \ell+1}}$, and define $f_{\ref{Ramsey}}(\ell,n,d)=8(d+2)m_1+2$.
Assume for a contradiction that the lemma is false.
Let $M$ be a counterexample so that $M/C_0$ has a $\DG^-(m+8(d+2)m_1+2,\Gamma)$-restriction $G_0$, and $|C_0|$ is minimal over all counterexamples.
If there is some $e\in C_0-\cl_M(E(G_0))$, then the lemma holds for $M/e$ if and only if it holds for $M$, so $C_0\subseteq \cl_M(E(G_0))$.
This implies that $r_M(C_0)=r_M(E(G_0))-r(G_0)\le d$.
The minimality of $|C_0|$ also implies that $C_0$ is independent in $M$.

By Lemma \ref{jointless} there is a set $C_1\subseteq E(G_0)$ so that $|C_1|\le 2$ and $G_0/C_1$ has a $\DG(m+8(d+2)m_1,\Gamma)$-restriction $G$.
Let $C=C_0\cup C_1$, and note that $r_M(C)\le d+2$ and that $C$ is independent in $M$.
Let $B$ be a frame for $G$, and let $H$ be a $\Gamma$-labeled graph so that $V(H)=B$ and the associated frame matroid is $G\del B$.
We say that two circuits of $G$ are \emph{vertex-disjoint} if the corresponding cycles of $H$ are vertex-disjoint.
Let $\cC$ be a maximum-size collection of pairwise vertex-disjoint balanced circuits of $G$, each of size at most four, which are independent in $M$.
If $|\cC|$ is large, then we can find an $\LG^+(n,\Gamma)$-minor for some nontrivial group $\Gamma$.


\begin{claim}
$|\cC|\le (d+2)(m_1-1)$.
\end{claim}
\begin{proof}
Assume for a contradiction that $|\cC|> (d+2)(m_1-1)$.
Since each set $Y\in\cC$ is independent in $M$, each $Y$ satisfies $\sqcap_{M}(Y,C)=1$.
Since $C$ is independent in $M$, this implies that $(M/Y)|C$ contains a unique circuit.
For each $Y \in \cC$, let $c_Y$ be an element in the unique circuit of $(M/Y)|C$.
Then $\sqcap_{M}(Y,C-\{c_Y\})=0$ since $C - \{c_Y\}$ is independent in $M/Y$.
We now contract all but one element of $C$.
Since $|\cC|> (d+2)(m_1-1)$ and $|C|\le d+2$,  there is some $c\in C$ and some $\cC_1\subseteq \cC$ so that $|\cC_1|=m_1$ and each $Y\in \cC_1$ satisfies $\sqcap_{M}(Y,C-\{c\})=0$.
Let $M_1=(M/(C-\{c\}))|(E(G)\cup \{c\})$.
Then $M_1/c\cong G$, and each $Y\in \cC_1$ is independent in $M_1$.

For each $Y\in \cC_1$, let $Y'\subseteq Y$ so that $|Y'|=|Y|-2$.
Let $Z=\cup_{Y\in \cC_1} Y'$, and let $M_2$ be a simplification of $M_1/Z$. 
Let $B_1\subseteq B$ be a frame for $G/Z$.
Then $\cl_{G}(B_1)$ is a spanning Dowling-geometry restriction of $M_2/c$, and for each $Y\in\cC_1$ the set $(Y-Y')\cup \{c\}$ is a circuit of $M_2$, since $Y-Y'$ is a parallel pair of $G/Z$ but not of $M_2$.
Also, for each $Y\in\cC_1$ the parallel class of $Y-Y'$ in $G/Z$ is disjoint from $B_1$, since $Y$ is a balanced circuit of $G$.
So for each $Y\in \cC_1$ there is some $B_Y\subseteq B_1$ so that $|B_Y|=2$ and $Y-Y'\subseteq \cl_{M_2/c}(B_Y)$. 
Note that the sets $B_Y$ are pairwise disjoint subsets of $B_1$, since the circuits in $\cC_1$ are pairwise vertex-disjoint.
Thus, $M_2$ is a simple matroid with an element $c$ so that $M_2/c$ has a spanning $B_1$-clique restriction, and there is a collection $\cX$ of $m_1$ pairwise disjoint 2-subsets of $B_1$ for which each $X\in \cX$ spans a nontrivial parallel class of $M_2/c$ which contains neither element of $X$.
Then $M_2$ has a rank-$n^{2^{\ell}}$ doubled-clique minor, by Lemma \ref{find d clique} with $n=n^{2^{\ell}}$.
But then $M_2$ has an $\LG^+(n,\Gamma)$-minor with $|\Gamma|\ge 2$, by Proposition \ref{LG minor} with $m=n$, a contradiction.
\end{proof}

Let $B_1\subseteq B$ be a set of minimum size so that $(\cup_{Y\in\cC}Y)\subseteq \cl_G(B_1)$, and let $B_2=B-B_1$.
Then $|B_1|\le 8|\cC|\le 8(d+2)(m_1-1)$, so $|B_2|\ge m+1\ge 4$.
By the maximality of $|\cC|$, each balanced circuit of $\cl_G(B_2)$ of size at most four is a circuit of $M$.
Then Lemma \ref{jointless rest} shows that $G|(\cl_G(B_2)-B_2)$ is a restriction of $M$.
\end{proof}

\end{subsection}

\end{section}


\begin{section}{Exploiting Density} \label{crit chap}
In this section, we exploit the density of a matroid $M$ to find restrictions which are incompatible with either Dowling geometries or lift geometries.
Given a collection $\cY$ of sets, we will write $\cup\cY$ for $\cup_{Y\in\cY}Y$, for convenience.


\begin{subsection}{An Upgraded Connectivity Reduction} \label{upgrade reduction}
We first prove a strengthening of a theorem from \cite{DenserthanClique}, which finds extremal matroids which either have a spanning clique restriction, or large vertical connectivity.
The proof closely follows \cite{DenserthanClique}, but we separate out two lemmas which we will also use in Section \ref{exact chap}.

The first lemma essentially follows from the Growth Rate Theorem, and the fact that quadratic functions with positive leading coefficient are concave up.


\begin{lemma} \label{define nu}
There is a function $\nu_{\ref{define nu}}\colon \bR^7\to \bZ$ so that for all integers $\ell,k\ge 2$ and $r,s\ge 1$ and any real quadratic polynomial $p(x)=ax^2+bx+c$ with $a>0$, if $M\in\cU(\ell)$ has no rank-$k$ projective-geometry minor, $r(M)>0$, and $\elem(M)>p(r(M))+\nu_{\ref{define nu}}(a,b,c,\ell,k,r,s)\cdot r(M)$, then $M$ has a vertically $s$-connected minor $N$ such that 
\begin{itemize}
\item $r(N)\ge r$ and $\elem(N)>p(r(N))+\nu_{\ref{define nu}}(a,b,c,\ell,k,r,s)\cdot r(N)$, and 

\item $\elem(N)-\elem(N/e)>p(r(N))-p(r(N)-1)+\nu_{\ref{define nu}}(a,b,c,\ell,k,r,s)$ for each $e\in E(N)$. 
\end{itemize}
\end{lemma}
\begin{proof}
Fix $\ell,k,r,s$, and $p$.
Let $n_0$ be a positive integer such that $p(x)> p(x-1)\ge 0$ for all real $x\ge n_0$.
By Theorem \ref{GRT} there is a real number $\alpha>0$ such that $\elem(M)\le \alpha p(r(M))$ for all matroids $M\in\cU(\ell)$ with no rank-$k$ projective-geometry minor and $r(M)\ge n_0$.
Let $n_1\ge \max(r,s,n_0)$ be an integer so that  
\begin{align*}
a(\alpha+2s)(x+y)+((\alpha+1)b+\alpha|c|)s+c-as^2\le 2axy
\end{align*}
for all real $x,y\ge n_1$.
Finally, define $\nu_{\ref{define nu}}(a,b,c,\ell,k,r,s)=\nu=\max(-b,\ell^{n_1},\ell^{n_1}-\min_{x\in \bR} p(x))$.
Note that the polynomial $p(x)+\nu x$ satisfies $p(x)+\nu  x\ge \ell^{n_1}$ for all $x\in \bR$, and is nondecreasing for all $x>0$.

Let $M\in\cU(\ell)$ with no rank-$k$ projective-geometry minor, $r(M)>0$ and $\elem(M)>p(r(M))+\nu r(M)$.
Let $N$ be a minimal minor of $M$ such that $r(N)>0$ and $\elem(N)>p(r(N))+\nu r(N)$.
Note that $N$ is simple.
Since $r(N)>0$ we have $\elem(N)\ge \ell^{n_1}$.
This implies that $r(N)\ge n_1$, since $N\in\cU(\ell)$. 
By the minor-minimality of $N$, each $e\in E(N)$ satisfies $\elem(N)-\elem(N/e)>p(r(N))-p(r(N)-1)+\nu$.
Since $r(N)\ge n_1\ge \max(r,s)$, $N$ is not vertically $s$-connected or else the lemma holds. 

Let $(A,B)$ be a partition of $E(N)$ so that $r_N(A)\le r_N(B)<r(N)$ and $r_N(A)+r_N(B)<r(N)+s-1$.
Let $r_A=r_N(A)$ and $r_B=r_N(B)$ and $r_N=r(N)$.
If $r_A<n_1$, then $|A|<\ell^{n_1}$, so 
\setcounter{equation}{0}
\begin{align}
|B|&=|N|-|A|\\
&>p(r_N)+\nu r_N-\ell^{n_1}\\
&\ge p(r_N-1)+\nu (r_N-1)\\
&\ge p(r_B)+\nu r_B,
\end{align}
which contradicts the minor-minimality of $N$.
Line (3) holds because $\nu\ge \ell^{n_1}$ and $p(x)\ge p(x-1)$ for all $x\ge n_1$, and
line (4) holds because $a>0$ and $\nu+b\ge 0$.
Thus, $r_B\ge r_A\ge n_1$.

We now show that $x=r_A$ and $y=r_B$ contradicts the definition of $n_1$.
Since $r_N\ge n_1\ge n_0$ we have $p(r_N)\ge 0$, so $\nu r_N<|N|\le \alpha p(r_N)$. 
Since $r_N\ge 1$, this implies that 
$$\nu\le \alpha\Big(ar_N+b+\frac{c}{r_N}\Big)\le \alpha\Big(a(r_A+r_B)+b+|c|)\Big).$$
Using the partition $(A,B)$, we have
\begin{align*}
p(r_A+r_B-s)+\nu(r_A+r_B-s)&\le p(r_N)+\nu r_N<|A|+|B|\\
&\le p(r_A)+\nu r_A+p(r_B)+\nu r_B,
\end{align*}
where the first inequality holds because $r_A+r_B-s\le r_N$ and $a>0$ and $\nu+b\ge 0$, and the third inequality holds by the minor-minimality of $N$.
Then expanding $p$ and simplifying gives
$$s(\nu+b)+c-as^2+2as(r_A+r_B)>2r_Ar_B.$$
Combining this with our upper bound for $\nu$ gives 
$$a(\alpha+2s)(r_A+r_B)+((\alpha+1)b+\alpha|c|)s+c-as^2>2r_Ar_B,$$
which contradicts that $r_B\ge r_A\ge n_1$.
\end{proof}


The second lemma is mostly a property of Dowling geometries, and relies on the fact that the closure of any subset of the frame of a Dowling geometry is itself a Dowling geometry.
The proof is taken almost verbatim from Claim 6.1.1 in \cite{DenserthanClique}.

\begin{lemma} \label{spanning Dowling}
There is a function $n_{\ref{spanning Dowling}}\colon \bR^6\to \bZ$ so that for all integers $\ell\ge 2$ and $r,s\ge 1$ and any real quadratic polynomial $q(x)=ax^2+bx+c$ with $a>0$, if $M\in \cU(\ell)$ satisfies $\elem(M)>q(r(M))$ and has a $\DG(n_{\ref{spanning Dowling}}(a,b,c,\ell,r,s),\Gamma)$-minor, then $M$ has a minor $N$ such that $r(N)\ge r$ and $\elem(N)>q(r(N))$, and either
\begin{enumerate}[$(a)$]
\item $N$ has a $\DG(r(N),\Gamma)$-restriction, or
\item $N$ has an $s$-element independent set $S$ so that each $e\in S$ satisfies $\elem(N)-\elem(N/e)>q(r(N))-q(r(N)-1)$.
\end{enumerate}
\end{lemma}
\begin{proof}
Define $n_{\ref{spanning Dowling}}(a,b,c,\ell,r,s)=(s(s-1)+1)n_2$, where $n_2\ge r+1$ is an integer such that $q(x)-q(y)\ge \ell^s$ for all real $x,y$ with $x\ge n_2$ and $x-1\ge y\ge 0$.

Let $M\in\cU(\ell)$ satisfy $\elem(M)>q(r(M))$ and have a $\DG(n_{\ref{spanning Dowling}}(a,b,c,\ell,r,s),\Gamma)$-minor $N_1$. 
Let $M_1$ be a minimal minor of $M$ so that $\elem(M_1)>q(r(M_1))$ and $N_1$ is a minor of $M_1$, and let $C$ be an independent set in $M_1$ so that $N_1$ is a spanning restriction of $M_1/C$.
We may assume that $|C|<s$ or else $M_1$ and $C$ satisfy (b), by the minimality of $M_1$.

Let $i\ge 0$ be minimal so that there is a minor $M_2$ of $M_1$ for which $\elem(M_2)>q(r(M_2))$, and there exists $X\subseteq E(M_2)$ such that $r_{M_2}(X)\le i$ and $M_2/X$ has a $\DG((is+1)n_2,\Gamma)$-restriction $N_2$.
Note that $(i,M_2, X)=(s-1,M_1,C)$ is a candidate since $|C|\le s-1$, so this choice is well-defined. 
We consider two cases depending on whether $i=0$. 

Suppose that $i>0$ and let $Y_1,Y_2,\dots,Y_s,Z$ be mutually skew sets in $N_2$ so that $N_2|Y_j\cong \DG(n_2,\Gamma)$ for each $j\in [s]$ and $N_2|Z\cong \DG(((i-1)s+1)n_2,\Gamma)$; these sets can be chosen to be the closures in $N_2$ of disjoint subsets of a frame for $N_2$. 
If $M_2|Y_j=N_2|Y_j$ for some $j\in [s]$, then $M_2$ has a $\DG(n_2,\Gamma)$-restriction, which contradicts that $i>0$ and $i$ is minimal.
Thus, $M_2|Y_j\ne N_2|Y_j$ for each $j$, implying that $r_{M_2/Y_j}(X)\le r_{M_2}(X)-1\le i-1$ for each $j$.

Let $Y=Y_1\cup\dots\cup Y_s$ and let $J$ be a maximal subset of $Y$ such that $\elem(M_2/J)>q(r(M_2/J))$. 
Let $M_3=M_2/J$.
If $Y_j\subseteq J$ for some $j$, then $r_{M_3}(X)\le i-1$ and $(M_3/X)|Z=N_2|Z\cong \DG((i-1)s+1)n_2,\Gamma)$, contradicting the minimality of $i$.
Therefore, $Y-J$ contains a transversal $T$ of $(Y_1,\dots,Y_s)$.
Note that $T$ is an $s$-element independent set of $N_2/J$ and therefore of $M_2/J=M_3$. Moreover, by the maximality of $J$, each $e\in T$ satisfies $\elem(M_3)-\elem(M_3/e)>q(r(M_3))-q(r(M_3)-1)$.
Since $r(M_3)\ge r(N_2|Z)\ge n_2-1\ge r$, (b) holds for $M_3$ and $T$.

If $i=0$, then $N_2$ is a $\DG(n_2,\Gamma)$-restriction of $M_2$.
Let $M_4$ be a minimal minor of $M_2$ such that $\elem(M_4)>q(r(M_4))$ and $N_2$ is a restriction of $M_4$.
If $N_2$ is spanning in $M_4$ then (a) holds.
Otherwise, by minimality we have $\elem(M_4|\cl_{M_4}(E(N_2))\le q(r(N_2))$, so since $r(M_4)\ge n_2$ we have 
\begin{align*}
\elem(M_4\del \cl_{M_4}(E(N_2)))&>q(r(M_4))-q(r(N_2))\ge \ell^s.
\end{align*}
Therefore, there is an $s$-element independent set $S$ of $M_4$ which is disjoint from $\cl_{M_4}(E(N_2))$. 
Since $N_2$ is a restriction of $M_4/e$ for each $e\in S$, it follows from the minor-minimality of $M_4$ that $M_4$ and $S$ satisfy (b).
\end{proof}

We now combine Lemmas \ref{define nu} and \ref{spanning Dowling} to prove a result which is vital for the proof of Theorem \ref{main corollary}.

\begin{theorem}\label{new reduction}
There is a function $r_{\ref{new reduction}}\colon \bR^6\to \bZ$ so that for all integers $\ell,k\ge 2$ and $r,s\ge 1$ and any real polynomial $p(x)=ax^2+bx+c$ with $a>0$, if $M\in\cU(\ell)$ satisfies $r(M)\ge r_{\ref{new reduction}}(a,b,c,\ell,r,s)$ and $\elem(M)>p(r(M))$, then $M$ has a minor $N$ with $\elem(N)>p(r(N))$ and $r(N)\ge r$ such that either
\begin{enumerate}[(1)]
\item $N$ has a spanning clique restriction, or
\item $N$ is vertically $s$-connected and has an $s$-element independent set $S$ so that $\elem(N)-\elem(N/e)>p(r(N))-p(r(N)-1)$ for each $e\in S$.
\end{enumerate}
\end{theorem}
\begin{proof}
We first define the function $r_{\ref{new reduction}}$.
Let $\nu=\nu_{\ref{define nu}}(a,b,c,\ell,\max(r,s),r,s)$, and define $\hat r_1$ to be an integer so that 
$$(2s+1)a(x+y)+s(\nu+b)+c-as^2\le 2axy$$
and $p(x-s)\le p(x-s+1)$ for all real $x,y\ge \hat r_1$. 
Let $f$ be a function which takes in an integer $m$ and outputs an integer $f(m)\ge \max(r,2m,2\hat r_1)$ such that $p(x)-p(x-1)\ge ax+\ell^{\max(m,\hat r_1)}$ for all real $x\ge f(m)$.
Define $r_{\lceil\nu/a\rceil}=1$, and for each $i\in \{0,1,2,\dots,\lceil\nu/a\rceil-1\}$ recursively define 
$r_i$ to be an integer so that $p(x)>\alpha_{\ref{clique minor}}(\ell,n_{\ref{spanning Dowling}}(a,b,c,\ell,f(r_{i+1}),s))\cdot x$ for all real $x\ge r_i$.
Finally, define $r_{\ref{new reduction}}(a,b,c,\ell,r,s)=r_0$.

Let $M\in \cU(\ell)$ such that $r(M)\ge r_0$ and $\elem(M)>p(r(M))$. 
We may assume that $M$ has no rank-$\max(r,s)$ projective-geometry minor; otherwise outcome (2) holds.
Let $\cM$ denote the class of minors of $M$.
We may assume that $h_{\cM}(n)\le p(n)+\nu n$ for all $n\ge 1$, or else (2) holds by Lemma \ref{define nu}.
The following claim essentially finds some $\nu'$ so that the coefficient of the linear term of $h_{\cM}(n)$ is in the interval $[\nu'+b-a,\nu'+b+a]$.

\begin{claim}\label{nu''}
There is some $0\le \nu'< \nu$ and $i\ge 0$ such that $h_{\cM}(n)>p(n)+\nu' n$ for some $n\ge r_i$, and $h_{\cM}(n)\le p(n)+(\nu'+a)n$ for all $n\ge r_{i+1}$. 
\end{claim}
\begin{proof}
We will break up the real interval $[0,\nu]$ into subintervals of size $a$.
Define $\nu_i=ai$ for all $i\in \{0,1,2,\dots,\lceil\frac{\nu}{a}\rceil\}$.
Let $i\ge 0$ be minimal so that $h_{\cM}(n)\le p(n)+\nu_{i+1}n$ for all $n\ge r_{i+1}$. 
This choice of $i$ is well-defined, because $i=\lceil \nu/a\rceil-1$ is a valid choice since $\nu_{\lceil \nu/a\rceil}\ge \nu$ and $h_{\cM}(n)\le p(n)+\nu n$ for all $n\ge 1=r_{\lceil \nu/a\rceil}$. 

If $i>0$, then $h_{\cM}(n)>p(n)+\nu_i$ for some $n\ge r_i$ by the minimality of $i$.
If $i=0$, then $M$ certifies that $h_{\cM}(n)>p(n)$ for some $n\ge r_0$.
Thus, there is some $i\ge 0$ such that $h_{\cM}(n)>p(n)+\nu_i n$ for some $n\ge r_i$, and $h_{\cM}(n)\le p(n)+\nu_{i+1}n$ for all $n\ge r_{i+1}$. 
Since $\nu_i+a=\nu_{i+1}$, we may choose $\nu'=\nu_i$.
Note that $\nu_i=ai< \nu$ since $i\le \lceil \frac{\nu}{a}\rceil-1$.
\end{proof}

By  \ref{nu''}, $M$ has a minor $M_1$ with $r(M_1)\ge r_i$ and $\elem(M_1)>p(r(M_1))+\nu' r(M_1)$. 
Since $r(M_1)\ge r_i$, we have $p(r(M_1))>\alpha_{\ref{clique minor}}(\ell,n_{\ref{spanning Dowling}}(a,b,c,\ell,f(r_{i+1}),s))\cdot r(M_1)$, so $M_1$ has a 
$\DG(n_{\ref{spanning Dowling}}(a,b,c,\ell,f(r_{i+1}),s),\{1\})$-minor by Theorem \ref{clique minor}.
Then by Lemma \ref{spanning Dowling} with $r=f(r_{i+1})$ and $q=p+\nu'$, $M_1$ has a minor $N$ such that $r(N)\ge f(r_{i+1})$ and $\elem(N)>p(r(N))+\nu'r(N)$, and $N$ either has a spanning clique restriction or an $s$-element independent set so that $\elem(N)-\elem(N/e)>p(r(N))-p(r(N)-1)+\nu'r(N)$ for each $e\in S$. 
We may assume that $N$ is simple.
Since $f(r_{i+1})\ge r$ and $\nu'\ge 0$ we may assume that $N$ is not vertically $s$-connected, or else the theorem holds.

Let $(A,B)$ be a partition of $E(N)$ so that $r_N(A)\le r_N(B)< r(N)$ and $r_N(A)+r_N(B)-r(N)<s-1$. 
Let $r_N=r(N)$ and $r_A=r_N(A)$ and $r_B=r_N(B)$.
We first show that $r_A\ge \max(\hat r_1,r_{i+1})$.
If not, then $r_B\ge r_N-r_A\ge \max(r_{i+1},\hat r_1)$, using that $r_N\ge f(r_{i+1})\ge \max(2r_{i+1},2\hat r_1)$. 
\setcounter{equation}{0}
Also, 
\begin{align}
|B|=|N|-|A|&>p(r_N)+\nu'r_N-\ell^{\max(\hat r_1,r_{i+1})}\\
&\ge p(r_N-1)+(\nu'+a)r_N\\
&\ge p(r_B)+(\nu'+a)r_B.
\end{align}
Line (1) holds because $r_A<\max(\hat r_1,r_{i+1})$ and $M\in \cU(\ell)$, and line (2) holds because $r_N\ge f(r_{i+1})$.
Line (3) holds because $r_B\ge \hat r_1$, so $p(r_B)\le p(r_N-1)$ since $r_B\le r_N-1$.
But then $r_B\ge r_{i+1}$ and $|B|>p(r_B)+(\nu'+a)r_B$, which contradicts  \ref{nu''} and the choice of $\nu'$.
Thus, $r_B\ge r_A\ge \max(\hat r_1,r_{i+1})$. 
Then 
\begin{align*}
p(r_A+r_B-s)+\nu'(r_A+r_B-s)&\le p(r_N)+\nu' r_N\\
&<|A|+|B|\\
&\le p(r_A)+p(r_B)+(\nu'+a)(r_A+r_B),
\end{align*}
where the first inequality holds because $r_A+r_B-s\le r_N$ and $p(x-s)\le p(x-s+1)$ for all $x\ge \hat r_1$, and the last inequality holds by  \ref{nu''} because $r_B\ge r_A\ge r_{i+1}$.
Expanding $p(x)=ax^2+bx+c$ and simplifying, we have 
$$(2s+1)a(r_A+r_B)+s(\nu'+b)+c-as^2>2ar_Ar_B,$$
which contradicts that $r_A\ge \hat r_1$, since $\nu'< \nu$.
\end{proof}
\end{subsection}


\begin{subsection}{Porcupines and Stacks} \label{porcs and stacks}
We now define three structures that arise from the second outcome of Theorem \ref{new reduction}. 
The first structure is a collection of bounded-size restrictions, each of which is not a restriction of any $\DG(k,\Gamma)$ with $|\Gamma|<t$. 
For integers $b\ge 2$ and $h\ge 1$ and a collection $\mathcal O$  of matroids, a matroid $M$ is an \emph{$(\mathcal O,b,h)$-stack} if there are disjoint sets $P_1,P_2,\dots,P_h\subseteq E(M)$ such that 
\begin{itemize}
\item $\cup_{i\in [h]} P_i$ spans $M$, and

\item for each $i\in [h]$, the matroid $(M/(P_1\cup\dots\cup P_{i-1}))|P_i$ has rank at most $b$ and is not in $\mathcal O$. 
\end{itemize}  
Stacks were used to find the extremal functions for exponentially dense minor-closed classes in \cite{Dense Exp} and \cite{Densest Exp} with $\cO$ as the class of $\GF(q)$-representable matroids; our definition generalizes the original definition from \cite{Hales-Jewett}.
We say that a matroid $M$ is an \emph{$\cO$-stack} if there are integers $b\ge 2$ and $h\ge 1$ so that $M$ is an $(\cO,b,h)$-stack.

We will always take $\mathcal O$ to be $\cF\cap \cU(t)$ where $\cF$ is the class of frame matroids; thus $\cF\cap \cU(t)$ is the class of frame matroids with no $U_{t+2,2}$-minor.
Stacks are helpful because every matroid with a spanning clique restriction and a large enough $(\cF\cap \cU(t))$-stack restriction either has a $\DG(k,\Gamma)$-minor with $|\Gamma|\ge t$, or is not a bounded distance from a frame matroid, as we prove in Section \ref{main proof chap}. 


The second structure is a collection of nearly skew spikes with a common tip, which we define more precisely as
a collection $\cS$ of mutually skew sets in $M/e$ so that for each $R\in\cS$, the matroid $M|(R\cup \{e\})$ is a spike with tip $e$. 
We prove in Section \ref{main proof chap} that such a restriction in the span of a clique admits an $\LG^+(n,\Gamma)$-minor for some nontrivial group $\Gamma$.
This structure also has the useful property that if $|\cS|\ge 2$, then $M|(\cup\cS\cup\{e\})$ is not a frame matroid, which we prove later in this section.


The third structure is a large independent set for which each element is the tip of many large spike restrictions. More precisely, for each integer $g\ge 3$ a \emph{$g$-preporcupine} is a matroid $P$ with an element $f$ such that 
\begin{itemize}
\item each line of $P$ through $f$ has at least three points, and 

\item $\si(P/f)$ has girth at least $g$. 
\end{itemize}
(Recall that the girth of a matroid is the size of its smallest circuit.)
We say that $f$ is the \emph{tip} of $P$, and we write $d(P)$ for $r^*(\si(P/f))$, the corank of $\si(P/f)$.
(We comment that while there may be multiple choices for an element $f$ certifying that $P$ is a $g$-preporcupine, we take the convention that $f$ is fixed, and we explicitly name $f$ whenever possible to avoid confusion.
In other words, a $g$-preporcupine $P$ with tip $f$ is technically an ordered pair $(P,f)$, but for readability we implicitly assume that the choice of $f$ is fixed.)

It will often be convenient to work with the following more specific matroid: a \emph{$g$-porcupine} $P$ is a simple $g$-preporcupine with tip $f$ such that each line of $P$ through $f$ has exactly three points, and $\si(P/f)$ has no coloops.
Clearly every $g$-preporcupine has a $g$-porcupine restriction, but porcupines are more restricted; every $g$-porcupine $P$ with $d(P)=0$ consists of a single element, and every $g$-porcupine $P$ with $d(P)=1$ is a spike. 

More generally, a matroid $P$ is a \emph{(pre)porcupine} if there is an integer $g\ge 3$ so that $P$ is a $g$-(pre)porcupine.
Porcupines are helpful because a large enough collection of porcupines with independent tips is not a bounded distance from a frame matroid, as we prove in Section \ref{main proof chap}.
\end{subsection}


\begin{subsection}{Nearly Skew Spikes with a Common Tip} \label{spikes sec}
We now prove a result (Proposition \ref{skew spikes}) which provides sufficient conditions for a matroid to have a collection of nearly skew small spike restrictions with a common tip.
We start by introducing some notation which we will use for the remainder of this section.

For a matroid $M$ and an element $f\in E(M)$, let $\delta(M,f)=\elem(M)-\elem(M/f)$, and  
let $\cL_M(f)$ denote the set of long lines of $M$ through $f$.
It is not hard to show that
$$\delta(M,f)=1+\sum_{L\in \cL_M(f)}\Big(\elem(M|L)-2\Big),$$
so $\delta(M,f)=\delta\big(M|(\cup \cL_M(f)),f\big)$.
We use this to prove the following lemma, which gives a formula for $\delta(M,f)-\delta(M/C,f)$ in terms of the long lines of $M/C$ through $f$.


\begin{lemma} \label{criticality formula}
If $M$ is a matroid with $C\subseteq E(M)$, and $f\in E(M)-\cl_M(C)$ so that $\cl_M(C)\cup \{f\}$ is a flat of $M$, then 
$$\delta(M,f)-\delta(M/C,f)=\sum_{L\in \cL_{M/C}(f)}\Big(\delta(M|L,f)-\elem\big((M/C)|L\big)+1\Big).$$
\end{lemma}
\begin{proof}
We may assume that $M$ is simple, by the definition of $\delta(M,f)$. 
Since $\cl_M(C)\cup\{f\}$ is a flat of $M$, each long line of $M$ through $f$ is disjoint from $C$ and is therefore contained in a long line of $M/C$ through $f$.
Since $M$ is simple and $\cl_M(C)\cup\{f\}$ is a flat of $M$, for each $L,L'\in \cL_{M/C}(f)$ with $L\ne L'$ we have $L\cap L'=\{f\}$. 
Thus, each long line of $M$ through $f$ is contained in precisely one set in $\cL_{M/C}(f)$, so
\begin{align*}
\delta(M,f)&=1+\sum_{L\in \cL_{M/C}(f)}\Big(\elem(M|L)-\elem\big((M|L)/f\big)-1\Big)\\
&=1+\sum_{L\in \cL_{M/C}(f)}\Big(\delta(M|L,f)-1\Big).
\end{align*}
Combining this with the fact that $$\delta(M/C,f)=1+\sum_{L\in \cL_{M/C}(f)}\big(\elem\big((M/C)|L\big)-2\big)$$ gives the desired result.
\end{proof}


We now use Lemma \ref{criticality formula} to find a large collection of small spikes with common tip $e$ such that each spike spans some other element $f$.

\begin{lemma} \label{criticality drop}
For all integers $b\ge 1$ and $\ell\ge 2$, if $\{e,f\}$ is a $2$-element independent set of a matroid $M\in \cU(\ell)$, and $\delta(M,f)-\delta(M/e,f)\ge \ell^{2b+4}$, then there is a collection $\cS$ of $b$ mutually skew sets in $M/\{e,f\}$ so that for each $S\in \cS$, the matroid $M|(S\cup \{e\})$ is a spike of rank at most four with tip $e$.
\end{lemma}
\begin{proof}
We may assume that $M$ is simple. Let $M_1=M\del (\cl_M(\{e,f\})-\{e,f\})$. 
Then $\{e,f\}$ is a flat of $M_1$, and $\delta(M_1,f)-\delta(M_1/e,f)\ge \ell^{2b+3}$. 
Let $\cL=\cL_{M_1/e}(f)$. 
We first prove a claim to help find long lines through $e$.

\begin{claim} \label{long lines}
For each $L\in \cL$, if $\delta(M_1|L,f)-\elem\big((M_1/e)|L\big)+1>0$, then $M_1|\big((L\cup\{e\})-\{f\}\big)$ has at least two long lines through $e$. 
\end{claim}
\begin{proof}
Assume for a contradiction that $M_1|(L\cup\{e\})$ has at most one long line through $e$, and let $m\ge 2$ be the length of the longest line of $M_1|(L\cup\{e\})$ through $e$. 
Since $\{e,f\}$ is a flat of $M_1$, each line of $M_1|(L\cup\{e\})$ through $e$ contains at most one element of each line of $M_1|L$ through $f$.
This implies that there are at least $m-1$ lines of $M_1|L$ through $f$, so $\elem\big((M_1|L)/f\big)\ge m-1$, and thus $\delta(M_1|L,f)\le \elem(M_1|L)-(m-1)$. 
Also, $\elem\big((M_1/e)|L\big)=\elem(M_1|L)-(m-2)$ since there is only one long line of $M_1|(L\cup\{e\})$ through $e$. 
Combining these facts gives 
\begin{align*}
\delta(M_1|L,f)-\elem\big((M_1/e)|L\big)+1&\le \big(\elem(M_1|L)-(m-1)\big)\\
&\hspace{0.5cm}-\big(\elem(M_1|L)-(m-2)\big)+1\\
&=0,
\end{align*}
a contradiction.
\end{proof}

Each $L\in \cL$ satisfies 
$$\delta(M_1|L,f)-\elem\big((M_1/e)|L\big)+1\le \elem(M_1|L)\le \ell^3$$ by Theorem \ref{l-Kung}, since $r_{M_1}(L)\le 3$ and $M_1\in\cU(\ell)$. 
Since $\delta(M_1,f)-\delta(M_1/e,f)\ge \ell^3\ell^{2b}$, by Lemma \ref{criticality formula} with $C=\{e\}$ there are at least $\ell^{2b}$ sets $L\in \cL$ for which $\delta(M_1|L,f)-\elem\big((M_1/e)|L\big)+1>0$. 
Since each $L\in\cL$ is a point of $M_1/\{e,f\}$ and $M_1\in \cU(\ell)$, there is some $\cL'\subseteq \cL$ such that $|\cL'|=2b$, $r_{M_1/\{e,f\}}(\cup \cL')=2b$, and each $L\in\cL'$ satisfies $\delta(M_1|L,f)-\elem\big((M_1/e)|L\big)+1>0$. 

For each $L\in \cL'$, the matroid $M_1|((L\cup\{e\})-\{f\})$ contains at least two long lines through $e$, by  \ref{long lines}. 
Let $\cL'=\{L_1,\dots,L_{2b}\}$. 
Then each $i\in [b]$ satisfies $r_{M_1}(L_{2i-1}\cup L_{2i})= 4$, and there are at least four long lines through $e$ in $M_1|((L_{2i-1}\cup L_{2i}\cup\{e\})-\{f\})$. 
For each $i\in [b]$, let $P_i$ denote the union of the long lines through $e$ in $M_1|((L_{2i-1}\cup L_{2i}\cup\{e\})-\{f\})$.  
Then $r((M_1|P_i)/e)\le 3$ and $\elem((M_1|P_i)/e)\ge 4$, so $\si((M_1|P_i)/e)$ contains a circuit. 
Since $M_1|P_i$ is simple and each parallel class of $(M_1|P_i)/e$ has size at least two, for each $i\in [b]$ there is some $S_i\subseteq P_i-\{e\}$ so that $M_1|(S_i\cup\{e\})$ is a spike with tip $e$. 
Since $r_{M_1/\{e,f\}}((\cup\cL')-\{f\})=2b$, the sets $\{S_i\colon i\in [b]\}$ are mutually skew in $M_1/\{e,f\}$. 
\end{proof}


If Lemma \ref{criticality drop} applies for enough elements $f$, then we can find many nearly skew small spikes with common tip $e$.
This is one of the three structures which arise from the second outcome of Theorem \ref{new reduction}.

\begin{proposition} \label{skew spikes}
Let $m\ge 1$ and $\ell\ge 2$ be integers, and let $M\in \cU(\ell)$ be a matroid with a $3m$-element  independent set $X$.
Let $e\in E(M)-X$ such that $M|(X\cup\{e\})$ is simple, and each $f\in X$ satisfies $\delta(M,f)-\delta(M/e,f)\ge \ell^{6m}$.
 Then there is a collection $\mathcal S$ of $m$ mutually skew sets in $M/e$ so that for each $S\in \cS$, the matroid $M|(S\cup \{e\})$ is a spike of rank at most four with tip $e$.
\end{proposition}
\begin{proof}
Let $\cS$ be a maximum-size collection of mutually skew sets in $M/e$ so that for each $S\in \cS$, the matroid $M|(S\cup \{e\})$ is a spike of rank at most four with tip $e$, and assume for a contradiction that $k<m$. 
Then $r_M(\cup \cS)\le 3m-2$, so there is some $f\in X-\cl_M(\cup \cS)$. Let $B\cup\{e\}$ be a basis for $M|(\cup \cS)$, so $|B|\le 3m-3$ and $B\cup\{e,f\}$ is independent. 

By Lemma \ref{criticality drop} with $b=3m-2$, there is a collection $\cP$ of $3m-2$ mutually skew sets in $M/\{e,f\}$ so that for each $P\in\cP$, the matroid $M|(P\cup \{e\})$ is a spike of rank at most four with tip $e$.
Recall that for sets $A_1$ and $A_2$ of a matroid $N$, we write $\sqcap_N(A_1,A_2)$ for $r_N(A_1)+r_N(A_2)-r_N(A_1\cup A_2)$.
Since the sets in $\cP$ are mutually skew in $M/\{e,f\}$, we can show using submodularity of the rank function of $M$ and induction on $|\mathcal P|$ that
$$\sqcap_{M/\{e,f\}}(B, \cup\cP)\ge \sum_{P\in \cP}\sqcap_{M/\{e,f\}}(B,P).$$
Thus, there is some $P\in\cP$ such that $\sqcap_{M/\{e,f\}}(B,P)=0$, or else 
$\sqcap_{M/\{e,f\}}(B,\cup \cP)\ge |\cP|>|B|$, a contradiction.
Then 
\begin{align*}
\sqcap_{M/e}(\cup \cS,P)&=\sqcap_{M/e}(B,P)\le \sqcap_{M/\{e,f\}}(B,P)=0,
\end{align*}
because $B$ is a basis for $(M/e)|(\cup\cS)$ and $f\notin\cl_{M/e}(B)$. But then $\cS\cup\{P\}$ contradicts the maximality of $|\cS|$.
\end{proof}
\end{subsection}


\begin{subsection}{Star-Partitions} \label{P-partitions sec}
Recall that for any integer $g\ge 3$, a \emph{$g$-preporcupine} is a matroid $P$ with an element $f$ so that each line of $P$ through $f$ has at least three points, and  $\si(P/f)$ has girth at least $g$. 
A \emph{$g$-porcupine} $P$ is a simple $g$-preporcupine such that each line of $P$ through $f$ has exactly three points, and $\si(P/f)$ has no coloops.

We now explore the notion that preporcupines are incompatible with frame matroids.
We show that frame matroids do not have spike restrictions of rank at least five; this implies that frame matroids can only have preporcupine restrictions with very specific structure.
This is described by Proposition \ref{g-star}, which roughly says that every preporcupine either has a bounded-size restriction which is not a frame matroid, or has structure similar to that of a frame matroid. 

The basic idea is that if a matroid $M$ is both a frame matroid and a preporcupine with tip $f$, then $\si(M/f)$ is a restriction of a \emph{star}, which is a matroid $S$ with a basis $B\cup \{t\}$ such that $E(S)=\cup_{b\in B}\cl_S(\{t,b\})$.
The element $t$ is the \emph{tip} of $S$.


\begin{lemma} \label{star}
If $M$ is both a frame matroid and a preporcupine with tip $f$ such that $\si(M/f)$ contains a circuit, then all but at most one line through $f$ has length three, and $\si(M/f)$ is a restriction of a star. Moreover, if $M$ has a line $L$ of length at least four through $f$ then for each $e\in L-\cl_M(\{f\})$, the matroid $\si(M/\{e,f\})$ contains no circuit.
\end{lemma}
\begin{proof}
We may assume that $M$ is simple.
Let $N$ be a matroid framed by $B=\{b_1,b_2,\dots,b_r\}$ so that $N\del B=M$.
We may assume without loss of generality that $f\in\cl_{N}(\{b_1,b_2\})$, by relabeling $B$. Since $\si(M/f)$ contains a circuit, $f$ is not parallel to $b_1$ or $b_2$ by Lemma \ref{frame obvious} (iii). Then each element of $M$ is in $\cl_N(\{b_1,b_i\})$ for some $i\ge 2$ or $\cl_N(\{b_2,b_i\})$ for some $i\ge 3$, or else that element is not on a long line through $f$. 

Let $(M/f)|T$ be a simplification of $M/f$. Each element of $T$ is in $\cl_{N/f}(\{b_2,b_i\})$ for some $i\ge 3$. Since $\{b_2,b_3,\dots,b_r\}$ is independent in $N/f$, the matroid $(M/f)|T$ is a restriction of a star with tip $b_2$. If $M$ has a line $L$ of length at least four through $f$, then $L\subseteq \cl_{N}(\{b_1,b_2\})$ by Lemma \ref{frame obvious} (i). 
Then for each $e\in L-\{f\}$, $f$ is parallel to $b_1$ or $b_2$ in $N/e$, so $\si(M/\{e,f\})$ contains no circuit by Lemma \ref{frame obvious} (iii). 
\end{proof}


To make use of Lemma \ref{star} we will need some properties of restrictions of stars.
Stars have the property that each circuit is contained in the union of at most two lines through the tip.
We would also like to describe this property for restrictions of stars, but it is trickier if the restriction does not contain the tip.
To deal with this, we define a \emph{star-partition} of a matroid $M$ to be a pair $(X,\cL)$ such that 
\begin{itemize}
\item $r_M(X)\le 1$, 
\item $\cL$ partitions $E(M)-X$,
\item $L\cup X$ is a flat of $M$ of rank at most two for each $L\in \cL$, and 
\item $r_{M}(L\cup L')\le 3$ for all distinct $L,L'\in \cL$.
\end{itemize}
The idea is that any matroid with a star-partition $(X,\cL)$ so that each small circuit is contained in $L\cup L'\cup X$ for some $L,L'\in\cL$ is `star-like'.
The set $X$ is either empty or contains the tip of the star, and each set $L$ is a line through the tip (excluding the tip).
We cannot define a star-partition to actually be a partition, since it may be the case that $X$ is empty.

Lemmas \ref{starfact1} and \ref{starfact2} combine to show that any restriction of a star has a star-partition $(X,\cL)$ so that each circuit is contained in $L\cup L'\cup X$ for some $L,L'\in\cL$.
We first show that stars have this property.


\begin{lemma} \label{transversals}
Let $S$ be a simple star with tip $t$ and basis $B\cup\{t\}$, and let $$\cL=\{L\subseteq E(S)-\{t\}\colon L\cup\{t\} \textrm{ is a line of $S$}\}.$$
Then $(\{t\},\cL)$ is a star-partition of $S$ so that for each circuit $C$ of $S$ there are sets $L,L'\in \cL$ so that $C\subseteq \{t\}\cup L\cup L'$.
In particular, $S$ has no circuit of size at least five.
\end{lemma}
\begin{proof}
Clearly $(\{t\},\cL)$ is a star-partition of $S$, since each line in $\cL$ spans $t$.
Let $C$ be a circuit of $S$, and note that $C-\{t\}$ is a union of circuits of $S/t$.
Since $\si(S/t)$ is independent, the only circuits of $S/t$ are parallel pairs, so $C-\{t\}$ is a union of parallel pairs of $S/t$.
Then since $r_{S/t}(C)\ge |C|-2$, the set $C-\{t\}$ intersects at most two parallel classes of $S/t$.
\end{proof}

Lemma \ref{transversals} and Lemma \ref{star} together imply that spikes of rank at least five are not frame matroids, since every spike is a preporcupine.


If a restriction of a star contains the tip of the star, then we can easily find a star-partition.

\begin{lemma} \label{starfact1}
Let $M$ be a simple restriction of a star so that $x\in E(M)$ is on two long lines of $M$, and let  $$\cL=\{L\subseteq E(M)-\{x\}\colon L\cup \{x\} \textrm{ is a line of } M\}.$$ 
Then $(\{x\},\cL)$ is a star-partition of $M$ such that for each circuit $C$ of $M$ there are sets $L,L'\in\cL$ so that $C\subseteq \{x\}\cup L\cup L'$. 
\end{lemma}
\begin{proof}
Let $S$ be a star with tip $t$ such that $S|A=M$. 
If $x\in A$ is on two long lines of $S|A$, then $x=t$ since each element of $S$ other than $t$ is on at most one long line of $S$.
Thus, the result holds by applying Lemma \ref{transversals} to $S$.
\end{proof}


If a restriction of a star does not contain the tip of the star then we can still find a star-partition.
Moreover, for each circuit there are two elements of that circuit which define a star-partition.

\begin{lemma} \label{starfact2}
Let $M$ be a simple restriction of a star such that no element of $M$ is on two long lines of $M$, and let $C$ be a circuit of $M$. 
Then for each $e\in C$ there is some  $e'\in C-\{e\}$ such that 
$$\Big(\varnothing, \{\cl_M(\{e,e'\})\}\cup\{P\colon P \textrm{ is a parallel class of $M/\{e,e'\}$}\}\Big)$$
 is a star-partition of $M$, and for each circuit $C'$ of $M$ there are sets  $L,L'\in\{\cl_M(\{e,e'\})\}\cup\{P\colon P \textrm{ is a parallel class of $M/\{e,e'\}$}\}$ so that $C'\subseteq L\cup L'$. 
\end{lemma}
\begin{proof}
Let $S$ be a star with tip $\{t\}$ such that $S|A=M$, and let $e\in C$. 
Since $C-\{t\}$ is a union of circuits of $S/t$, there is some $e'\in C$ so that $\{e,e'\}$ is a parallel pair of $S/t$.
Let $\cP$ denote the collection of parallel classes of $M/\{e,e'\}$. 
If $(\varnothing, \{\cl_M(\{e,e'\})\}\cup\cP)$ is not a star-partition of $M$, then there is some $P\in \cP$ which is not a flat of $M$, so $\cl_M(P)\cap \cl_M(\{e,e'\})\ne \varnothing$.
But then there is an element of $M$ which is on two long lines of $M$, a contradiction. 

Since $t\in \cl_S(\{e,e'\})$ and $\si(S/t)$ is independent, each set in $\{\cl_M(\{e,e'\})\}\cup\cP$ is contained in a set in 
$$\{L\subseteq E(S)-\{t\}\colon L\cup\{t\} \textrm{ is a line of $S$}\}.$$
Thus, for each circuit $C'$ of $M$ there exist $L,L'\in\{\cl_M(\{e,e'\})\}\cup\cP$ so that $C'\subseteq L\cup L'$, by applying Lemma \ref{transversals} to $S$.
\end{proof}


We now use these properties of stars to prove an extension of Lemma \ref{star} to matroids for which each bounded-rank restriction is a frame matroid. 
The first outcome will help build a stack for which each piece is not a frame matroid, and the second outcome will help find a $g$-porcupine.

\begin{proposition} \label{g-star}
Let $M$ be a preporcupine with tip $f$, and let $(M/f)|T$ be a simplification of $M/f$. 
If $(M/f)|T$ has a circuit of size at most four, then for each integer $g\ge 4$, either 
\begin{enumerate}[(1)]
\item $M$ has a restriction of rank at most $3g$ which is not a frame matroid, or 

\item $(M/f)|T$ has a star-partition $(X,\cL)$ so that each line of $M$ through $f$ and an element of $T-X$ has length three, and for each circuit $C'$ of $(M/f)|T$ of size less than $g$, there are sets $L,L'\in \cL$ so that $C'\subseteq (X\cup L\cup L')$.
\end{enumerate}
\end{proposition}
\begin{proof}
Suppose for a contradiction that (1) and (2) do not hold for $M$.
Then $M$ satisfies the following property:
\begin{enumerate}[(1')]
\item For each set $X\subseteq E(M)$ of rank at most $3g$, the matroid $M|X$ is a frame matroid.
\end{enumerate}
We may assume that $M$ is simple. 
Let $M'=(M/f)|T$ and let $C$ be a circuit of $M'$ of size at most four.


\begin{claim} \label{length three}
Each line of $M$ through $f$ has length three.
\end{claim}
\begin{proof}
Suppose for a contradiction that a line $L_1$ of $M$ through $f$ has length at least four, let $x\in L_1\cap T$. 
Let $X=\{x\}$ and $\cL=\{L\subseteq E(M')-\{x\}\colon L\cup\{x\} \textrm{ is a line of $M'$}\}$.
Then $(X,\cL)$ is a star-partition of $M'$, since each line in $\cL$ spans $x$.
We show that $(X,\cL)$ satisfies (2).
If there is some $e\in T-X$ such that the line $L_2$ of $M$ through $e$ and $f$ has length at least four, then $\cl_M(C\cup\{f,x,e\})$ is not a frame matroid by Lemma \ref{star}.
Since $r_M(C\cup\{f,x,e\})\le 7\le 3g$, this contradicts (1').
Thus, each line of $M$ through $f$ and an element of $T-X$ has length three.

If $M'$ has a circuit $C'$ of size less than $g$ which is not contained in $X\cup L\cup L'$ for any $L,L'\in\cL$, then $C'-\{x\}$ intersects at least three parallel classes of $M'/x$.
Since $r_{M'/x}(C'-\{x\})\ge |C'|-2$, some parallel class of $(M'/x)|(C'-\{x\})$ has size one.
Then since $C'-\{x\}$ is a union of circuits of $M'/x$, the matroid $(M'/x)|(C'-\{x\})$ has a circuit which intersects at least three parallel classes of $M'/x$, so $\si(M'/x)$ has a circuit of size less than $g$.
Since $|L_1|\ge 4$ and $x\in L_1$, the matroid $M|\cl_M(C\cup C'\cup \{f,x\})$ is not frame by Lemma \ref{star}.
Since $r_M(C\cup C'\cup \{f,x\})\le g+6\le 3g$, this contradicts (1').
Therefore, $(X,\cL)$ is a star-partition of $M'$ which satisfies (2), a contradiction.
\end{proof}


The next claim allows us to apply Lemma \ref{starfact2}.

\begin{claim} \label{X is empty}
No element of $M'$ is on two long lines of $M'$.
\end{claim}
\begin{proof}
Suppose for a contradiction that $x\in E(M')$ is on two long lines of $M'$.
Let $X=\{x\}$ and $\cL=\{L\subseteq E(M')-\{x\}\colon L\cup\{x\} \textrm{ is a line of $M'$}\}$.
Then $(X,\cL)$ is a star-partition of $M'$, and there are distinct $L_1,L_2\in\cL$ such that $|L_1|\ge 2$ and $|L_2|\ge 2$.
By  \ref{length three}, each line of $M$ through $f$ and an element of $\cup\cL$ has length three.
Since $(X,\cL)$ does not satisfy (2) there is some circuit $C'$ of $M'$ of size less than $g$ which is not contained in $X\cup L\cup L'$ for any $L,L'\in\cL$.

Let $\cL'=\{L\in \cL\colon L\cap C'\ne \varnothing\}\cup \{L_1,L_2\}$, and note that $C'\subseteq \cup(\cL'\cup\{X\})$.
If $M'|(\cup (\cL'\cup \{X\}))$ is a restriction of a star, then each circuit is contained in $X\cup L\cup L'$ for some $L,L'\in\cL'$ by Lemma \ref{starfact1}.
But this contradicts the existence of $C'$, so $M'|(\cup (\cL'\cup \{X\}))$ is not a restriction of a star.
But then the set of elements of $M$ on a long line through $f$ and an element of $\{x\}\cup (\cup\cL')$ is not a frame matroid by Lemma \ref{star}.
Since $r_{M}(\{f,x\}\cup (\cup\cL'))\le g+2\le 3g$, this contradicts (1').
\end{proof}


Fix some element $e\in C$.
For each $e'\in C-\{e\}$, let $\cP_{e,e'}$ denote the collection of parallel classes of $M'/\{e,e'\}$.

\begin{claim} \label{stack}
There is some $e'\in C-\{e\}$ so that the pair $(\varnothing,\{\cl_{M'}(\{e,e'\})\}\cup\cP_{e,e'})$ is a star-partition of $M'$ such that for each circuit $C'$ of $M'$ of size less than $g$, there are $L,L'\in \{\cl_{M'}(\{e,e'\})\}\cup\cP_{e,e'}$ for which $C'\subseteq L\cup L'$. 
\end{claim}
\begin{proof}
Assume for a contradiction that the claim is false.
For each $e'\in C-\{e\}$, the pair $(\varnothing,\{\cl_{M'}(\{e,e'\})\}\cup\cP_{e,e'})$ is a star-partition of $M'$ by  \ref{X is empty}, since these sets are nonempty pairwise coplanar flats of rank at most two.
So for each $e'\in C-\{e\}$ there is a circuit $C_{e'}$ of $M'$ of size less than $g$ not contained in $L\cup L'$ for any $L,L'\in \{\cl_{M'}(\{e,e'\})\}\cup\cP_{e,e'}$.
Let $J=\cl_{M'}(C)\cup (\cup_{e'\in C-\{e\}}C_{e'})$, and for each $e'\in C-\{e\}$ let $\cP'_{e,e'}$ denote the collection of parallel classes of $(M'|J)/\{e,e'\}$.

If $M'|J$ is a restriction of a star, then by  \ref{X is empty} and Lemma \ref{starfact2} there is some $e'\in C-\{e\}$ for which the pair ($\{\cl_{M'}(\{e,e'\})\}\cup\cP'_{e,e'},\varnothing)$ is a star-partition of $M'|J$ such that each circuit of $M'|J$ is contained in $L\cup L'$ for some $L,L'\in \{\cl_{M'}(\{e,e'\})\}\cup\cP'_{e,e'}$.
This contradicts the existence of $C_{e'}$, so $M'|J$ is not a restriction of a star.
Then $\cl_M(J\cup\{f\})$ is not a frame matroid, by Lemma \ref{star}. Since $|C-\{e\}|\le 3$ we have $r_{M'}(J)\le 3(g-2)+3$.
But then $r_M(J\cup \{f\})\le 3g$, which contradicts (1'). 
\end{proof}

Let $e'$ be given by  \ref{stack}. Then by \ref{length three} and \ref{stack}, outcome (2) holds with $X=\varnothing$ and $\cL=\{\cl_{M'}(\{e,e'\})\}\cup \cP_{e,e'}$, a contradiction. 
\end{proof}
\end{subsection}


\begin{subsection}{Porcupines and Frame Matroids} \label{porcs and frames}
Recall that $\cF\cap \cU(t)$ is the class of frame matroids with no $U_{2,t+2}$-minor.
We now apply Proposition \ref{g-star} to a collection of preporcupines with independent tips, and show that we can either find a small matroid which is not in $\cF\cap \cU(t)$, or a large independent set for which each element is the tip of a $g$-porcupine.

Before proving this, we need two straightforward lemmas to help find a $U_{2,t+2}$-minor. The first deals with a rank-4 frame matroid.


\begin{lemma} \label{m=1}
For each integer $t\ge 2$, if $M$ is both a frame matroid and a porcupine with tip $f$ such that $\si(M/f)$ has rank three and is the disjoint union of a line of length $t$ and a line of length two, then $M$ has a $U_{2,t+2}$-minor. 
\end{lemma}
\begin{proof}
We may assume that $M$ is simple.
Let $N$ be a matroid framed by $B=\{b_1,b_2,b_3,b_4\}$ such that $N\del B=M$.
Let $(M/f)|T$ be a simplification of $M/f$. 
We may assume that $f\in\cl_{N}(\{b_1,b_2\})$, by relabeling $B$.
Note that $f$ is not parallel to $b_1$ or $b_2$, by Lemma \ref{frame obvious} (iii), since $f$ is the tip of a spike restriction of $M$.

Each element of $N/f$ is spanned by $\{b_2,b_3\}$ or $\{b_2,b_4\}$, or else that element is not on a long line of $M$ through $f$. 
Since $M$ is a porcupine, the matroid $(M/f)|T$ has no coloops, and therefore $\cl_{(N/f)|T}(\{b_2,b_3\})$ and $\cl_{(N/f)|T}(\{b_2,b_4\})$ each contain at least two elements in $E(M)$ which are not parallel to $b_2$ in $N/f$. 
Let $L_1$ and $L_2$ be the lines of $(M/f)|T$ of length $t$ and two, respectively, whose union is $T$.

If $t\ge 3$, assume without loss of generality that $L_1\subseteq \cl_{(N/f)|T}(\{b_2,b_3\})$, since $\{b_2,b_3\}$ and $\{b_2,b_4\}$ span the only long lines of $(N/f)|T$. 
Then $L_2\subseteq \cl_{(N/f)|T}(\{b_2,b_4\})$, and no element of $(N/f)|T$ is parallel to $b_2$ since $L_1$ and $L_2$ are flats of $(N/f)|T$. 
Then since $|L_1|= t$ and each element of $L_1$ is on a long line of $M\del b_2$ through $f$, there are $t$ long lines of $M$ through $f$ so that each contains one element in $\cl_N(\{b_1,b_3\})$ and one element in $\cl_N(\{b_2,b_3\})$.
Similarly, $\{b_1,b_4\}$ and $\{b_2,b_4\}$ each span two elements of $E(M)$ since $|L_2|=2$. 

Let $e,e'\in \cl_{N}(\{b_2,b_4\})\cap E(M)$ with $e\ne e'$, and let $x\in \cl_{N}(\{b_2,b_3\})\cap E(M)$. 
Let $Y=\cl_{N}(\{b_1,b_3\})\cap E(M)$, so $|Y|=t$.
Note that $\{f,b_1\}$ and $\{x,b_3\}$ are parallel pairs in $N/\{e,e'\}$ since $b_2\in\cl_{N}(\{e,e'\})$, and $f$ and $x$ are not parallel in $N/\{e,e'\}$. 
Thus, 
$\{f,x\}\subseteq \cl_{M/\{e,e'\}}(Y)-Y,$ so $(M/\{e,e'\})|(Y\cup\{f,x\})\cong U_{2,t+2}$, as desired.
If $t=2$ then $\{b_2,b_3\}$ and $\{b_2,b_4\}$ each span two points of $M/f$ which are not parallel to $b_2$, and the same argument applies.
\end{proof}

The second lemma deals with a rank-3 frame matroid.

\begin{lemma}\label{rank-3}
For each integer $t\ge 2$, if $M$ is both a frame matroid and a preporcupine with tip $f$ such that $\si(M/f)\cong U_{2,t+1}$, and $M$ has a line of length $t+1$ through $f$, then $M$ has a $U_{2,t+2}$-minor.
\end{lemma}
\begin{proof}
We may assume that $M$ is simple.
If $t=2$ then $M$ is a rank-3 spike, and if $M$ is binary then $M$ is isomorphic to the Fano plane, which is a contradiction since the Fano plane is not a frame matroid.
Thus, $M$ is not binary, so $M$ has a $U_{2,4}$-minor and the lemma holds for $t=2$.
Therefore, we may assume that $t\ge 3$.
Let $N$ be a matroid framed by $B=\{b_1,b_2,b_3\}$ so that $N\del B=M$.
Assume without loss of generality that $f\in\cl_{N}(\{b_1,b_2\})$, and note that $b_3$ is not parallel to an element of $E(M)$ by Lemma \ref{frame obvious} (iii) since $b_3$ is not on any long line with $f$.
Since $t\ge 3$ and $M$ has a line of length $t+1$ through $f$ we have $|E(M)\cap \cl_N(\{b_1,b_2\})|=t+1$ by Lemma \ref{frame obvious} (i).

If $b_i$ is not parallel to any element in $E(M)$, then for each $e\in \cl_M(\{b_i,b_3\})$ the matroid $M/e$ has a $U_{2,t+2}$-restriction. 
So $b_1$ and $b_2$ are each parallel in $N$ to an element in $E(M)$, say $b_1'$ and $b_2'$, respectively. 
Then since $f$ is on $t$ long lines of $N$ other than $\cl_N(\{b_1,b_2\})$, we have $|E(M)\cap \cl_N(\{b_i,b_3\})|= t+1$ for each $i\in \{1,2\}$.
Since $b_3$ is not parallel to an element in $E(M)$, we see that $M/b_1'$ has a $U_{2,t+2}$-restriction.
\end{proof}


We now prove a result which does the bulk of the work in the proof of Theorem \ref{critical structure}, the main result of this section.
The proof of this proposition relies on Proposition \ref{g-star} to find star-partitions.
Recall that $\cF\cap \cU(t)$ is the class of frame matroids with no $U_{2,t+2}$-minor.
Also, if $P$ is a porcupine with tip $f$, then we write $d(P)$ for the corank of $\si(P/f)$.

\begin{proposition}  \label{hard lemma}
For all integers $t\ge 2,h\ge 0$ and $g\ge 4$, if $M$ is a matroid with a size-$(2h+1)$ independent set $S$ so that $\elem(M)-\elem(M/f)>(t-1)(r(M)-1+h)+1$ for each $f\in S$, then $M$ has either
\begin{enumerate}[(1)]
\item a restriction of rank at most $3g$ which is not in $\cF\cap\cU(t)$, or

\item a size-$(h+1)$ independent set such that each element is the tip of a $g$-porcupine restriction $P$ of $M$ with $d(P)=h+1$. 
\end{enumerate} 
\end{proposition}
\begin{proof}
Assume that (1) and (2) do not hold for $M$, and that $M$ is simple.
Since (1) does not hold, $M$ has no $U_{2,t+2}$-restriction.
Also, $M$ has no spike restriction of rank at least five and at most $3g$, since Lemma \ref{transversals} and Lemma \ref{star} together imply that spikes of rank at least five are not frame matroids. 
For convenience, we say that a $g$-(pre)porcupine $P$ is an \emph{$(s,g)$-(pre)porcupine} if  $d(P)=s$.

If $f\in S$ is not the tip of a spike of rank less than $g$, then $f$ is on at least $r(M)+h$ long lines, or else 
$$\delta(M,f)=\elem(M)-\elem(M/f)\le 1+(t-1)(r(M)+h-1).$$
But if $f\in S$ is not the tip of a spike of rank less than $g$ and is on at least $r(M)+h$ long lines of $M$, then $f$ is the tip of an $(h+1,g)$-porcupine. Since (2) does not hold, there is some $S_1\subseteq S$ such that $|S_1|\ge |S|-h\ge h+1$ and each element of $S_1$ is the tip of a spike of rank less than $g$ in $M$. 
Each of these spikes has rank at most four, or else (1) holds.

For each $f\in S_1$, let $P_f=M|(\cup_{L\in \cL_M(f)}L)$, and let $(P_f/f)|T_f$ be a simplification of $P_f/f$.
Note that $\delta(P_f,f)=\delta(M,f)$. 


\begin{claim}\label{tipless}
There is some $f\in S_1$ for which $(P_f/f)|T_f$ has a star-partition $(X,\cL)$ such that each line of $P_f$ through $f$ and an element of $T_f-X$ has length three, and $|\cL|\le r(M)+h-2$. 
\end{claim}
\begin{proof}
By Proposition \ref{g-star}, for each $f\in S_1$, the matroid $(P_f/f)|T_f$ has a star-partition $(X,\cL)$ so that each circuit of $(P_f/f)|T_f$ of size less than $g$ is contained in $X\cup L\cup L'$ for some $L,L'\in\cL$, and each line of $P_f$ through $f$ and an element of $T_f-X$ has length three. 

Since $(P_f/f)|T_f$ has a circuit of size at most four, there is some $L\in \cL$ such that $|L|\ge 2$.
Let $T'\subseteq T_f$ such that $|T'|=|\cL|+1$, while $T'$ contains a transversal of $\cL$ and $|T\cap L|=2$.
If $|T'|\ge r(M)+h$, then the set of elements of $P_f$ on a long line of $P_f$ through $f$ and an element of $T'$ is an $(h+1,g)$-preporcupine with tip $f$, and thus contains an $(h+1,g)$-porcupine with tip $f$. 
Since (2) does not hold and $|S_1|-h\ge 1$ the claim holds. 
\end{proof}


Note that $(t-1)(r(M)-1+h)+1=(t-1)(r(M)+h-2)+t$, and let $f\in S_1$ and $(X,\cL)$ be given by  \ref{tipless}.
Let $m=\sum_{x\in X}(|\cl_{M}(\{x,f\})|-2)$, and note that $m\le t-1$ since $M$ has no $U_{2,t+2}$-restriction.
Then 
\setcounter{equation}{0}
\begin{align*}
(t-1)(r(M)+h-2)+t&\le \delta(P_f/f)-1\\
&=\sum_{e\in T_f}\big(|\cl_{M}(\{e,f\})|-2\big)\\
&= m+\sum_{L\in \cL}|L|\\
&\le m+(t-1)(r(M)+h-2)\\
& \hspace{0.75cm} -(t-1)|\cL|+\sum_{L\in\cL}|L|.
\end{align*}
The third line holds because each line of $P_f$ through $f$ and an element of $T_f-X$ has length three by  \ref{tipless}, and the last inequality holds since $|\cL|\le r(M)+h-2$.
Thus, 
\begin{alignLetter}
t-m+(t-1)|\cL|\le \sum_{L\in \cL}|L|
\end{alignLetter}
so there is some $L_1\in \cL$ such that $|L_1|\ge t$ since $m\le t-1$.
By Lemma \ref{m=1} each $L\in \cL-\{L_1\}$ satisfies $|L|=1$, or else (1) holds.
Then we have $\sum_{L\in\cL}|L|=|L_1|+|\cL|-1$.

We have two cases to consider. 
If $m>0$, then $X\ne\varnothing$ and thus $|L_1|=t$ or else $(P_f/f)|(L_1\cup X)$ has a $U_{2,t+2}$-restriction by the definition of a star-partition. 
Then $m=t-1$ or else (i) does not hold. 
But then by Lemma \ref{rank-3}, $M$ has a rank-3 restriction which is not in $\cU(t)$ and (1) holds. 

Now assume that $m=0$.
Since $P_f/f$ has no $U_{2,t+2}$-restriction, $|L_1|\le t+1$.
Then $t=2$ and $|L_1|=3$, or else (i) does not hold. 
But then the set of elements of $P_f$ on a line through $f$ and an element of $L_1$ is a rank-3 spike with no $U_{2,4}$-minor.
Thus, $P_f$ has an $F_7$-restriction, and (1) holds since $F_7$ is not a frame matroid.
\end{proof}
\end{subsection}


\begin{subsection}{The Proof} \label{critical proof sec}
We now prove the main result of this section, which refines the structure found in the second outcome of Theorem \ref{new reduction}. 
The case $h=0$ is particularly important and has a more exact flavor; it says that if $M$ has an element with density loss greater than any element of a Dowling geometry with group size less than $t$, then $M$ has a highly structured restriction which is not contained in any Dowling geometry with group size less than $t$.
The proof only uses Proposition \ref{skew spikes} and Proposition \ref{hard lemma}.


\begin{theorem}\label{critical structure}
Let $\ell\ge t\ge 2$ and $h\ge 0$ be integers, and let $M\in \cU(\ell)$ be a matroid with a size-$2^{15h}$ independent set $S$ so that each $f\in S$ satisfies $\elem(M)-\elem(M/f)>(t-1)(r(M)-1)+\ell^{28h}$.
Then there is a set $C\subseteq E(M)$ of rank at most $h2^{h+7}$ so that $M/C$ has either
\begin{enumerate}[(1)]
\item an $(\cF\cap \cU(t),15\cdot 2^{h},h+1)$-stack restriction,

\item an element $e$ and a collection $\cS$ of $h+2$ mutually skew sets in $M/(C\cup \{e\})$ so that for each $R\in \cS$, the matroid $(M/C)|(R\cup \{e\})$ is a spike of rank at most four with tip $e$, or 

\item a size-$(h+1)$ independent set so that each element is the tip of a $5\cdot 2^{h}$-porcupine restriction $P$ of $M/C$ with $d(P)=h+1$.

\end{enumerate} 
\end{theorem}
\begin{proof}
The constants in the theorem statement are larger than we need, to make the statement more readable. 
Using crude estimates for $h\ge 1$ we have $2^{15h}\ge 2h+h(15\cdot 2^h)(3(h+2)+1)$, and 
$$\ell^{28h}\ge (t-1)h+1+h(15\cdot 2^{h})\ell^{6(h+2)},$$
and these inequalities also hold for $h=0$. 
Let $M$ be a matroid, and assume that there is no set $C$ of rank at most $h^{h+7}$ for which (1) or (2) holds; we will show that there is such a set for which (3) holds. 
Let $k\ge 0$ be maximal so that $M$ has an $(\cF\cap \cU(t),15\cdot 2^{h},k)$-stack restriction $C_0$. 
Then $k\le h$, since (1) does not hold, so $r_M(C_0)\le h(15\cdot2^{h})$. 
The following claim uses Proposition \ref{skew spikes} to reduce to a situation in which we can apply Proposition \ref{hard lemma}.


\begin{claim}
There is some $S_1\subseteq S$ with $|S_1|\ge 2h+1$ so that $S_1$ is independent in $M/C_0$,  and each $f\in S_1$ satisfies $\delta(M/C_0,f)> (t-1)(r(M)-1+h)+1$. 
\end{claim}
\begin{proof}
Let $c_1,c_2,\dots,c_b$ be a basis of $C_0$, and let $m=h+2$. 
Let $S_1\subseteq S$ be a maximum-size independent set in $M/C_0$ for which each $f\in S_1$ satisfies $\delta(M/C_0,f)>(t-1)(r(M)-1+h)+1$, and assume for a contradiction that $|S_1|\le 2h$.
Then there is some $S'\subseteq S$ with $|S'|\ge |S|-b-2h\ge b(3m)$ so that $S'$ is independent in $M/C_0$, and each $f\in S'$ satisfies
\begin{align*}
\delta(M,f)-\delta(M/C_0,f)\ge b\ell^{6m}.
\end{align*}
 Then for each $f\in S'$ there is some $i\in [b]$ such that 
$$\delta(M/\{c_1,\dots,c_{i-1}\},f)-\delta(M/\{c_1,\dots,c_i\},f)\ge \ell^{6m},$$
since $r_M(C_0)=b$.
Since $|S'|\ge b(3m)$, there is some $i\in [b]$ which is chosen for at least $3m$ elements of $S'$. 
Let $S''\subseteq S'$ denote this set of $3m$ elements.
By Proposition \ref{skew spikes} with $(M,X,e)=(M/\{c_1,\dots,c_{i-1}\},S'',c_i)$, outcome (2) holds with $C=\{c_1,\dots,c_{i-1}\}$ and $e=c_i$, a contradiction. 
\end{proof}

By Proposition \ref{hard lemma} with $(M,S,g)=(M/C_0,S_1,5\cdot 2^{h})$ and the maximality of $k$, outcome (3) holds with $C=C_0$.
\end{proof}
\end{subsection}
\end{section}


\begin{section}{Exploiting Connectivity} \label{conn chap}
In this section we prove the following result, which contracts a collection of porcupines onto a Dowling-geometry minor.
Recall that if $P$ is a porcupine with tip $f$, then we write $d(P)$ for the corank of $\si(P/f)$.


\begin{theorem} \label{main porcupines}
There are functions $s_{\ref{main porcupines}}\colon \bZ^3\to \bZ$ and $r_{\ref{main porcupines}}\colon \bZ^6\to \bZ$, so that for all integers $\ell\ge 2$ and $k,s\ge 1$ and $g,m,n\ge 3$ and each finite group $\Gamma$,
if $M\in \cU(\ell)$ is vertically $s_{\ref{main porcupines}}(k,s,g)$-connected with no $\LG^+(n,\Gamma')$-minor with $|\Gamma'|\ge 2$, and with a $\DG(r_{\ref{main porcupines}}(\ell,m,k,n,s,g),\Gamma)$-minor $G$, and a size-$k$ independent set such that each element is the tip of a $g$-porcupine restriction $P$ of $M$ with $d(P)=s$, then $M$ has a minor $N$ of rank at least $m$ so that 
\begin{itemize}
\item $N$ has a $\DG(r(N),\Gamma)$-restriction, and

\item $N$ has a size-$k$ independent set such that each element is the tip of a $g$-porcupine restriction $P$ with $d(P)=s$.
\end{itemize}
\end{theorem}


\begin{subsection}{A Generalization of Tutte's Linking Theorem} \label{Tutte}
We first prove a generalization of Tutte's linking theorem.
The following lemma was proved in \cite{STLT lemma}, and follows from the submodularity of the connectivity function.

\begin{lemma} \label{separation}
Let $e$ be an element of a matroid $M$, and let $X,Y$ be subsets of $M-\{e\}$. Then 
$$\lambda_{M\del e}(X)+\lambda_{M/e}(Y)\ge \lambda_M(X\cap Y)+\lambda_M(X\cup Y\cup \{e\})-1.$$
\end{lemma}

We use this to prove a lemma about a collection of nested sets in a matroid.

\begin{lemma} \label{contraction GTLT}
Let $M$ be a matroid, $m\ge 1$ be an integer, and $Y_1\subseteq Y_2\subseteq \dots \subseteq Y_m\subseteq E(M)-X$. 
If $E(M)\ne \cl_M(X)\cup\cl_M(Y_m)$, then there is some $e\in E(M)-(\cl_M(X)\cup\cl_M(Y_m))$ so that $\kappa_{M/e}(X,Y_i)=\kappa_{M}(X,Y_i)$ for each $i\in [m]$.
\end{lemma}
\begin{proof}
Let $|E(M)-(\cl_M(X)\cup \cl_M(Y_m))|$ be minimal so that the claim is false, and let $e\in E(M)-(\cl_M(X)\cup \cl_M(Y_m))$.
Since the claim is false for $M$, there is some $k\in [m]$ so that $\kappa_{M/e}(X,Y_k)<\kappa_M(X,Y_k)$.
Let $(E-e-Z^k,Z^k)$ be a partition of $E(M)-\{e\}$ so that $(\cl_M(X)-Y_m)\subseteq E-Z^k$ and $Y_k\subseteq Z^k$, while $\lambda_{M/e}(Z_k)=\kappa_{M/e}(X,Y_k)$; such a partition exists because each set $A\subseteq E(M/e)$ satisfies $\lambda_{M/e}(A)\ge \lambda_{M/e}(\cl_{M/e}(A))$.
Note that $e\in \cl_M(Z^k)$, or else $\lambda_M(Z^k)=\lambda_{M/e}(Z^k)<\kappa_M(X,Y_k)$, which contradicts the defintion of $\kappa_M(X,Y_k)$.  

Now, if $|E(M)-(\cl_M(X)\cup \cl_M(Y_m))|=1$, then $Z^k\subseteq \cl_M(Y_m)$, 
since $(\cl_M(X)-Y_m)\subseteq E-Z^k$.
But then since $e\in \cl_M(Z^k)$, we have $e\in \cl_M(Y_m)$, which contradicts the choice of $e$.
Thus, $|E(M)-(\cl_M(X)\cup \cl_M(Y_m))|>1$, so by minimality there is some $j\in [m]$ so that $\kappa_{M\del e}(X,Y_j)<\kappa_M(X,Y_j)$.

First assume that $j\le k$, so $Y_j\subseteq Y_k$. Let $(E-e-Z_j,Z_j)$ be a partition of $E(M)-\{e\}$ such that $X\subseteq E-e-Z_j$ and $Y_j\subseteq Z_j$, while $\lambda_{M\del e}(Z_j)<\kappa_M(X,Y_j)$. 
Note that $X\subseteq (E-Z_j)\cap (E-Z^k)=E-(Z_j\cup Z^k\cup\{e\})$ and $Y_j\subseteq Z_j\cap Z^k$. 
Then 
\begin{align*}\setcounter{equation}{0}
\kappa_M(X,Y_j)-1+\kappa_M(X,Y_k)-1&\ge \lambda_{M\del e}(Z_j)+\lambda_{M/e}(Z^k)\\
&\ge \lambda_M(Z_j\cap Z^k)\\
&\hspace{0.5cm}+\lambda_M(Z_j\cup Z^k\cup \{e\})-1\\
&\ge \kappa_M(X,Y_j)+\kappa_M(X,Y_k)-1,
\end{align*}
a contradiction. 
The second inequality follows from Lemma \ref{separation} applied to $Z_j$ and $Z^k$, and the third holds because $Y_j\subseteq Z_j\cap Z^k$ and $X\subseteq E(M)-(Z_j\cup Z^k\cup \{e\})$. 
If $k\le j$, then we apply the same argument with $M^*$, since $\kappa_{M/e}(A,B)=\kappa_{M^*\del e}(A,B)$ for each $e\in E(M)$ and disjoint $A,B\subseteq E(M)-\{e\}$.
\end{proof}

We now prove a generalization of Tutte's linking theorem, which allows us find a minor which preserves connectivity between a set and each set in a nested family of sets.

\begin{theorem} \label{GTLT}
Let $M$ be a matroid, $m\ge 1$ be an integer, and $Y_1\subseteq Y_2\subseteq \dots \subseteq Y_m\subseteq E(M)-X$. Then $M$ has a minor $N$ with ground set $X\cup Y_m$ such that $\kappa_N(X,Y_i)=\kappa_{M}(X,Y_i)$ for each $i\in [m]$, while $N|X=M|X$ and $N|Y_1=M|Y_1$. 
\end{theorem}
\begin{proof}
Let $N$ be a minimal minor of $M$ such that $X\cup Y_m\subseteq E(N)$, while $\kappa_N(X,Y_i)=\kappa_{M}(X,Y_i)$ for each $i\in [m]$, $N|X=M|X$, and $N|Y_1=M|Y_1$. 
Assume for a contradiction that $E(N)\ne X\cup Y_m$. 
Take $e\in E(N)-(X\cup Y_m)$.
If $E(N)=\cl_N(X)\cup\cl_N(Y_1)$, then $\kappa_{N\del e}(X,Y_i)=\kappa_N(X,Y_i)$ for each $i\in [m]$ by Theorem \ref{STLT} applied with $Y=Y_1$, which contradicts the minor-minimality of $N$.

Let $j\in [m]$ be maximal so that $E(N)\ne \cl_N(X)\cup \cl_N(Y_j)$.
By Lemma \ref{contraction GTLT} applied to $X$ and $Y_1,\dots,Y_j$, there is some $e\in E(N)-(\cl_N(X)\cup \cl_N(Y_j))$ so that $\kappa_{N/e}(X,Y_i)=\kappa_N(X,Y_i)$ for each $1\le i\le j$.
There is some integer $i\in [m]$ so that $\kappa_{N/e}(X,Y_i)<\kappa_N(X,Y_i)$, or else $N/e$ contradicts the minimality of $N$.
Then $i>j$, so $E(N)=\cl_N(X)\cup \cl_N(Y_i)$ by the maximality of $j$.
Then each set $Z$ with $X\subseteq Z\subseteq E(N)-Y_i$ and $\lambda_N(Z)=\kappa_N(X,Y_i)$ satisfies $Z\subseteq \cl_N(X)$.
Since $e\notin \cl_N(X)$ this implies that $\lambda_{N/e}(Z-\{e\})=\lambda_N(Z)$ for each such set $Z$, and thus $\kappa_{N/e}(X,Y_i)=\kappa_N(X,Y_i)$, a contradiction.
Thus, $E(N)=X\cup Y$, and $N$ is the desired minor of $M$. 
\end{proof}

We also need a lemma concerning $\kappa_M$.

\begin{lemma} \label{sep union}
Let $Y$ and $J$ be disjoint sets of elements of a matroid $M$, and let
$$D=\{e\in E(M)-(J\cup Y)\colon \kappa_{M/e}(A,J)< \kappa_M(A,J) \text{ for some } A\subseteq Y\}.$$
Then $\kappa_M(Y\cup D,J)= \kappa_M(Y,J)$.
\end{lemma}
\begin{proof}
Let $E=E(M)$, and let $\kappa=\kappa_M$.
Clearly $\kappa(Y\cup D,J)\ge \kappa(Y,J)$; assume for a contradiction that $\kappa(Y\cup D,J)>\kappa(Y,J)$.
Let $D_1$ be a maximal subset of $D$ so that $\kappa(Y\cup D_1,J)= \kappa(Y,J)$, and let $e\in D-D_1$.
Since $e\in D$, there is some $A\subseteq Y$ so that $\kappa_{M/e}(A, J)<\kappa_M(A,J)$, which implies that $\kappa(A\cup\{e\},J)=\kappa(A,J)$.
Then there is some $Z_1\subseteq E$ so that $A\cup\{e\}\subseteq Z_1\subseteq E-J$ and $\lambda(Z_1)=\kappa(A,J)$.
Similarly, there is some $Z_2\subseteq E$ so that $Y\cup D_1\subseteq Z_2\subseteq E-J$ and $\lambda(Z_2)=\kappa(Y\cup D_1,J)=\kappa(Y,J)$.
Then $(Y\cup D_1\cup \{e\})\subseteq Z_1\cup Z_2$, and 
\begin{align*}
\lambda(Z_1\cup Z_2)&\le \lambda(Z_1)+\lambda(Z_2)-\lambda(Z_1\cap Z_2)\\
& = \kappa(A,J)+\kappa(Y,J)-\lambda(Z_1\cap Z_2)\\
&\le \kappa(A,J)+\kappa(Y,J)-\kappa(A,J)\\
&=\kappa(Y,J),
\end{align*}
where the third lines holds because $A\subseteq Z_1\cap Z_2$ and $J\subseteq E-(Z_1\cap Z_2)$.
Thus, $\kappa(Y\cup D_1\cup \{e\},J)= \kappa(Y,J)$, which contradicts the maximality of $D_1$.
\end{proof}
\end{subsection}


\begin{subsection}{Prickles} \label{porc props}
We now define some terminology related to porcupines, and then prove several lemmas which we will use in the proof of Theorem \ref{main porcupines}.
Recall that a $g$-porcupine is a simple matroid $P$ with an element $t$ so that each line of $P$ through $t$ has length three, and $\si(P/t)$ has no coloops and girth at least $g$.
We write $d(P)$ for the corank of $\si(P/t)$, and we take the convention that the choice of $t$ is fixed.
Note that a porcupine is a spike if and only if $d(P)=1$, and that if $d(P)=0$ then $P$ consists of only the tip $t$.
If $d(P)=0$ then $P$ is a \emph{trivial porcupine}.


For porcupines $P$ and $P'$, we say that $P'$ is a \emph{subporcupine} of $P$, and write $P'\preceq P$, if $E(P')\subseteq E(P)$ and $P'$ and $P$ have the same tip $t$.
If $E(P')$ is a proper subset of $E(P)$ then we write $P'\prec P$; this implies that $d(P')<d(P)$, since $\si(P'/t)$ has no coloops.
Note that $P'$ is a restriction of $P$, and that each line of $P'$ through $t$ is a line of $P$ through $t$.
Also, note that every porcupine $P$ has a unique trivial subporcupine.

If $P_1$ and $P_2$ are subporcupines of $P$, then $P|(E(P_1)\cup E(P_2))$ is a porcupine, which we denote by $P_1\cup P_2$.
Note that if $P_1\npreceq P_2$ and $P_2\npreceq P_1$ and neither $P_1$ nor $P_2$ is equal to $P$, then $d(P_1\cup P_2)>\max(d(P_1),d(P_2))$, since $\si(P_1/t)$ and $\si(P_2/t)$ have no coloops.
 We say that $P'$ is a \emph{retract} of $P$ if $P'$ and $P$ have the same tip, $E(P')\subseteq E(P)$, and $d(P')=d(P)$.
Whenever we work with retracts it will be the case that $P$ is a restriction of a matroid $M$ and $P'$ is a restriction of a minor $N$ of $M$.
Since $d(P')=d(P)$, one can think of $P'$ as being a copy of $P$ which we recover in the minor $N$.


This notation for porcupines extends to collections of porcupines.
In the animal kingdom, a collection of porcupines is called a prickle, and we use the same terminology here.
A \emph{prickle} is a pair ${\bf R}=(R,\cP)$ where $R$ is a matroid and $\cP$ is a collection of pairwise disjoint porcupine restrictions of $R$ such that $\cup_{P\in \cP} E(P)=E(R)$. 
If each porcupine in $\cP$ is a $g$-porcupine, then $\R$ is a \emph{$g$-prickle}.
For a matroid $M$, we say that ${\bf R}$ is a \emph{prickle of $M$} if $R$ is a restriction of $M$, and we write $E({\bf R})$ for $E(R)$.
We define $d({\bf R})=\sum_{P\in \cP} d(P)$; note that if $d({\bf R})=0$, then each porcupine in $\cP$ is simply a tip. 

We say that a prickle $(R',\cP')$ is a \emph{subprickle} of $(R,\cP)$, and write $(R',\cP')\preceq (R,\cP)$, if there is a bijection $\psi\colon \cP\to \cP'$ so that for each $P\in \cP$, the porcupine $\psi(P)$ is a  subporcupine of $P$.
If there is some $P\in \cP$ so that $\psi(P)\ne P$, then we write $(R',\cP')\prec (R,\cP)$; this implies that $d(R',\cP')<d(R,\cP)$ since $d(\psi(P))<d(P)$.
Since the porcupines in $\cP$ are pairwise disjoint, this also implies that $E(R')$ is a proper subset of $E(R)$.

It is important to note that each subprickle of $(R,\cP)$ contains the tip of each porcupine in $\cP$, which implies that the unique subprickle $(R',\cP')$ of $(R,\cP)$ with $d(R',\cP')=0$ is simply the collection of tips of porcupines in $\cP$; we say that $(R',\cP')$ is the \emph{trivial subprickle} of $(R,\cP)$.
If $(R_1,\cP_1)$ and $(R_2,\cP_2)$ are subprickles of a prickle $(R,\cP)$, then $(R_1\cup R_2,\cP_1\cup \cP_2)$ denotes the subprickle of $(R,\cP)$ such that each porcupine $P$ in $\cP_1\cup\cP_2$ is the union of the porcupines $P_1$ of $\cP_1$ and $P_2$ of $\cP_2$ with the same tip as $P$.

Finally, we need terminology for recovering a prickle after applying projections and deletions.
We say that a prickle $(R',\cP')$ is a \emph{retract} of a prickle $(R,\cP)$ if there is a bijection $\psi\colon \cP\to \cP'$ so that for each $P\in \cP$, the porcupine $\psi(P)$ is a retract of $P$.
Just as for porcupines, we will only work with retracts when $(R,\cP)$ is a prickle of a matroid $M$ and $(R',\cP')$ is a prickle of a minor $N$ of $M$.
Note that the collection of tips of porcupines in $\cP'$ is equal to the collection of tips of porcupines in $\cP$, and  that $d({\bf R}')=d({\bf R})$; these properties allow us to think of ${\bf R}'$ as a copy of ${\bf R}$ which we recover in the minor $N$.


We now provide some properties of prickles which we will need in order to prove Theorem \ref{main porcupines}.
Our first lemma deals with a maximal proper subprickle, which will be important for using inductive arguments with prickles.

\begin{lemma} \label{one smaller}
Let ${\bf A}$ be a subprickle of a prickle ${\bf R}$ so that $d({\bf A})=d({\bf R})-1$.
Then
\begin{enumerate}[(i)]
\item each ${\bf Z}\preceq {\bf R}$ with ${\bf Z}\npreceq {\bf A}$ satisfies ${\bf A}\cup {\bf Z}={\bf R}$, and

\item there is a prickle ${\bf C}\preceq {\bf R}$ so that $d({\bf C})=1$ and ${\bf A}\cup {\bf C}={\bf R}$.
\end{enumerate}
\end{lemma}
\begin{proof}
Let ${\bf A}=(A,\cA)$ and ${\bf R}=(R,\cR)$.
Since $d({\bf A})=d({\bf R})-1$, there is a unique porcupine in $P\in\cR$ so that $P\notin \cA$.
Since ${\bf A}\preceq {\bf R}$, there is some $P_1\in \cA$ so that $P_1$ is a subporcupine of $P$ and $d(P_1)=d(P)-1$.
Let $t$ denote the common tip of $P$ and $P_1$.

We first prove (i).
Let ${\bf Z}=(Z,\cP)$, and let $P_2$ denote the porcupine in $\cP$ with tip $t$.
Let $(P/t)|X$ be a simplification of $P/t$, and for each $i\in \{1,2\}$ let $(P/t)|X_i$ be a simplification of $P_i/t$ so that $X_i\subseteq X$.
Since $(P/t|X_2)$ has no coloops and $X_2$ is not contained in $X_1$, we have $r^*((P/t)|(X_1\cup X_2))>r^*((P/t)|X_1)$.
Since $(P/t)|X$ has no coloops and $r^*((P/t)|X)=r^*((P/t)|X_1)+1$, it follows that $X_1\cup X_2=X$. 
Thus, $P_1\cup P_2=P$, and so ${\bf A}\cup {\bf Z}={\bf R}$.

We now prove (ii).
Since $\si(P/t)$ has no coloops, there is a subporcupine $P'$ of $P$ so that $d(P')=1$ ($P'$ is a spike) and $P'$ is not a subporcupine of $P_1$.
Let ${\bf C}=(C,\cC)$ be the unique subprickle of ${\bf R}$ so that $P'\in \cC$ and each other porcupine $P''\in \cC$ satisfies $d(P'')=0$.
Then $d({\bf C})=1$, and ${\bf C}\npreceq {\bf A}$, and thus ${\bf A}\cup {\bf C}={\bf R}$ by (i)
\end{proof}


Our second lemma finds a specific collection of maximal proper subprickles of a prickle, again for use in inductive arguments.
It is perhaps more natural to write the lemma using set intersection instead of union; we use union because we will apply this lemma to bound $\kappa_M(J,E({\bf R})-E({\bf Y}))$ for some set $J$.

\begin{lemma} \label{subprickles}
Let $M$ be a matroid with a prickle  ${\bf R}$ such that $d({\bf R})>0$, and let ${\bf Y}\preceq {\bf R}$.
Then there is a collection $\cA$ of subprickles of ${\bf R}$ so that each ${\bf A}\in \cA$ satisfies $d({\bf A})=d({\bf R})-1$ and ${\bf Y}\preceq {\bf A}$, while $|\cA|\le d({\bf R})-d({\bf Y})$ and $E({\bf R})-E({\bf Y})\subseteq \cup_{{\bf A}\in\cA}(E({\bf R})-E({\bf A}))$.
\end{lemma}
\begin{proof}
Let $\cA$ be a minimal collection of subprickles of ${\bf R}$ such that each ${\bf A}\in \cA$ satisfies $d({\bf A})=d({\bf R})-1$ and ${\bf Y}\preceq {\bf A}$, while $E({\bf R})-E({\bf Y})\subseteq \cup_{{\bf A}\in\cA}(E({\bf R})-E({\bf A}))$.
Some choice for $\cA$ exists because for each $e\in E({\bf R})-E({\bf Y})$ there is some ${\bf A}\preceq {\bf R}$ such that $d({\bf A})=d({\bf R})-1$ and ${\bf Y}\preceq {\bf A}$ and $e\notin E({\bf A})$, using that the porcupines of ${\bf R}$ are pairwise disjoint.
For each ${\bf A}\in\cA$ there is some ${\bf C}_A\preceq {\bf R}$ such that $d({\bf C}_A)=1$ and ${\bf A}\cup {\bf C}_A={\bf R}$, by Lemma \ref{one smaller} (ii).
Since $\cA$ is minimal, for each ${\bf A}\in\cA$ there is an element of $E({\bf R})-E({\bf Y})$ which is in $E({\bf C}_A)-E((\cup_{{\bf A}'\in\cA-\{{\bf A}\}}{\bf C}_{A'}))$.
This implies that $d({\bf Y}\cup (\cup_{{\bf A}\in \cA}{\bf C}_A))\ge d({\bf Y})+|\cA|$, since $\si(P/t)$ has no coloops for any porcupine $P$ with tip $t$.
However, since ${\bf Y}\cup(\cup_{{\bf A}\in\cA}{\bf C}_A)\preceq {\bf R}$ we have $d({\bf Y}\cup (\cup_{{\bf A}\in \cA}{\bf C}_A))\le d({\bf R})$.
Thus, $d({\bf Y})+|\cA|\le d({\bf R})$, so $|\cA|\le d({\bf R})-d({\bf Y})$.
\end{proof}


Before proving our next lemma about prickles, we need a straightforward lemma about the corank of sets in a matroid.
Again, we state this lemma using union instead of intersection for ease of application.

\begin{lemma} \label{circ coll}
Let $M$ be a matroid, and let $X\subseteq E(M)$ so that $M$ and $M|X$ have no coloops.
Then there is a collection $\cY$ of subsets of $E(M)$ so that $|\cY|\le r^*(M)-r^*(M|X)$ and $E(M)-X\subseteq \cup_{Y\in\cY}(Y-X)$, while each $Y\in \cY$ satisfies $X\subseteq Y$ and $r^*(M|Y)=r^*(M|X)+1$, and $M|Y$ has no coloops.
\end{lemma}
\begin{proof}
Let $E=E(M)$.
This statement is easier to prove in $M^*$.
Since $(M|X)^*=M^*/(E-X)$, we have $r^*(M|X)=r^*(M)-r^*(E-X)$.
Since $M|X$ has no coloops, the set $E-X$ is a flat of $M^*$.
Let $B$ be a basis of $M^*|(E-X)$, and let $\cH=\{\cl_{M^*}(B-\{e\})\colon e\in B\}$.
Then each set $H\in \cH$ is a hyperplane of $M^*|(E-X)$, and $|\cH|=|B|=r^*(E-X)=r^*(M)-r^*(M|X)$.
Moreover, for each $e\in E-X$ there is some $H\in \cH$ so that $e\notin H$, which implies that $\cap \cH=\varnothing$.

Let $\cY=\{E-H\colon H\in \cH\}$, and note that $\cap \cY=X$, so $E-X\subseteq \cup_{Y\in\cY}(Y-X)$.
Fix some $Y\in\cY$, and let $H\in \cH$ so that $Y=E-H$.
Since $H$ is a hyperplane of $M^*|(E-X)$ and $E-X$ is a flat of $M^*$, it follows that $H$ is a flat of $M^*$, so $M^*/H$ has no loops.
Thus, $M|Y$ has no coloops.
Finally, since $H$ is a hyperplane of $M^*|(E-X)$ we have $r(M^*/H)=r(M^*/(E-X))+1=r^*(M|X)+1$.
Since $M^*/H=(M|(E-H))^*$, this implies that $Y$ satisfies $r(M|Y)^*=r^*(M|X)+1$.
\end{proof}


The following lemma describes projections of a prickle from the perspective of a subprickle.
It generalizes the fact that a projection of a spike by an element not parallel to the tip is the union of at most two spikes. 
Just as for Lemma \ref{subprickles}, we will apply this lemma to bound $\kappa_M(J,E({\bf A})-E({\bf A}_0))$ for some set $J$, so we write the statement using set union instead of intersection.

\begin{lemma} \label{two spikes}
Let $M$ be a matroid with prickles ${\bf A}_0$ and ${\bf A}$ such that ${\bf A}_0\preceq {\bf A}$ and $d({\bf A}_0)=d({\bf A})-1$.
Let $C\subseteq E(M)$ so that $M/C$ has a retract ${\bf A}_0'$ of ${\bf A}_0$.
Then there is a non-empty set $C'$ and a collection $\cA$ of prickles of $M/C$ such that each ${\bf A}'\in\cA$ is a retract of ${\bf A}$ with ${\bf A}_0'\preceq {\bf A}'$, while $|\cA|+|C'|\le r_M(C)+2$ and $E({\bf A})-E({\bf A}_0)\subseteq \cl_M(C\cup C' \cup (\cup_{{\bf A}'\in\cA}(E({\bf A}')-E({\bf A}_0'))))$. 
\end{lemma}
\begin{proof}
We may assume that $C$ is independent, since the lemma holds for $C$ if and only if it holds for a basis of $C$.
Since $d({\bf A}_0)=d({\bf A})-1$ and ${\bf A}_0\preceq {\bf A}$, there are porcupines $P$ of ${\bf A}$ and $P_0$ of ${\bf A}_0$ such that $P_0$ is a subporcupine of $P$ and $d(P_0)=d(P)-1$.
Then $M/C$ contains a retract $P_0'$ of $P_0$.
Let $t$ denote the common tip of $P$, $P_0$, and $P_0'$.
Note that $t\notin\cl_M(C)$ since $P_0'$ is a retract of $P_0$, and thus each long line of $P$ through $t$ contains at most one element of $\cl_M(C)$.

Let $(P/t)|X$ be a simplification of $P/t$, let $(P/t)|X_0$ be a simplification of $P_0/t$, and let $M/(C\cup \{t\})|X_0'$ be a simplification of $P_0'/t$ such that $X_0'\subseteq X_0\subseteq X$ and $\cl_M(C)\cap X=\varnothing$.
Such a set $X$ exists because each long line of $P$ through $t$ contains at most one element of $\cl_M(C)$.
Let $N=M/(C\cup\{t\})$.
Then $N|X$ has no coloops, and 
\begin{align*}
r^*(N|X)-r^*(N|X_0')&\le \big(r^*((M/t)|X)+|C|\big)-r^*(N|X_0')\\
&=|C|+d(P)-d(P_0')\\
&=|C|+1.
\end{align*} 
In particular, this implies that $N|X$ has at most $|C|+1$ nontrivial parallel classes since $N|X_0'$ is simple.

We first find a suitable simplification of $N|X$.
Let $N|X'$ be a simplification of $N|X$ so that $X_0'\subseteq X'$, and let $L_2$ denote the set of coloops of $N|X'$.
Since $N|X$ has no coloops, each element of $L_2$ is in a nontrivial parallel pair of $N|X$.
Also, $L_2\cap X_0'=\varnothing$ since $N|X_0'$ has no coloops.
Let $X''=X'-L_2$, and let $T$ be a transversal of the nontrivial parallel classes of $N|X$ so that $L_2\cap T=\varnothing$.
Note that $|T|\le |C|+1$.
Then $X-X''\subseteq \cl_M(C\cup \{t\}\cup T)$, since each element in $X-X''$ is a loop of $N$ or is parallel in $N$ to an element of $T$.
We also see that $X_0'\subseteq X''$, the matroid $N|X''$ is simple and has no coloops, and 
\begin{align*}
r^*(N|X'')&=|X''|-r(N|X'')\\
&=|X'|-|L_2|-(r(N|X')-|L_2|)\\
&=|X'|-r(N|X')\\
&=|X'|-r(N|X)\\
&\le |X|-|T|-r(N|X)=r^*(N|X)-|T|, 
\end{align*}
where the inequality holds because we delete at least one element from each parallel class of $N|X$ to obtain $N|X'$.
Then we have 
\begin{align*}
r^*(N|X'')-r^*(N|X_0')&\le \big(r^*(N|X)-|T|\big)-r^*(N|X_0')\\
&\le |C|+1-|T|.
\end{align*}

Since $N|X''$ is simple and $r^*(N|X'')-r^*(N|X_0')\le |C|+1-|T|$, by Lemma \ref{circ coll} with $M=N|X''$ and $X=X_0'$, there is a collection $\cY$ of subsets of $X''$ such that $|\cY|\le |C|+1-|T|$  and $X''-X_0'\subseteq \cup_{Y\in\cY}(Y-X_0')$, and each $Y\in\cY$ satisfies $X_0'\subseteq Y$ and $r^*(N|Y)=r^*(N|X_0')+1$, and the matroid $N|Y$ is simple and has no coloops.

For each $Y\in\cY$, define ${\bf A}_Y$ to be the prickle  of $M/C$ obtained from ${\bf A}_0'$ by replacing $P_0'$ by a porcupine $P_Y$ with tip $t$ such that $E(P_Y)\subseteq E(P)$, each line through $t$ contains an element of $Y$, and each element of $Y$ is on a line through $t$.
Such a porcupine $P_Y$ exists since $N|Y$ is simple. 
Since $X_0'\subseteq Y$ and $r^*(N|Y)=r^*(N|X_0')+1$, we have $P_0'\preceq P_Y$ and $d(P_Y)=d(P_0')+1=d(P)$.
Since $d(P_Y)=d(P)$ and $P_Y$ and $P$ have the same tip and $E(P_Y)\subseteq E(P)$ it follows that $P_Y$ is a retract of $P$, and thus ${\bf A}_Y$ is a retract of ${\bf A}$.
Since $X''-X_0'\subseteq \cup_{Y\in\cY}(Y-X_0')$ and $X-X''\subseteq \cl_M(C\cup \{t\}\cup T)$ we have 
$$E({\bf A})-E({\bf A}_0)\subseteq \cl_M(C\cup \{t\}\cup T\cup (\cup_{{Y}\in \cY}(E({\bf A}_Y)-E({\bf A}_0')))).$$
Thus, $\cA=\{{\bf A}_Y\colon Y\in\cY\}$ and $C'=T\cup \{t\}$ satisfies the lemma statement.
\end{proof}


Given a $g$-porcupine $P$ with tip $t$, it is often easier to work with a simplification of $P/t$, which is a matroid with girth at least $g$ and no coloops.
The following lemma shows that we can piece together two minors with large girth and corank to find a restriction with large girth and corank.

\begin{lemma} \label{uncontraction}
Let $M$ be a matroid with sets $Z,K_1,K_3\subseteq E(M)$ such that $Z\subseteq K_1\subseteq E(M)-K_3$, and let $s\ge s'\ge 0$ and $g\ge 3$ be integers.
If $M|Z$ has corank $s'$ and girth at least $g$ and $(M/K_1)|K_3$ has corank $s-s'$ and girth at least $g$, then $M|(K_1\cup K_3)$ has a restriction with corank $s$ and girth at least $g$.
\end{lemma}
\begin{proof}
Let $B$ be a basis of $(M/Z)|(K_1-Z)$. 
If $M|(Z\cup B\cup K_3)$ has a circuit of size less than $g$, it is not contained in $Z\cup B$ since $M|Z$ has girth at least $g$, and $B$ is independent in $M/Z$.
But then $(M/K_1)|K_3$ has a circuit of size less than $g$, a contradiction. 
Also, 
\begin{align*}
|Z\cup B\cup K_3|-r_M(Z\cup B\cup K_3)&=\big(|Z\cup B|-r_M(Z\cup B)\big)\\
&\hspace{0.5cm}+\big(|K_3|-r_{M/(Z\cup B)}(K_3)\big)\\
&=s'+(s-s')=s,
\end{align*}
so $r^*(M|(Z\cup B\cup K_3))=s$ and the lemma holds.
\end{proof}


The final lemma essentially shows that a projection of a matroid with corank $s$ and girth at least $g$ has a restriction with corank $s$ and girth at least $g/2$.

\begin{lemma} \label{projections1}
For all integers $g\ge 3$, $s\ge 1$, and $m\ge 0$, if $M$ is a matroid with $C,S\subseteq E(M)$ such that $M|S$ has corank $s$ and girth at least $g2^m$ and $\sqcap_M(S,C)\le m$, then $(M/C)|(S-C)$ has a restriction with corank $s$ and girth at least $g$.
\end{lemma}
\begin{proof}
Since $g\ge 3$ we may assume that $M$ is simple.
Let $|C|$ be minimal so that the lemma is false. 
Then $r_M(C)\ge 1$ and $C\subseteq \cl_M(S)$.
If $e\in C\cap S$ then $(M/e)|(S-\{e\})$ has corank $s$ and girth at least $g2^m-1$.
Since $g2^m-1\ge g2^{m-1}$, the result holds by the minimality of $|C|$.

If $e\in C-S$ then $(M/e)|S$ has corank $s+1$, and has at most one circuit of size less than $g2^{m-1}$.
To see this, if $(M/e)|S$ has distinct circuits $C_1$ and $C_2$ of size less than $g2^{m-1}$, then $C_1\cup\{e\}$ and $C_2\cup\{e\}$ are circuits of $M|(S\cup\{e\})$.
But then $M|((C_1\cup C_2)-\{e\})$ has a circuit of size at most $2(g2^{m-1})-1<g2^m$, which contradicts that $M|S$ has girth at least $2g^m$.
Thus, for each $f$ in the smallest circuit of $(M/e)|S$, the matroid $(M/e)|(S-\{f\})$ has corank $s$ and girth at least $g2^{m-1}$, and the result holds by the minimality of $|C|$. 
\end{proof}

\end{subsection}


\begin{subsection}{The Dowling-Geometry-Restriction Case} \label{hardest sec}
We now prove Theorem \ref{main porcupines} in the case that $G$ is a Dowling-geometry restriction of $M$.
Since vertical connectivity is not preserved under taking minors, we strengthen Theorem \ref{main porcupines} enough so that we no longer need vertical connectivity in the statement.
To accomplish this, we identify a prickle ${\bf Q}$ which contains all prickles whose connectivity to $E(G)$ is `too small', and we maintain the property that $\kappa_M(E(G),E({\bf Q}))=r_M(E({\bf Q}))$.
The lemma below makes no reference to Dowling geometries, although we only invoke this lemma when $M|J$ is a Dowling geometry.
This is our most technical result, and it uses the results from Sections \ref{Tutte} and \ref{porc props}.

\begin{lemma}\label{hard}
Let $g\ge 3$ and $k,s\ge 1$ be integers, and define integers $n_0=k+1$, and $n_{i}=2kg2^{2in_{i-1}+1}$ for $i\ge 1$. Let $M$ be a matroid with $J\subseteq E(M)$ such that $M\del J$ has a $g$-prickle  ${\bf R}=(R,\cP)$ for which
\begin{itemize}
\item $|\cP|\le k$ and each $P\in\cP$ satisfies $d(P)\le s$, and

\item there is a subprickle ${\bf Q}$ of ${\bf R}$ such that $\kappa_M(J,E({\bf Q}))=r_M(E({\bf Q}))<n_{d(\Q)}$ and each ${\bf A}\preceq {\bf R}$ with $\kappa_M(J,E({\bf A}))< n_{d({\bf A})}$ satisfies ${\bf A}\preceq {\bf Q}$.
\end{itemize}
Then $M$ has a contract-minor $M'$ which has $M|J$ as a spanning restriction, and has a $g$-retract $\R'$ of $\R$ so that $\Q\preceq \R'$.
\end{lemma}
\begin{proof}
Note that $n_i>2kin_{i-1}+kg2^{2in_{i-1}+1}$ for each $i\ge 1$.
Let $M$ with prickle ${\bf R}=(R,\cP)$ and ${\bf Q}\preceq \R$ be a counterexample with $d({\bf R})$ minimum.
If ${\bf R}={\bf Q}$, then by Theorem \ref{STLT} applied to $J$ and $E({Q})$, the result holds with $\R'={\bf R}$. 
If $d({\bf R})=0$, then $\R\preceq \Q$ and so ${\bf R}= {\bf Q}$ and the result holds.
Thus, ${\bf R}\neq {\bf Q}$ and $d({\bf R})>0$.

We will take a minimal minor $M_1$ of $M$ with a retract $\R_1$ of $\R$ so that $\R_1$ satisfies a slightly weaker property than the second condition of the lemma statement, so that we can exploit the minimality of $d(\R)$.
We first develop some notation for the collection of all such minors and retracts.
Let $\cM$ denote the collection of pairs $(M',\R')$ where $M'$ is a minor of $M$ with $M'|J=M|J$, and $\R'$ is a retract of $\R$ contained in $M'$ with $\Q\preceq \R'$ so that
\begin{enumerate}[(i)]
\item $\kappa_{M'}(J,E(\Q))=r_M(E(\Q))<n_{d(\Q)}$, and 

\item each $\A\prec \R'$ with $\kappa_{M'}(J,E(\A))<n_{d(\A)}$ satisfies $\A\preceq \Q$.
\end{enumerate} 


Note that condition (ii) only applies for proper subprickles of $\R'$.
Also, we do not require $\R'$ to be a $g$-prickle; we will show later that this follows from (i) and (ii).
The idea behind conditions (i) and (ii) is that all prickles with connectivity `too small' to $J$ are subprickles of $\Q$ by (ii), and are thus `safe' by (i).
Our first claim shows that if $(M',\R')\in \cM$, then proper subprickles of $\R'$ with `small' connectivity to $J$ in $M'$ which are not subprickles of $\Q$ have a nested structure.

\begin{claim}\label{nested1}
Let $(M',\R')\in \cM$.
If ${\bf A}_1\prec {\bf R}'$ and ${\bf A}_2\prec {\bf R}'$ such that $\kappa_{M'}(J,E({\bf A}_i))\le 2(d({\bf A}_i)+1)n_{d({\bf A}_i)}$ for each $i\in \{1,2\}$, then either ${\bf A}_1\preceq {\bf A}_2$ or ${\bf A}_2\preceq {\bf A}_1$, or $({\bf A}_1\cup {\bf A}_2)\preceq {\bf Q}$. 
\end{claim}
\begin{proof}
If the claim is false, then ${\bf A}_1\cup {\bf A}_2$ is a prickle such that $d({\bf A}_1\cup {\bf A}_2)> \max(d({\bf A}_1),d({\bf A}_2))$ and $({\bf A}_1\cup {\bf A}_2)\npreceq {\bf Q}$. 
Assume that $d({\bf A}_1)\ge d({\bf A}_2)$, without loss of generality. Then 
\begin{align*}
\kappa_{M'}(J,E({\bf A}_1\cup {\bf A}_2))&\le 2(d({\bf A}_1)+1)n_{d({\bf A}_1)}+2(d({\bf A}_2)+1)n_{d({\bf A}_2)}\\
&\le 4(d({\bf A}_1)+1)n_{d({\bf A}_1)}\\
&<n_{d({\bf A}_1\cup {\bf A}_2)},
\end{align*}
so $({\bf A}_1\cup {\bf A}_2)\preceq {\bf Q}$, a contradiction.
The last inequality holds since $4(i+1)n_i<n_{i+1}$ for all $i\ge 1$.
\end{proof}

Our second claim shows that if $(M',\R')\in\cM$, then each maximal proper subprickle of $\R'$ with `small' connectivity to $J$ in $M'$ contains each other subprickle with `small' connectivity to $J$ in $M'$.

\begin{claim} \label{Y unique}
Let $(M',\R')\in \cM$, and let $\Y'\prec \R'$ with $d(\Y')$ maximal so that $\kappa_{M'}(J,E(\Y'))\le 2(d(\Y')+1)n_{d(\Y')}$.
Then each $\A\prec \R'$ with $\kappa_{M'}(J,E(\A))\le 2(d(\A)+1)n_{d(\A)}$ satisfies $\A\preceq \Y'$.
\end{claim}
\begin{proof}
Assume for a contradiction that there is a prickle ${\bf A}\prec {\bf R}'$ so that $\kappa_{M'}(J,E({\bf A}))\le 2(d({\bf A})+1)n_{d({\bf A})}$ and ${\bf A}\npreceq {\bf Y}'$.
Then in particular ${\bf A}\ne {\bf Y}'$, so by the maximality of $d({\bf Y}')$ we have ${\bf Y}'\npreceq {\bf A}$.
By \ref{nested1} with  ${\bf A}_1={\bf A}$, and ${\bf A}_2={\bf Y}'$ we have ${\bf A}\cup {\bf Y}'\preceq {\bf Q}$.
If $\Y'\ne \Q$ then $d(\Y')<d(\Q)$, which contradicts the maximality of $d(\Y')$ since $\Q$ is a valid choice for $\Y'$ since $\kappa_{M'}(J,E(\Q))<n_{d(\Q)}\le 2(d(\Q)+1)n_{d(\Q)}$.
Thus, $\Y'=\Q$, and so $\A\preceq \Y'$ since $\A\cup \Y'\preceq \Q$, which contradicts that ${\bf A}\npreceq {\bf Y}'$.
\end{proof}


We now take a minimal minor of $M$.
Since ${\bf R}\neq {\bf Q}$ we have $\kappa_M(J,E({\bf R}))\ge n_{d({\bf R})}$, so there is some porcupine $P_k\in \cP$ such that $\kappa_M(J,E(P_k))\ge \frac{1}{k}n_{d({\bf R})}$. 
Let $M_1$ be a minimal minor of $M$ for which there is a retract $\R_1=(R_1,\cP_1)$ of $\R$ so that  $(M_1,\R_1)\in\cM$, and 
\begin{enumerate}[$(*)$]
\item each ${\bf Y}\prec {\bf R}_1$ with $d({\bf Y})$ maximal such that $\kappa_{M_1}(J,E({\bf Y}))\le 2(d({\bf Y})+1)n_{d({\bf Y})}$ satisfies $\kappa_{M_1}(J,E({\bf Y})\cup E(S_k))\ge \frac{1}{k}n_{d({\bf R}_1)}$, where $S_k\in \cP_1$ has the same tip as $P_k\in \cP$. 
\end{enumerate}
Note that $(M_1,\R_1)$ exists because $(M,\R)$ is a valid choice.
Also, note that $(*)$ is a relaxation of the connectivity requirement that $\kappa_{M_1}(J,E(\R_1))\ge n_{d(\R_1)}$.
This allows us to find a proper subprickle of $\R_1$  for which we can exploit the minimality of $d({\bf R}_1)$.


Let ${\bf Y}$ be a proper subprickle of ${\bf R}_1$ with $d({\bf Y})$ maximal such that $\kappa_{M_1}(J,E({\bf Y}))\le 2(d({\bf Y})+1)n_{d({\bf Y})}$. 
Note that ${\bf Y}$ exists because the trivial subprickle of ${\bf R}_1$ is a choice for ${\bf Y}$ since $k<2n_0$.
One can show that $\Y$ is unique, although we will not use this fact.

The main idea of this proof is that there are three types of prickles ${\bf A}\preceq {\bf R}_1$, depending on the connectivity to $J$.
If $\kappa_{M_1}(J,E({\bf A}))$ is `too small' to be useful, then ${\bf A}\preceq {\bf Q}$.
If $\kappa_{M_1}(J,E({\bf A}))$ is `in danger' of becoming too small, then ${\bf A}\preceq {\bf Y}$.
Finally, if ${\bf A}\npreceq {\bf Y}$ then $\A$ is `safe'.
The following claim makes this idea more precise.
When applied with ${\bf A}={\bf R}_1$, it shows that for each $e\in E({M_1})$ which does not `interact' with ${\bf Y}$, the matroid ${M_1}/e$ has a retract $\R_1'$ so that $(M_1/e,\R_1')\in \cM$ and ${\bf Y}\preceq \R_1'$. 
We only apply this claim with ${\bf A}={\bf R}_1$, but we state it more generally so that we can prove it using induction.

\begin{claim} \label{new R}
Let $e\in E({M_1})-(\cl_{M_1}(J)\cup\cl_{M_1}(E({\bf Y})))$ such that each ${\bf Y}'\preceq {\bf Y}$ with $\kappa_{{M_1}/e}(J,E({\bf Y}'))< n_{d({\bf Y}')}$ satisfies ${\bf Y}'\preceq {\bf Q}$. Then for each prickle ${\bf A}\preceq {\bf R}_1$ with ${\bf Y}\preceq {\bf A}$, the matroid ${M_1}/e$ has a retract ${\bf A}'$ of ${\bf A}$ so that ${\bf Y}\preceq {\bf A}'$, and each ${\bf Z}\preceq {\bf A}'$ with $\kappa_{{M_1}/e}(J,E({\bf Z}))< n_{d({\bf Z})}$ and ${\bf Z}\ne {\bf R}_1$ satisfies ${\bf Z}\preceq {\bf Q}$.
\end{claim}
\begin{proof}
Assume that the claim is false for some prickle ${\bf A}$ with $d({\bf A})$ minimum and $e\in E({M_1})$. 
Then $d({\bf A})>d({\bf Y})$ or else ${\bf A}={\bf Y}$ and the claim holds with $\A'=\Y$, since $e\notin \cl_{M_1}(E({\bf Y}))$ and each ${\bf Y}'\preceq {\bf Y}$ with $\kappa_{{M_1}/e}(J,E({\bf Y}'))< n_{d({\bf Y}')}$ satisfies ${\bf Y}'\preceq {\bf Q}$.
By Lemma \ref{subprickles} applied to ${\bf A}$ and ${\bf Y}$, there is a collection $\cA$ of subprickles of ${\bf A}$ such that each ${\bf A}_0\in \cA$ satisfies $d({\bf A}_0)=d({\bf A})-1$ and ${\bf Y}\preceq {\bf A}_0$, while $|\cA|\le d({\bf A})-d({\bf Y})$ and $E({\bf A})-E({\bf Y})\subseteq \cup_{{\bf A}_0\in\cA}(E({\bf A})-E({\bf A}_0))$.

Let ${\bf A}_0\in\cA$. 
Let $j=\min(d({\bf A}),d({\bf R}_1)-1)$; we will show that $\kappa_{M_1}(J,E({\bf A})-E({\bf A}_0))\le 2n_j$.
Since $d({\bf A}_0)<d({\bf A})$, by the minimality of $d({\bf A})$ the matroid ${M_1}/e$ has a retract $\A_0'$ of ${\bf A}_0$ so that ${\bf Y}\preceq {\bf A}'_0$, and each ${\bf Z}\preceq {\bf A}'_0$ with $\kappa_{{M_1}/e}(J,E({\bf Z}))< n_{d({\bf Z})}$ and ${\bf Z}\ne {\bf R}_1$ satisfies  ${\bf Z}\preceq {\bf Q}$.

By Lemma \ref{two spikes} with $C=\{e\}$ there is a non-empty set $C'\subseteq E(M)-\{e\}$ and a collection $\cK$ of prickles of $M_1/e$ such that $|\cK|+|C'|\le 3$, each ${\bf K}\in\cK$ is a retract of ${\bf A}$ with ${\bf A}_0'\preceq {\bf K}$,  and 
\setcounter{equation}{0}
\begin{align}
E({\bf A})-E({\bf A}_0)\subseteq \cl_{M_1}\big(\big(\cup_{{\bf K}\in \cK}(E({\bf K})-E({\bf A}_0'))\big)\cup C'\cup\{e\}\big).
\end{align}
Note that for each ${\bf K}\in \cK$ we have ${\bf Y}\preceq {\bf A}_0'\preceq {\bf K}$ and $d({\bf A}_0')=d(\A)-1=d({\bf K})-1$, by the definition of a retract.

Since the claim is false for ${\bf A}$, for each ${\bf K}\in\cK$ there is some ${\bf Z}\preceq {\bf K}$ with ${\bf Z}\npreceq {\bf Q}$ and ${\bf Z}\ne {\bf R}_1$ such that $\kappa_{M_1/e}(J,E({\bf Z}))<n_{d({\bf Z})}\le n_j$. 
Since $d({\bf A}_0')=d({\bf K})-1$ and ${\bf Z}\npreceq {\bf A}_0'$, by Lemma \ref{one smaller} (i) we have $E({\bf K})-E({\bf A}_0')\subseteq E({\bf Z})$, so $\kappa_{M_1/e}(J,E({\bf K})-E({\bf A}_0'))<n_j$. 
Thus, by (1) and the fact that $|\cK|+|C'|\le 3$, we have
\begin{align*}
\kappa_{M_1/e}(J,E({\bf A})-E({\bf A}_0))&\le \kappa_{M_1/e}(J,\cup_{{\bf K}\in \cK}(E({\bf K})-E({\bf A}_0')))+|C'|\\
&\le |\cK|(n_j-1)+|C'|\\
&\le (|\cK|+|C'|-1)(n_j-1)+1\\
&\le 2n_j-1.
\end{align*}
Thus, $\kappa_{M_1}(J,E({\bf A})-E({\bf A}_0))\le 2n_j$, as claimed.

We now use that $|\cA|\le d({\bf A})-d({\bf Y})$, and $E({\bf A})-E({\bf Y})\subseteq \cup_{{\bf A}_0\in\cA}(E({\bf A})-E({\bf A}_0))$.
We have 
\begin{align*}
\kappa_{M_1}(J,E({\bf A}))&\le \kappa_{M_1}(J,E({\bf Y}))+\kappa_{M_1}(J,E({\bf A})-E({\bf Y}))\\
&\le 2(d({\bf Y})+1)n_{d({\bf Y})}+\sum_{{\bf A}_0\in \cA}\kappa_{M_1}(J,E({\bf A})-E({\bf A}_0))\\
&\le 2(d({\bf Y})+1)n_{d({\bf Y})}+|\cA|2n_{j}\\
&\le 2(d({\bf Y})+1)n_{d({\bf Y})}+2(d({\bf A})-d({\bf Y}))n_{j}\\
&\le 2(d({\bf A})+1)n_{j},
\end{align*}
since $d(\Y)\le d(\A)-1\le j$.
If ${\bf A}\ne {\bf R}_1$, then $j=d({\bf A})$ and $\kappa_{M_1}(J,E({\bf A}))\le 2(d({\bf A})+1)n_{d({\bf A})}$, and so ${\bf A}={\bf Y}$ by \ref{nested1} and the definition of ${\bf Y}$, so ${\bf A}$ is not a counterexample.
If ${\bf A}={\bf R}_1$, then $j=d({\bf R}_1)-1$ and 
$$\kappa_{M_1}(J,E({\bf R}_1))\le 2(d({\bf R}_1)+1)n_{d({\bf R}_1)-1}<\frac{1}{k}n_{d({\bf R}_1)},$$
where the last inequality holds since $2(i+2)n_i<\frac{1}{k}n_{i+1}$ for all $i\ge 0$.
This contradicts (iii) since $E({\bf Y})\cup E(S_k)\subseteq E({\bf R}_1)$.
\end{proof}


The following claim uses Lemma \ref{contraction GTLT} to show that we can contract any suitable element in $E({M_1})-\cl_{M_1}(E(S_k))$, and contradict the minimality of ${M_1}$.

\begin{claim} \label{linking}
$E({M_1})=\cl_{M_1}(E({\bf Y})\cup E(S_k))\cup\cl_{M_1}(J)$.
\end{claim}
\begin{proof}
Assume for a contradiction that $E({M_1})-(\cl_{M_1}(E({\bf Y})\cup E(S_k))\cup\cl_{M_1}(J))\ne \varnothing$, and let 
\begin{align*}
\mathcal N=\{E({\bf Q})\}&\cup \{E({\bf A})\colon {\bf A}\prec {\bf R}_1, {\bf A}\npreceq {\bf Q} \textrm{ and } \kappa_{M_1}(J,E({\bf A}))=n_{d({\bf A})}\}\\ 
&\cup\{E({\bf Y})\cup E(S_k)\}.
\end{align*}
Other than $E(\Q)$ and $E(\Y)\cup E(S_k)$, the set $\mathcal N$ contains the ground set of each prickle $\A$ with the minimum connectivity to $J$ necessary to  satisfy (ii).
By  applying \ref{nested1} to each pair of sets in $\mathcal N$ we see that $\mathcal N$ is a nested collection of sets. 
In particular, if $\A_i\npreceq \Q$ and $E(\A_i)\in \mathcal N$ for each $i\in \{1,2\}$, then either $\A_1\preceq \A_2$ or $\A_2\preceq \A_1$. 
Then by Lemma \ref{contraction GTLT} with $X=J$ there is some $e\in E({M_1})-(\cl_{M_1}(E({\bf Y})\cup E(S_k))\cup\cl_{M_1}(J))$ so that $\kappa_{{M_1}/e}(J,Z)=\kappa_{M_1}(J,Z)$ for each $Z\in \mathcal N$. 
This means that we may apply \ref{new R} with the element $e$.

By  \ref{new R} with $\A=\R_1$, the matroid ${M_1}/e$ has a retract $\R'$ so that $(M_1/e,\R')\in \cM$ and ${\bf Y}\preceq \R'$.
Note that $({M_1}/e,\R')$  satisfies (i) since $E({\bf Q})\in\mathcal N$. 
Since ${\bf R}'$ is a retract of ${\bf R}_1$ and $e\notin\cl_{M_1}(E(S_k))$, the porcupine of ${\bf R}'$ with the same tip as $S_k$ is $S_k$ itself.
We will show that $(M_1/e,\R')$ satisfies $(*)$.
Let ${\bf Y}'\prec {\bf R}'$ be a prickle with $d({\bf Y}')$ maximal such that $\kappa_{{M_1}/e}(J,E({\bf Y}'))\le 2(d({\bf Y}')+1)n_{d({\bf Y}')}$.
By \ref{Y unique} with $M'=M_1/e$ and $\A=\Y$ we have ${\bf Y}\preceq {\bf Y}'$.
Then we have 
\begin{align*}
\kappa_{{M_1}/e}(J,E({\bf Y}')\cup E(S_k))&\ge \kappa_{{M_1}/e}(J,E({\bf Y})\cup E(S_k))\\
&=\kappa_{M_1}(J,E({\bf Y})\cup E(S_k))\\
&\ge \frac{1}{k}n_{d({\bf R}_1)},
\end{align*}
where the equality holds since $E({\bf Y})\cup E(S_k)\in \mathcal N$.
Thus, $(M_1/e,{\bf R}')$ satisfies $(*)$, so ${M_1}/e$ contradicts the minimality of ${M_1}$.
\end{proof}


Define $D_Q=\{e\in E({M_1})\colon \kappa_{{M_1}/e}(J,E({\bf Q}))<\kappa_{M_1}(J,E({\bf Q}))\}$ and
\begin{align*}
D_Y=\{e\in E({M_1})\colon \kappa_{{M_1}/e}(J,E({\bf A}))<n_{d({\bf A})} \\ \textrm{ for some } {\bf A}\preceq {\bf Y} \textrm{ with } {\bf A}\npreceq {\bf Q}\},
\end{align*}
and $D=\cl_{M_1}(E({\bf Y}))\cup D_Q\cup D_Y$.
This is essentially the set of elements $e$ so that there is no retract $\R'$ of $\R$ in $M_1/e$ so that $(M_1/e,\R')\in\cM$ and $\Y\preceq \R'$.
Then $\kappa_{M_1}(J,D)= \kappa_{M_1}(J,E({\Y}))$, by Lemma \ref{sep union} with $Y=E(\Y)$. 
Let $t$ denote the tip of $S_k$.
The next claim partitions $E(S_k)-\{t\}$ into a subset of $D$ and a subset of $\cl_{M_1}(J)$.

\begin{claim} \label{partition}
$E(S_k)-\{t\}\subseteq D\cup\cl_{M_1}(J)$. 
\end{claim}
\begin{proof}
If there is some $e\in E(S_k)-(D\cup \cl_{M_1}(J))$ with $e\ne t$, then by  \ref{new R} and the definition of $D$, the matroid ${M_1}/e$ has a retract $\R'$ so that $(M_1/e,\R')\in \cM$ and ${\bf Y}\preceq \R'$.
Note that $\kappa_{{M_1}/e}(J,E({\bf Q}))=\kappa_{M_1}(J,E({\bf Q}))$ by the definition of $D_Q$. 
Since ${\bf R}'$ is a retract of ${\bf R}_1$ and $e\in E(S_k)-\{t\}$, the porcupine $S_k'$ of ${\bf R}'$ with tip $t$ is a simplification of $S_k/e$.

We will show that $(M_1/e,\R')$ satisfies $(*)$.
Let ${\bf Y}'\prec {\bf R}'$ be a prickle with $d({\bf Y}')$ maximal such that $\kappa_{{M_1}/e}(J,E({\bf Y}'))\le 2(d({\bf Y}')+1)n_{d({\bf Y}')}$.
By \ref{Y unique} with $M'=M_1/e$ and $\A=\Y$ we have ${\bf Y}\preceq {\bf Y}'$.
Since $E(M_1)=\cl_{M_1}(E({\bf Y})\cup E(S_k))\cup\cl_{M_1}(J)$ by \ref{linking}, each set $Z$ for which $J\subseteq Z\subseteq E(M_1)-(E({\bf Y})\cup E(S_k))$ and $\lambda_{M_1}(Z)=\kappa_{M_1}(J,E({\bf Y})\cup E(S_k))$ satisfies $Z\subseteq \cl_{M_1}(J)$.
Then since $e\notin\cl_{M_1}(J)$, this implies that $\lambda_{M_1/e}(Z-\{e\})=\lambda_{M_1}(Z)$ for each $Z\subseteq E(M_1)$ with $\lambda_{M_1}(Z)=\kappa_{M_1}(J,E({\bf Y})\cup E(S_k))$, and so
$$\kappa_{{M_1}/e}\big(J,(E({\bf Y})\cup E(S_k))-\{e\}\big)=\kappa_{M_1}\big(J,E({\bf Y})\cup E(S_k)\big).$$ 
Combining this with the facts that $\Y\preceq \Y'$ and $S_k'$ is a simplification of $S_k/e$ gives
\begin{align*}
\kappa_{{M_1}/e}(J,E({\bf Y}')\cup E(S_k'))&\ge \kappa_{{M_1}/e}(J,E({\bf Y})\cup E(S_k'))\\
&=\kappa_{M_1}(J,E({\bf Y})\cup E(S_k))\\
&\ge \frac{1}{k}n_{d({\bf R}_1)},
\end{align*}
so $(M_1/e,{\bf R}')$ satisfies $(*)$. 
Thus, ${M_1}/e$ contradicts the minimality of ${M_1}$.
\end{proof}


We now use the minimality of $d({\bf R})=d({\bf R}_1)$.
Let $Y_k$ denote the porcupine of ${\bf Y}$ which has the same tip as $S_k$, and
let ${\bf R}_2$ denote the subprickle of ${\bf R}_1$ obtained from $(R_1,\cP_1)$ by replacing $S_k\in \cP_1$ with $Y_k$.
If $d(Y_k)=d(S_k)$, then $Y_k=S_k$, and so $E(S_k)\subseteq E({\bf Y})$.
But then the facts that $\Y\ne \R_1$ and $2(i+1)n_i<\frac{1}{k}n_{i+1}$ for all $i\ge 0$ together imply that $\kappa_{M_1}(J,E({\bf Y})\cup E(S_k))= \kappa_{M_1}(J,E({\bf Y}))<\frac{1}{k}n_{d({\bf R})}$, which contradicts $(*)$. 
Thus, $d(Y_k)<d(S_k)$, so $d({\bf R}_2)<d({\bf R}_1)$.
Also, ${\bf R}_2$ satisfies the conditions of the lemma statement using the same prickle ${\bf Q}\preceq {\bf R}_2$ and the matroid $M_1$. 

By the minimality of $d({\bf R}_1)$, the matroid $M_1$ with the prickle ${\bf R}_2$ is not a counterexample.
Thus, there is some $X\subseteq E({M_1})$ such that $M_1/X$ has $M_1|J$ as a spanning restriction and $\bf{Q}$ as a prickle, and has a $g$-retract ${\bf R}_3$ of ${\bf R}_2$. 
Let $Y_k'$ denote the porcupine of ${\bf R}_3$ which is a retract of $Y_k$.

Let $(S_k/t)|K$ be a simplification of $S_k/t$. 
To finish the proof, it suffices to show that $({M_1}/t/X)|K$ has a restriction with corank $d(S_k)$, girth at least $g$, and no coloops.
This is because such a restriction implies that there is a $g$-porcupine $S_k'$ of $M_1/X$ which is a retract of $S_k$, and we can take $M'=M_1/X$ and obtain the desired prickle $\R'$ from the lemma statement by replacing $Y_k'$ with $S_k'$ in the prickle ${\bf R}_3$.
Let $K_1=K\cap D$ and $K_2=K-D$, and note that $K_2\subseteq \cl_{M_1}(J)$ by  \ref{partition}. 
Also note that there is some $K_Y\subseteq K_1$ so that $(S_k/t)|K_Y$ is a simplification of $Y_k/t$, since $E(\Y)\subseteq D$.


\begin{claim}
There is some $K_3\subseteq K_2$ such that $({M_1}/t/(X\cup K_1))|K_3$ has corank $d(S_k)-d(Y_k)$ and girth at least $g$.
\end{claim}
\begin{proof}
Let $m=2(d({\bf Y})+1)n_{d({\bf Y})}$, so $\kappa_{M_1}(J,D)=\kappa_{M_1}(J,E(\Y))\le m$. We first show that $({M_1}/t/K_1)|K_2$ has corank $d(S_k)-d(Y_k)$. 
Since $E(Y_k)\subseteq D$, the matroid $({M_1}/t)|K_1$ has corank at least $d(Y_k)$.
If $({M_1}/t)|K_1$ has corank greater than $d(Y_k)$, then there is some $K_Y^+\subseteq K_1$ so that $K_Y\subseteq K_Y^+$ and $(M_1/t)|K_Y^+$ has corank $d(Y_k)+1$.
Define $Y_k^+$ to be the subporcupine of $S_k$ so that each element of $K_Y^+$ is on a line of $Y_k^+$ through $t$ and each line through $t$ contains an element of $K_Y^+$, so $d(Y^+_k)=d(Y_k)+1$.
Define $\A\preceq \R_1$ to be the prickle obtained from $\Y$ by replacing $Y_k$ with $Y^+_k$, so $d(\A)=d(\Y)+1$.
Since $E(\Y)\subseteq D$ and $K_Y^+\subseteq K_1\subseteq D$ we have $\kappa_{M_1}(J,E(\A))\le \kappa_{M_1}(J,D)\le m<\frac{1}{k}n_{d({\bf A})}$, using that $2(i+1)n_{i}<\frac{1}{k}n_{i+1}$ for all $i\ge 0$.
If ${\bf A}={\bf R}_1$ this contradicts $(*)$, so ${\bf A}\ne {\bf R}_1$.
But then ${\bf A}\preceq {\bf Y}$ by \ref{Y unique} with $(M',\R')=(M_1,\R_1)$ and $\Y'=\Y$, which contradicts that $d({\bf A})=d({\bf Y})+1$. 
Thus, $({M_1}/t)|K_1$ has corank $d(Y_k)$.
Then since $({M_1}/t)|K$ has corank $d(S_k)$ and $K$ is the disjoint union of $K_1$ and $K_2$, the matroid $({M_1}/t/K_1)|K_2$ has corank $d(S_k)-d(Y_k)$. 

We now show that $({M_1}/t/K_1)|K_2$ has girth at least $g2^{m+1}$.
If $({M_1}/t/K_1)|K_2$ has a circuit $C$ of size less than $g2^{m+1}$, then $(M_1/t)|(K_1\cup C)$ 
has a circuit $C'$ which is not contained in $(M_1/t)|K_1$, and is thus not contained in $(M_1/t)|K_Y$.
This implies that $(M_1/t)|(K_Y\cup C')$ has no coloops, and has corank greater than $d(Y_k)$, since $(M_1/t)|K_Y$ has corank $d(Y_k)$. 
Define $Y_k^+$ to be the subporcupine of $S_k$ so that each element of $K_Y\cup C'$ is on a line of $Y_k^+$ through $t$ and each line through $t$ contains an element of $K_Y\cup C'$, so $d(Y^+_k)>d(Y_k)$.
Define $\A\preceq \R_1$ to be the prickle obtained from $\Y$ by replacing $Y_k$ with $Y^+_k$, so $d(\A)>d(\Y)$.
Then $E(\A)-D\subseteq \cl_{M_1}((C'-K_1)\cup\{t\})\subseteq \cl_{M_1}(C\cup\{t\})$, so we have 
\begin{align*}
\kappa_{M_1}(J,E(\A))&\le \kappa_{M_1}(J,D)+\kappa_{M_1}(J,E(\A)-D)\\
&\le m+\kappa_{M_1}(J,C\cup\{t\})\\
&< m+g2^{m+1}\le \frac{1}{k}n_{d(\A)},
\end{align*}
since $|C|<g2^m$ and $2(i+1)n_i+g2^{2(i+1)n_i+1}\le \frac{1}{k}n_{i+1}$ for all $i\ge 0$. 
If ${\bf A}={\bf R}_1$ this contradicts $(*)$, so ${\bf A}\ne {\bf R}_1$.
But then ${\bf A}\preceq {\bf Y}$ by Lemma \ref{Y unique} with $(M',\R')=(M_1,\R_1)$ and $\Y'=\Y$, which contradicts that $d({\bf A})=d({\bf Y})+1$. 
Thus, $({M_1}/t/K_1)|K_2$ has girth at least $g2^{m+1}$.

Lastly, using $K_2\subseteq\cl_{M_1}(J)$ and $K_1\subseteq D$,
\begin{align}
\sqcap_{{M_1}/t/K_1}(K_2, X-K_1)&\le \sqcap_{{M_1}/t}(K_2,X\cup K_1)\\
&\le \sqcap_{{M_1}/t}(J,(X\cup D)-\{t\})\\
&= \sqcap_{M_1/t}(J,X)+\sqcap_{M_1/t/X}(J,D-(\{t\}\cup X))\\
&\le 1+\kappa_{M_1/t/X}(J,D-(\{t\}\cup X))\\
&\le 1+\kappa_{M_1}(J,D)\le 1+m.
\end{align}
where $\sqcap_{M_1/t}(J,X)\le 1$ since $J$ and $X$ are skew in $M_1$.
Since $({M_1}/t/K_1)|K_2$ has corank $d(S_k)-d(Y_k)$ and girth at least $g2^{m+1}$ and $\sqcap_{{M_1}/t/K_1}(X-K_1,K_2)\le m+1$, the claim holds by Lemma \ref{projections1} applied with ${M}={M_1}/t/K_1$ and $C=X-K_1$ and $S=K_2$. 
\end{proof}


Since $\cl_{M_1}(E({\bf Y}))\subseteq D$, there is some $Z\subseteq K_1$ such that $({M_1}/t/X)|Z$ is a simplification of $Y'_k/t$.
Then $({M_1}/t/X)|Z$ has corank $d(Y_k)$ and girth at least $g$, so
by Lemma \ref{uncontraction} applied to ${M_1}/t/X$, $Z$, $K_1-X$, and $K_3$, there is a restriction of $({M_1}/t/X)|K$ with corank $d(S_k)$ and girth at least $g$. 
Then $({M_1}/X)|E(S_k)$ has a $g$-porcupine $S_k'$ with tip $t$ so that $d(S'_k)=d(S_k)$. 
Thus, the prickle obtained from ${\bf R}_3$ by replacing $Y'_k$ with $S_k'$ is a retract of ${\bf R}_1$, and thus a retract of ${\bf R}$, and is a $g$-prickle.
This contradicts that ${\bf R}$ is a counterexample.  
\end{proof}

\end{subsection}


\begin{subsection}{Tangles} \label{tangle sec}
To prove Theorem \ref{main porcupines}  when the Dowling geometry is not a restriction, we use objects called `tangles' to maintain connectivity between the prickle and the ground set of the Dowling geometry as we take a minor.
We now introduce tangles, and define the tangle which we will use for the remainder of this section.
Tangles were first defined for graphs by Robertson and Seymour as part of the Graph Minors Project to describe the structure of graphs with no $K_t$-minor \cite{graph tangles}.
They were generalized explicitly to matroids in \cite{tangle matroids Dharma} and \cite{tangle matroids GGRW}.
Most of the following material can be found in \cite{tangles Nelson}, which in turn is based on \cite{tangles Geelen}.

Roughly speaking, a tangle on a matroid $M$ is the collection of the `less interesting' sides of the small separations of $M$.
Let $M$ be a matroid and let $\theta\ge 2$ be an integer.
We say that a set $Z\subseteq E(M)$ is \emph{$(\theta-1)$-separating in $M$} if $\lambda_M(Z)<\theta-1$.
A collection $\pT$ of subsets of $E(M)$ is a \emph{tangle of order $\theta$} if
\begin{enumerate}[(1)]
\item each set in $\pT$ is $(\theta-1)$-separating in $M$, and for each $(\theta-1)$-separating set $Z$, either $Z\in\pT$ or $E(M)-Z\in\pT$;

\item if $A,B,C\in\cT$ then $A\cup B\cup C\ne E(M)$; and

\item $E(M)-\{e\}\notin\pT$ for each $e\in E(M)$.
\end{enumerate}

The idea is that if $\lambda_M(Z)<\theta -1$, then either $Z$ or $E-Z$ does not contain much information about $M$, and that is the set which goes into the tangle.
If $Z\in\pT$ then we say that $Z$ is \emph{$\pT$-small}. 
Intuitively, if $Z\in \pT$ then $\cl_M(Z)\in \pT$, and this is indeed true.

\begin{lemma} \label{Z closure}
If $\pT$ is a tangle of order at least three on a matroid $M$, and $Z\in \pT$, then $\cl_M(Z)\in \pT$.
\end{lemma}
\begin{proof}
If the lemma is false, then there is some $Z\in \pT$ and some $e\in \cl_M(Z)-Z$ so that $Z\cup\{e\}\notin \pT$.
Since $\lambda_M(Z\cup\{e\})\le \lambda_M(Z)$, tangle axiom (1) implies that $E-(Z\cup\{e\})\in \pT$.
Since the order of $\pT$ is at least three we have $\{e\}\in \pT$, since $\lambda_M(\{e\})\le 1$ and $E-\{e\}\notin \pT$ by axiom (3).
But then $Z\cup (E-(Z\cup\{e\}))\cup\{e\}=E(M)$, which violates axiom (2).
\end{proof}

A key property of tangles is that they induce another matroid with ground set $E(M)$.
Given a tangle $\pT$ of order $\theta$ on a matroid $M$ and $X\subseteq E(M)$, define $r_{\pT}(X)=\theta-1$ if there is no set $Z\in \pT$ for which $X\subseteq Z$, and define $r_{\pT}(X)=\min \{\lambda_M(Z)\colon X\subseteq Z\in\pT\}$ otherwise.
The following appears as Lemma 4.3 in \cite{tangle matroids GGRW}.

\begin{lemma}
If $\pT$ is a tangle of order $\theta$ on a matroid $M$, then $r_{\pT}$ is the rank function of a rank-$(\theta-1)$ matroid on $E(M)$.
\end{lemma}


This matroid is a \emph{tangle matroid} of $M$, and is denoted $M(\pT)$.
We denote the closure function $\cl_{M(\pT)}$ by $\cl_{\pT}$ for readability.
Note that Lemma \ref{Z closure} implies that $r_{\pT}(X)=r_{\pT}(\cl_M(X))$  for each $X\subseteq E(M)$, and so $\cl_M(X)\subseteq \cl_{\pT}(X)$.
This also shows that $r_{\pT}(X)\le r_M(X)$, or else there is some $X'\subseteq X$ and $e\in X-X'$ so that $e\in \cl_{M}(X')-\cl_{\pT}(X')$.

The next lemma shows the connection between tangles and minors of $M$.

\begin{lemma}
If $N$ is a minor of a matroid $M$ and $\pT_N$ is a tangle of order $\theta$ on $N$, then $\{X\subseteq E(M)\colon \lambda_M(X)<\theta-1, X\cap E(N)\in \pT_N\}$ is a tangle of order $\theta$ on $M$.
\end{lemma}

This is the tangle on $M$ \emph{induced} by $\pT_N$. 
Here the minor $N$ will always be a Dowling geometry, which allows us to use a special tangle which behaves very nicely with Dowling geometries, and in fact with round matroids in general.


For each matroid $M$ and each integer $3\le k\le r(M)$, let $\pT_k(M)$ denote the collection of $(k-1)$-separating sets of $M$ which are neither spanning nor cospanning in $M$.

We say that $M$ is \emph{round} if $E(M)$ is not the union of two hyperplanes.
This implies that if $r(X)<r(M)$, then $E(M)-X$ is spanning in $M$.
It is not hard to see that every matroid with a spanning round restriction is itself round;
in particular, every matroid with a spanning clique restriction is round.
When $M$ is round, $\pT_k(M)$ is nearly a tangle.

\begin{lemma} \label{round tangle}
If $M$ is a round matroid and $3\le k\le r(M)$, then $\pT_k(M)$ is the collection of subsets of $E(M)$ of rank at most $k-2$. Moreover, $\pT_k(M)$ satisfies tangle axioms (1) and (3).
\end{lemma}
\begin{proof}
Since $M$ is round, each nonspanning $(k-1)$-separating set has rank at most $k-2$, so each set in $\pT_k(M)$ has rank at most $k-2$.
If $r_M(X)\le k-2$ and $X$ is cospanning, then some subset of $X$ is the complement of a basis of $M$ and has rank at most $k-2$. 
But then $E(M)$ is the union of two hyperplanes, which contradicts that $M$ is round.
Thus, each set of rank at most $k-2$ is in $\pT_k(M)$.

Clearly each set in $\pT_k(M)$ is $(k-1)$-separating in $M$, and $\pT_k(M)$ satisfies tangle axiom (3).
Since $M$ is round, if $\lambda_M(X)<k-1$, then either $r_M(X)\le k-2$ or $r_M(E(M)-X)\le k-2$.
Thus, $\pT_k(M)$ satisfies tangle axiom (1).
\end{proof}


This shows that when $M$ is round, $\pT_k(M)$ is a tangle if and only if $E(M)$ is not the union of three subsets of rank at most $k-2$.
In the case that $M$ is a rank-$n$ Dowling geometry, it is not hard to see that $E(M)$ is the union of three subsets of rank $m$ if and only if $m\ge \lceil 2n/3\rceil$.

\begin{lemma} \label{Dowling tangle}
Let $n\ge 3$ be an integer, let $M\cong \DG(n,\Gamma)$, and let $3\le k\le r(M)$.
Then $\pT_{k}(M)$ is a tangle of order $k$ in $M$ if and only if $3\le k\le \lceil 2n/3\rceil+1$.
\end{lemma}

If $M$ is a matroid with a minor $G$ so that $\pT_k(G)$ is a tangle, then we write $\pT_{k}(M,G)$ for the tangle of order $k$ in $M$ induced by $\pT_k(G)$.
We will work with this tangle for the remainder of this section, most often in the case that $G$ is a Dowling geometry.
\end{subsection}


\begin{subsection}{Tangles and Dowling Geometries} \label{tangle props}
We now prove some properties of the tangle $\pT_k(M,G)$ in the case that $G$ is round, although we will only apply these lemmas in the case that $G$ is a Dowling geometry.
We then use these lemmas to prove Theorem \ref{Dowling minor}, a key ingredient in the proof of Theorem \ref{main corollary}.

We first prove two lemmas which provide lower bounds on $r_{\pT_k(M,G)}(X)$ for any $X\subseteq E(M)$.

\begin{lemma} \label{intersect}
Let $M$ be a matroid with a round minor $G$, and let $k$ be an integer so that $\pT_k(G)$ is a tangle.
Then each $X\subseteq E(M)$ satisfies $r_{\pT_k(M,G)}(X)\ge \min(r_G(X\cap E(G)), k-1)$.
\end{lemma}
\begin{proof}
Let $\pT=\pT_k(M,G)$, and $m=r_{\pT}(X)$.
If $m<k-1$ then there is some $Z\in\pT$ so that $X\subseteq Z$ and $\lambda_M(Z)=m$.
Since $G$ is a minor of $M$, $\lambda_G(Z\cap E(G))\le \lambda_M(Z)\le m$.
Since $G$ is round, either $r_G(Z\cap E(G))\le m$ or $r_G((E(M)-Z)\cap E(G))\le m<k-1$. 
In the latter case, $(E(M)-Z)\cap E(G)\in \pT_k(G)$ by Lemma \ref{round tangle}.
But then $E(M)-Z\in \pT$ by the definition of $\pT$, which contradicts that $\pT$ is a tangle and $Z\in \pT$.
Thus, $r_G(Z\cap E(G)) \le m$, so $r_G(X\cap E(G))\le m$ since $X\subseteq Z$.
\end{proof}

The second lemma shows that tangles behave nicely with vertical connectivity; the proof is similar to the proof of Lemma \ref{intersect}.

\begin{lemma} \label{tangle conn}
Let $M$ be a vertically $s$-connected matroid with a round minor $G$, and let $k$ be an integer so that $\pT_k(G)$ is a tangle.
If $X\subseteq E(M)$ and $r_{\pT_k(M,G)}(X)<\min(k-1,s-1)$, then $r_{\pT_k(M,G)}(X)\ge r_M(X)$.
\end{lemma}
\begin{proof}
Let $\pT=\pT_k(M,G)$, and $m=r_{\pT}(X)$.
Since $m<k-1$ there is some $Z\in \pT$ so that $X\subseteq Z$ and $\lambda_M(Z)=m$.
Since $M$ is vertically $s$-connected and $m<s-1$, either $r_M(Z)<m$ or $r_M(E(M)-Z)<m$.
If $r_M(E(M)-Z)<m$, then $r_G((E(M)-Z)\cap E(G))\le k-2$ since $m<k-1$.
Then $(E(M)-Z)\cap E(G)\in \pT_k(G)$ by Lemma \ref{round tangle}, so $E(M)-Z\in \pT$ by the definition of $\pT$.
This contradicts that $Z\in \pT$ and $\pT$ is a tangle.
Thus, $r_M(X)\le r_M(Z)<m$, so the lemma holds.
\end{proof}

The following lemma shows the relationship between tangles and the connectivity between a pair of sets.
We apply this lemma in the case that $G$ is a restriction of $M$.


\begin{lemma}\label{tangle restriction}
Let $M$ be a matroid with a round minor $G$, and let $k$ be an integer so that $\pT_k(G)$ is a tangle. 
If $J\subseteq E(G)$ and $X\subseteq E(M)-J$, then $\kappa_M(J,X)\ge \min(r_{\pT_k(M,G)}(X),r_G(J))$. 
\end{lemma}
\begin{proof}
Let $\pT=\pT_k(M,G)$, and $E=E(M)$.
Assume for a contradiction that $\kappa_{M}(J,X)<\min(r_{\pT}(X),r_G(J))$ for some $J\subseteq E(G)$ and $X\subseteq E-J$. Let $(Z,E-Z)$ be a partition of $E$ such that $X\subseteq Z$ and $J\subseteq E-Z$, and $\lambda_{M}(Z)<\min(r_{\pT}(X),r_G(J))$. 
Either $Z$ or $E-Z$ is in $\pT$, by tangle axiom (1). If $Z\in\pT$, then 
$$r_{\pT}(X)\le \lambda_M(Z)<\min(r_{\pT}(X),r_G(J))\le r_{\pT}(X),$$ a contradiction. 
If $E-Z\in \pT$, then $E(G)\cap (E-Z)\in \pT_k(G)$ by the definition of $\pT_k(M,G)$.
Then 
\begin{align}\setcounter{equation}{0}
r_G(J)&\le r_G((E-Z)\cap E(G))\\
&= \lambda_G((E-Z)\cap E(G))\\
&\le \lambda_M(E-Z)\\
&=\lambda_M(Z)<r_G(J),
\end{align}
a contradiction.
Line (1) holds because $J\subseteq (E-Z)\cap E(G)$ and line (2) holds by Lemma \ref{round tangle} and the fact that $(E-Z)\cap E(G)\in \pT_k(G)$.
\end{proof}

Using tangles while taking a minor is tricky, because each time we contract an element we have a new tangle, and we must ensure that the rank of each set of interest does not decrease with respect to this new tangle.
To deal with this, we prove two lemmas which provide sufficient conditions for maintaining tangle connectivity of a set as we contract towards a minor. 

In the first case we contract an element which preserves the minor.
Note that there are two different tangles in the statement of the lemma.

\begin{lemma} \label{tangle closure}
Let $M$ be a matroid with a round minor $G$, and let $k$ be an integer so that $\pT_k(G)$ is a tangle. 
Let $e\in E(M)$ so that $G$ is a minor of $M/e$, and let $X\subseteq E(M)$.
Then $r_{\pT_{k}(M/e,G)}(X-\{e\})\ge r_{\pT_k(M,G)}(X)-1$.
Moreover, if $e\notin \cl_{\pT_k(M,G)}(X)$, then $r_{\pT_k(M/e,G)}(X)=r_{\pT_k(M,G)}(X)$ and $(M/e)|X=M|X$. 
\end{lemma}
\begin{proof}
Let $\pT=\pT_k(M,G)$.
If $r_{\pT_{k}(M/e,G)}(X-\{e\})<r_{\pT}(X)-1$, then there is some $Z\in \pT_{k}(M/e,G)$ so that $X-\{e\}\subseteq Z$ and $\lambda_{M/e}(Z)<r_{\pT}(X)-1$. Then $\lambda_M(Z\cup\{e\})<r(M(\pT))$, so either $Z\cup\{e\}\in \pT$ or $E(M)-(Z\cup\{e\})\in \pT$ by tangle axiom (1).

If $E(M)-(Z\cup\{e\})\in \pT$, then $E(G)\cap (E(M)-(Z\cup\{e\}))$ is not spanning in $G$, by the definition of $\pT_k(G)$. This means that $E(G)\cap (Z\cup\{e\})$ spans $G$, since $G$ is round. But then $E(G)\cap Z$ spans $G$ since $e\notin E(G)$, which contradicts that $Z\in \pT_k(M/e,G)$, by the definition of $\pT_k(M/e,G)$.
Thus, $Z\cup \{e\}\in \pT$. But $\lambda_M(Z\cup\{e\})< r_{\pT}(X)$ and $X\subseteq Z$, a contradiction. 
If $e\notin\cl_{\pT}(X)$, then Lemma \ref{Z closure} implies that $e\notin\cl_M(X)$, so $(M/e)|X=M|X$.
\end{proof}

In the second case we contract a subset of the minor itself.
Again, the lemma statement involves two different tangles.
The proof is very similar to the proof of Lemma \ref{tangle closure}.

\begin{lemma} \label{tangle minor}
Let $M$ be a matroid with a round minor $G$ and $C\subseteq E(G)$, and let $k$ be an integer so that $\pT_k(G/C)$ is a tangle.  
Then each $X\subseteq E(M)-C$ satisfies $r_{\pT_k(M/C,G/C)}(X-C)\ge r_{\pT_k(M,G)}(X)-|C|$.
Moreover, if $r_{\pT_k(M,G)}(X\cup C)=r_{\pT_k(M,G)}(X)+|C|$, then $r_{\pT_k(M/C,G/C)}(X)=r_{\pT_k(M,G)}(X)$ and $(M/C)|X=M|X$. 
\end{lemma}
\begin{proof}
Let $\pT=\pT_k(M,G)$.
If $r_{\pT_k(M/C,G/C)}(X)< r_{\pT_k(M,G)}(X)-|C|$, then there is some $Z\in \pT_{k}(M/C,G/C)$ such that $X-C\subseteq Z$ and $\lambda_{M/C}(Z)<r_{\pT}(X)-|C|$. 
Then $\lambda_M(Z\cup C)<r_{\pT}(X)$, so either $(Z\cup C)\in \pT$ or $E(M)-(Z\cup C)\in \pT$ by tangle axiom (1).

If $E(M)-(Z\cup C)\in \pT$, then $E(G)\cap (E(M)-(Z\cup C))$ is not spanning in $G$, by the definition of $\pT_k(G)$. This means that $E(G)\cap (Z\cup C)$ spans $G$, since $G$ is round. But then $E(G)\cap Z$ spans $G/C$, which contradicts that $Z\in \pT_k(M/C,G/C)$, by the definition of $\pT_k(M/C,G/C)$.
Thus, $(Z\cup C)\in \pT$.
But $X\subseteq Z\cup C$ and $\lambda_M(Z\cup C)<r_{\pT}(X)$, a contradiction. 

If $r_{\pT}(X\cup C)=r_{\pT}(X)+|C|$, then since $\cl_M(X\cup C')\subseteq \cl_{\pT}(X\cup C')$ for all $C'\subseteq C$ by Lemma \ref{Z closure}, we have 
$$r_M(X\cup C)-r_M(X)\ge r_{\pT}(X\cup C)-r_{\pT}(X)=|C|,$$
which implies that $(M/C)|X=M|X$.
\end{proof}

Tangles can also tell us how close a minor is to being a restriction.

\begin{lemma} \label{C_0}
Let $\ell\ge 2$ and $m\ge 0$ be integers, and let $M\in \cU(\ell)$ so that $M/C$ has a  simple round restriction $G$.
Let $k$ be an integer so that $\pT_k(G)$ is a tangle.
If $r_{\pT_k(M,G)}(C)\le m$, then $r_M(E(G))-r(G)\le \ell^{m+1}$. 
\end{lemma}
\begin{proof}
Since $r_{\pT_k(M,G)}(C)\le m$, there is a set $Z\in \pT_k(M,G)$ so that $C\subseteq Z$ and $\lambda_M(Z)\le m$.
Moreover, $r_G(Z\cap E(G))\le m$ by Lemma \ref{intersect}.
Thus, $r_M(Z\cap E(G))\le |Z\cap E(G)|\le \ell^m$ since $G$ is simple and $G\in\cU(\ell)$, by Theorem \ref{l-Kung}.
We have
\begin{align*}
r_M(E(G))-r(G)&=\sqcap_M(E(G),C)\\
&\le \kappa_M(E(G),C)\\
&\le \lambda_M(Z-E(G))\\
&\le \lambda_M(Z)+|Z\cap E(G)|=m+\ell^m\le \ell^{m+1},
\end{align*}
as desired.
\end{proof}


Finally, we prove a result which lets us move from a Dowling-geometry minor to a Dowling-geometry restriction in certain situations, while maintaining the connectivity of a set $X$ to the ground set of the Dowling geometry.
This theorem relies on Lemma \ref{Ramsey}.

\begin{theorem} \label{Dowling minor}
There is a function $r_{\ref{Dowling minor}}\colon \bZ^3\to \bZ$, so that
for all integers $\ell,s\ge 2$ and $n,m\ge 3$ and each finite group $\Gamma$, if $M\in \cU(\ell)$ is vertically $(s+2)$-connected with no $\LG^+(n,\Gamma')$-minor with $|\Gamma'|\ge 2$, and with a $\DG(r_{\ref{Dowling minor}}(\ell,m,n),\Gamma)$-minor $G$, then for each $X\subseteq E(M)$ for which $r_M(X)\le \min(s,m)$, there is minor $N$ of $M$ with a spanning $\DG(m,\Gamma)$-restriction so that $N|X=M|X$.
\end{theorem}
\begin{proof}
Define $r_{\ref{Dowling minor}}(\ell,m,n)=r=\max(3m+f_{\ref{Ramsey}}(\ell,n,\ell^{m+1}), 3(m+2))$.
Let $k=m+2$, and note that $\pT_k(G)$ is a tangle by Lemma \ref{Dowling tangle} since $k\le r/3\le \lceil 2r/3\rceil$.
Let $M_1$ be a minimal minor of $M$ so that $G$ is a minor of $M_1$, while $M_1|X=M|X$ and $r_{\pT_k(M_1,G)}(X)=r_M(X)$.
Since $M$ is vertically $(s+2)$-connected and $r_M(X)\le \min(s,k-2)$, the matroid $M$ is a valid choice for $M_1$ by Lemma \ref{tangle conn}.
Let $\pT=\pT_k(M_1,G)$.

Let $C_0\subseteq E(M_1)$ so that $G$ is a restriction of $M_1/C_0$.
Then $C_0\subseteq \cl_{\pT}(X)$ by the minimality of $M_1$ and Lemma \ref{tangle closure}, so $r_{\pT}(C_0)\le r_{\pT}(X)\le r_M(X)\le m$.
Then by Lemma \ref{C_0} we have $r_{M_1}(E(G))-r(G)\le \ell^{m+1}$.
By Lemma \ref{Ramsey} applied to $M_1$ with $d=\ell^{m+1}$, the matroid $M_1|E(G)$ has a  $\DG^-(3m,\Gamma)$-restriction $G_1$.
Then $G_1$ has a $\DG^-(m+2,\Gamma)$-restriction $G_2$ so that $X\cap E(G_2)=\varnothing$. By Lemma \ref{tangle restriction} we have 
$$\kappa_{M_1}(E(G_2),X)\ge \min(r_{\pT}(X),r_{G}(E(G_2)))=r_{\pT}(X)=r_M(X).$$
By Theorem \ref{STLT}, the matroid $M_1$ has a minor $N$ so that $N|E(G_2)=M_1|E(G_2)$ and $N|X=M_1|X$, while $E(N)=E(G_2)\cup X$ and $\lambda_{N}(E(G_2))=\kappa_{M_1}(E(G_2),X)=r_M(X)$.
Since $\lambda_N(X)=r_M(X)=r_N(X)$ and $E(N)=X\cup E(G_2)$, we have $r_N(X)+r_N(E(G_2))-r(N)=r_N(X)$, so $E(G_2)$ spans $N$.
By Lemma \ref{jointless} and the fact that $r_M(X)\le m$, there is a set $C$ of $N$ of size at most two so that $N/C$ has a $\DG(r(N/C),\Gamma)$-restriction, and $(N/C)|X=N|X$.
\end{proof}

When we apply this theorem, $X$ will either be a stack, or a collection of nearly skew spikes of rank at most four.
\end{subsection}


\begin{subsection}{The Proof} \label{porc proof}
We now prove a theorem which easily implies Theorem \ref{main porcupines}, the main result of this section.
Recall that if $P$ is a porcupine with tip $f$, then we write $d(P)$ for the corank of $\si(P/f)$.

\begin{theorem} \label{new main porcupines}
There are functions $s_{\ref{new main porcupines}}\colon \bZ^3\to \bZ$ and $r_{\ref{new main porcupines}}\colon \bZ^6\to \bZ$, so that for all integers $\ell\ge 2$ and $k,s\ge 1$ and $g,m,n\ge 3$ and each finite group $\Gamma$, if $M\in\cU(\ell)$ is vertically $s_{\ref{new main porcupines}}(k,s,g)$-connected with no $\LG^+(n,\Gamma')$-minor with $|\Gamma'|\ge 2$, 
and with a $\DG(r_{\ref{new main porcupines}}(\ell,m,k,n,s,g),\Gamma)$-minor $G$, and
 sets $T,T'$ with $T'\subseteq T$ so that $|T|\le k$ and each element of $T'$ is the tip of a $g$-porcupine restriction $P$ of $M$ with $d(P)=s$, then $M$ has a minor $N$ of rank at least $m$ so that 
\begin{itemize}
\item $N$ has a $\DG(r(N),\Gamma)$-restriction, 
\item $N|T=M|T$, and
\item each element of $T'$ is the tip of a $g$-porcupine restriction $P$ of $N$ with $d(P)=s$.
\end{itemize}
\end{theorem}
\begin{proof}
Define integers $n_0=k+1$, and $n_{i}=2kg2^{2in_{i-1}+1}$ for $i\ge 1$, which appear in the statement of Lemma \ref{hard}.
Let $m_1=\max(3(ks+1)n_{ks},m)$.
Define $s_{\ref{new main porcupines}}(k,g,s)=n_{ks}+1$ and $r_{\ref{new main porcupines}}(\ell,m,k,n,s,g)=r=\max(m_1+3+f_{\ref{Ramsey}}(\ell,n,\ell^{m_1+1}), 3m_1)$.
Let $M\in\cU(\ell)$ be a vertically $s_{\ref{new main porcupines}}(k,s,g)$-connected matroid with a $\DG(r_{\ref{new main porcupines}}(\ell,m,k,n,s,g),\Gamma)$-minor $G$, a set $T\subseteq E(M)$ with $|T|\le k$, and $T'\subseteq T$ so that each element of $T'$ is the tip of a $g$-porcupine restriction $P$ of $M$ with $d(P)=s$.
Let $t=|\Gamma|$.
We may assume that $T$ contains no loop of $M$.

Let $M_0$ be obtained from $M$ by performing parallel extensions so that these porcupine restrictions are pairwise disjoint, their union is disjoint from $E(G)$, and $T$ is disjoint from $E(G)$.
Say that $E(M_0)=E(M)\cup X$, where each element of $X$ is parallel to an element of $M_0|(E(M))$.
Note that $M_0$ is vertically $s_{\ref{new main porcupines}}(k,s,g)$-connected.
Let ${\bf R}_0=(R_0,\cP_0)$ be a $g$-prickle  of $M_0$ so that
\begin{itemize}
\item the tip of each porcupine of $\cP_0$ is in $T$,

\item each element of $T$ is the tip of a porcupine in $\cP_0$,

\item  each $P\in\cP_0$ with tip in $T'$ satisfies $d(P)=s$, and 

\item each $P\in \cP_0$ with tip in $T-T'$ satisfies $d(P)=0$ (so $P$ is simply a tip).
\end{itemize} 
We first show that there exists a subprickle ${\bf Q}$ of $\R_0$ with `very small' rank which contains all other subprickles with `very small' rank.

\begin{claim} \label{first Q claim}
There is a subprickle ${\bf Q}$ of ${\bf R}_0$ so that $r_{M_0}(E({\bf Q}))<n_{d({\bf Q})}$, and each ${\bf A}\preceq {\bf R}_0$ with $r_{M_0}(E({\bf A}))< n_{d({\bf A})}$ satisfies ${\bf A}\preceq {\bf Q}$.
\end{claim}
\begin{proof}
Let ${\bf Q}$ be a subprickle of ${\bf R}_0$ with $d({\bf Q})$ maximal such that $r_{M_0}(E({\bf Q}))< n_{d({\bf Q})}$.
The trivial subprickle of ${\bf R}_0$ is a candidate for ${\bf Q}$ since $k<n_0$, so ${\bf Q}$ exists.
If there is some ${\bf A}\npreceq {\bf Q}$ such that $r_{M_0}(E({\bf A}))<n_{d({\bf A})}$, then $({\bf Q}\cup {\bf A})\preceq {\bf R}_0$ and $d({\bf Q}\cup {\bf A})>d({\bf Q})$ and $r_{M_0}(E({\bf Q}\cup {\bf A}))<n_{d({\bf Q})}+n_{d({\bf A})}<n_{d({\bf Q}\cup {\bf A})}$.
We use that $2n_i<n_{i+1}$ for all $i\ge 1$.
This contradicts the maximality of $d({\bf Q})$.
\end{proof}


The main idea of this proof is that we will use Lemmas \ref{Ramsey} and \ref{jointless} to find a $\DG(m,\Gamma)$-restriction, and then apply Lemma \ref{hard}.
In order to apply Lemma \ref{Ramsey}, we need to find a minor of $M_0$ for which we can contract a set of bounded size and obtain $G$ as a restriction.
We also need to preserve the tangle rank of prickles so that we can apply Lemma \ref{hard}.

Let $M_1$ be a minimal minor of $M_0$ such that $G$ is a minor of $M_1$ and $M_1|T=M_0|T$, while $M_1$ has a retract  ${\bf R}_1$ of $\R$ with ${\bf Q}\preceq {\bf R}_1$ so that
\begin{enumerate}[(i)]

\item $r_{\pT_{\lceil2r/3\rceil-2}(M_1,G)}(E({\bf Q}))=r_{M_0}(E({\bf Q}))$, and

\item each ${\bf A}\preceq {\bf R}_1$ with $r_{\pT_{\lceil2r/3\rceil-2}(M_1,G)}(E({\bf A}))< n_{d({\bf A})}$ satisfies ${\bf A}\preceq {\bf Q}$.
\end{enumerate}
Since $M_0$ is vertically $(n_{d({\bf R})}+1)$-connected, each prickle $\A\preceq \R$ with $r_{\pT_{\lceil2r/3\rceil-2}(M_0,G)}(E(\A))<n_{d(\A)}$ satisfies $r_{\pT_{\lceil2r/3\rceil-2}(M_0,G)}(E(\A))=r_{M_0}(E(\A))$ by Lemma \ref{tangle conn}.
This implies that $M_0$ is a valid choice for $M_1$, using \ref{first Q claim}.

Let $\pT=\pT_{\lceil2r/3\rceil-2}(M_1,G)$, and let $C_0\subseteq E(M_1)$ such that $G$ is a restriction of $M_1/C_0$. 
Note that $\pT_{\lceil2r/3\rceil-2}(M_1/C,G/C)$ is a tangle for each $C\subseteq E(G)$ of rank at most two, by Lemma \ref{Dowling tangle}; this is why we work with $\pT$ instead of $\pT_{\lceil2r/3\rceil}(M_1,G)$.
We will show that there is a prickle $\Y\preceq \R_1$ with `small' tangle rank so that $C_0\subseteq \cl_{\pT}(E(\Y))$.
We first prove the following claim, which shows that the union of two subprickles of ${\bf R}_1$ with `small' tangle rank also has `small' tangle rank.


\begin{claim}\label{nested}
If ${\bf A}_1\preceq {\bf R}_1$ and ${\bf A}_2\preceq {\bf R}_1$ such that $r_{\pT}(E({\bf A}_i))\le 3(d({\bf A}_i)+1)n_{d({\bf A}_i)}$ for $i\in \{1,2\}$, then $r_{\pT}(E({\bf A}_1\cup {\bf A}_2))\le 3(d({\bf A}_1\cup {\bf A}_2)+1)n_{d({\bf A}_1\cup {\bf A}_2)}$. 
\end{claim}
\begin{proof}
This is clearly true if $\A_1=\A_2$ or $\A_i=\R_1$, so assume that $\A_1\ne \A_2$ and $\A_i\ne \R_1$ for each $i\in\{1,2\}$.
Then ${\bf A}_1\cup {\bf A}_2$ is a prickle such that $d({\bf A}_1\cup {\bf A}_2)> \max(d({\bf A}_1),d({\bf A}_2))$. 
Assume that $d({\bf A}_1)\ge d({\bf A}_2)$, without loss of generality. 
Then 
\begin{align*}
r_{\pT}(E({\bf A}_1\cup {\bf A}_2))&\le r_{\pT}(E({\bf A}_1))+r_{\pT}(E({\bf A}_2))\\
&\le 6(d({\bf A}_1)+1)n_{d({\bf A}_1)}\\
&\le 3(d({\bf A}_1\cup {\bf A}_2)+1)n_{d({\bf A}_1\cup {\bf A}_2)},
\end{align*}
as desired.
The last inequality holds because $2(i+1)n_i\le (i+1)n_{i+1}$ for all $i\ge 0$.
\end{proof}


Let ${\bf Y}$ be a subprickle of ${\bf R}_1$ with $d({\bf Y})$ maximal such that $r_{\pT}(E({\bf Y}))\le 3(d({\bf Y})+1)n_{d({\bf Y})}$. 
Note that ${\bf Y}$ exists because the trivial subprickle of  ${\bf R}_1$ is a choice for ${\bf Y}$.
One can show that $\Y$ is unique by  \ref{nested}, but we only need the existence of $\Y$.
Also note that ${\bf Q}\preceq {\bf Y}$, by applying \ref{nested} with ${\bf A}_1={\bf Q}$ and ${\bf A}_2={\bf Y}$.


The following claim shows that we can contract certain elements of $E(M_1)-\cl_{\pT}(E({\bf Y}))$ and recover a prickle which satisfies (i) and (ii) and has ${\bf Y}$ as a subprickle.
We apply this claim with ${\bf A}={\bf R}_1$ to show that $C_0\subseteq \cl_{\pT}(E(\Y))$ so that we can apply Theorem \ref{Ramsey} and reduce to the case with a $\DG^-(m+2,\Gamma)$-restriction.
We then apply this claim to show that we can contract two elements and recover a $\DG(m,\Gamma)$-restriction, and still have a prickle with connectivity properties which satisfy the hypotheses of Lemma \ref{hard}.
While we only apply this claim with ${\bf A}={\bf R}_1$, we state it more generally so that we can prove it using induction on $d(\A)$.

\begin{claim} \label{tangle safe}
Let $C\subseteq C_0$ or $C\subseteq E(G)$ so that $|C|\le 2$ and $r_{\pT}(E({\bf Y})\cup C)=r_{\pT}(E({\bf Y}))+|C|$. 
Let $\pT/C=\pT_{\lceil2r/3\rceil-2}(M_1/C,G)$ or $\pT/C=\pT_{\lceil2r/3\rceil-2}(M_1/C,G/C)$ if $C\subseteq C_0$ or $C\subseteq E(G)$, respectively. 
Then for each prickle ${\bf A}\preceq {\bf R}_1$ with ${\bf Y}\preceq {\bf A}$,
the matroid $M_1/C$ has a retract ${\bf A}'$ of ${\bf A}$ so that ${\bf Y}\preceq {\bf A}'$, and each ${\bf Z}\preceq {\bf A}'$ with $r_{\pT/C}(E({\bf Z}))< n_{d({\bf Z})}$ satisfies ${\bf Z}\preceq {\bf Q}$.
\end{claim}
\begin{proof}
Assume that the claim is false for some prickle ${\bf A}$ with $d({\bf A})$ minimum. 
Then $d({\bf A})>d({\bf Y})$ or else ${\bf A}={\bf Y}$ and the claim holds, since Lemma \ref{tangle minor} applied with $X=E(\Y)$ implies that $(M/C)|E(\Y)=M|E(\Y)$ and that $r_{\pT/C}(Y')=r_{\pT}(Y')$ for all $Y'\subseteq E(\Y)$.
Then by Lemma \ref{subprickles} applied to ${\bf A}$ and ${\bf Y}$, there is a collection $\cA$ of subprickles of ${\bf A}$ such that each ${\bf A}_0\in \cA$ satisfies $d({\bf A}_0)=d({\bf A})-1$ and ${\bf Y}\preceq {\bf A}_0$, while $|\cA|\le d({\bf A})-d({\bf Y})$ and $E({\bf A})-E({\bf Y})\subseteq \cup_{{\bf A}_0\in\cA}(E({\bf A})-E({\bf A}_0))$.

Let ${\bf A}_0\in\cA$. We will show that $r_{\pT}(E({\bf A})-E({\bf A}_0))\le 3n_{d({\bf A})}$.
Since $d({\bf A}_0)<d({\bf A})$ and $d({\bf A})$ is minimum, the matroid $M_1/C$ has a retract ${\bf A}_0'$ of ${\bf A}_0$ so that ${\bf Y}\preceq {\bf A}'_0$, and each ${\bf Z}\preceq {\bf A}'_0$ with $r_{\pT/C}(E({\bf Z}))< d({\bf Z})$ satisfies ${\bf Z}\preceq {\bf Q}$.

By Lemma \ref{two spikes} there is a non-empty set $C'\subseteq E(M)-C$ and a collection $\cK$ of prickles of $M_1/C$ such that $|\cK|+|C'|\le 4$, each ${\bf K}\in\cK$ is a retract of ${\bf A}$ with ${\bf A}_0'\preceq {\bf K}$,  and 
\setcounter{equation}{0}
\begin{align}
E({\bf A})-E({\bf A}_0)\subseteq \cl_{M_1}\big(\big(\cup_{{\bf K}\in \cK}(E({\bf K})-E({\bf A}_0'))\big)\cup C'\cup C\big).
\end{align}
Note that for each ${\bf K}\in \cK$ we have ${\bf Y}\preceq {\bf A}_0'\preceq {\bf K}$ and $d({\bf A}_0')=d({\bf K})-1$, by the definition of a retract.
Since the claim is false for ${\bf A}$, for each ${\bf K}\in\cK$ there is some ${\bf Z}\preceq {\bf K}$ with ${\bf Z}\npreceq {\bf Q}$ such that $r_{\pT/C}(E({\bf Z}))<n_{d({\bf Z})}\le n_{d({\bf A})}$. 
Since $d({\bf A}_0')=d({\bf K})-1$ and ${\bf Z}\npreceq {\bf A}_0'$, by Lemma \ref{one smaller} (i) we have $E({\bf K})-E({\bf A}_0')\subseteq E({\bf Z})$, so $r_{\pT/C}(E({\bf K})-E({\bf A}_0'))<n_{d({\bf A})}$. 
Thus, by (1) and the fact that $|\cK|+|C'|\le 4$, we have
\begin{align*}
r_{\pT/C}(E({\bf A})-E({\bf A}_0))&\le r_{\pT/C}(\cup_{{\bf K}\in \cK}(E({\bf K})-E({\bf A}_0')))+|C'|\\
&\le |\cK|(n_{d({\bf A})}-1)+|C'|\\
&\le (|\cK|+|C'|-1)(n_{d({\bf A})}-1)+1\\
&\le 3n_{d({\bf A})}-2.
\end{align*}
Then since $|C|\le 2$ we have
$r_{\pT}(E({\bf A})-E({\bf A}_0))\le 3n_{d({\bf A})}$, using Lemma \ref{tangle closure} if $C\subseteq C_0$ and Lemma \ref{tangle minor} if $C\subseteq E(G)$.

We now use that $|\cA|\le d({\bf A})-d({\bf Y})$ and $E({\bf A})-E({\bf Y})\subseteq \cup_{{\bf A}_0\in\cA}(E({\bf A})-E({\bf A}_0))$.
We have 
\begin{align*}
r_{\pT}(E({\bf A}))&\le r_{\pT}(E({\bf Y}))+r_{\pT}(E({\bf A})-E({\bf Y}))\\
&\le 3(d({\bf Y})+1)n_{d({\bf Y})}+\sum_{{\bf A}_0\in \cA}r_{\pT}(E({\bf A})-E({\bf A}_0))\\
&\le 3(d({\bf Y})+1)n_{d({\bf Y})}+|\cA|3n_{d({\bf A})}\\
&\le 3(d({\bf Y})+1)n_{d({\bf Y})}+3(d({\bf A})-d({\bf Y}))n_{d({\bf A})}\\
&\le 3(d({\bf A})+1)n_{d({\bf A})}.
\end{align*}
But then ${\bf A}\preceq {\bf Y}$ by \ref{nested} and the definition of ${\bf Y}$, so ${\bf A}$ is not a counterexample.
\end{proof}


If there is some $e\in C_0-\cl_{\pT}(E({\bf Y}))$, then $r_{\pT_{\lceil2r/3\rceil-2}(M_1/e,G)}(E({\bf Q}))=r_{M_0}(E({\bf Q}))$ by Lemma \ref{tangle closure} since $E(\Q)\subseteq E(\Y)$.
Then by  \ref{tangle safe} applied with $C=\{e\}$ and ${\bf A}={\bf R}_1$, the matroid $M_1/e$ contradicts the minor-minimality of $M_1$. 
Thus, $C_0\subseteq\cl_{\pT}(E({\bf Y}))$ and so $r_{\pT}(C_0)\le 3(ks+1)n_{ks}\le m_1$. 
So by Lemma \ref{C_0} we have $r_{M_1}(E(G))-r(G)\le \ell^{m_1+1}$.
Then since $r(G)\ge m_1+3+f_{\ref{Ramsey}}(\ell,n,\ell^{m_1+1})$, by Lemma \ref{Ramsey} with $m=m_1+3$ and $d=\ell^{m_1+1}$ there is a $\DG^-(m_1+3,\Gamma)$-restriction $G_1$ of $M_1|E(G)$.

Since $r_{\pT}(E(G_1))\ge r_G(E(G_1))\ge m_1+2$ by Lemma \ref{intersect},  by Lemma \ref{jointless} there is a set  $C\subseteq E(G_1)$ of size at most two so that $r_{\pT}(E(\Y)\cup C)=r_{\pT}(E(\Y))+|C|$ and $G_1/C$ has a $\DG(m_1,\Gamma)$-restriction $G_2$.
Let $\pT/C=\pT_{\lceil2r/3\rceil-2}(M_1/C,G/C)$.
Since $r_{\pT}(E({\bf Y})\cup C)=r_{\pT}(E({\bf Y}))+|C|$, 
by  \ref{tangle safe} applied to ${\bf R}_1$ and $C$, the matroid $M_1/C$ has a retract  ${\bf R}_2$ of ${\bf R}_1$ so that ${\bf Y}\preceq {\bf R}_2$, and each ${\bf A}\preceq {\bf R}_2$ with $r_{\pT/C}(E({\bf A}))< n_{d({\bf A})}$  satisfies ${\bf A}\preceq {\bf Q}$. 
This last condition implies that ${\bf R}_2$ is a $g$-prickle.

\begin{claim}
${\bf R}_2$ is a $g$-prickle.
\end{claim}
\begin{proof}
If not, then there is some porcupine $P$ of $\R_2$ which satisfies $d(P)=1$ and $r_{M_1/C}(P)\le g$ (so $P$ is a spike of rank at most $g$). 
Then the subprickle $\A$ of $\R_2$ consisting of $P$ and otherwise only trivial porcupines
 is not a $g$-prickle, and satisfies $d({\bf A})=1$. 
 But then $r_{M_1/C}(E({\bf A}))< |T|+g<n_1$, and so $r_{\pT/C}(E({\bf A}))\le r_{M_1/C}(E({\bf A}))<n_{d({\bf A})}$, so ${\bf A}\preceq {\bf Q}$.
But then ${\bf A}\preceq {\bf R}_0$ since ${\bf Q}\preceq {\bf R}_0$, which contradicts that ${\bf R}_0$ is a $g$-prickle, since every subprickle of a $g$-prickle is also a $g$-prickle.
\end{proof}

We now show that $M_1/C$ and ${\bf R}_2$ satisfy the connectivity conditions of Lemma \ref{hard} with $J=E(G_2)$.
Since $r_{\pT}(E({\bf Y})\cup C)=r_{\pT}(E({\bf Y}))+|C|$ we have $r_{\pT/C}(E(\Q))=r_{\pT}(E(\Q))$ by Lemma \ref{tangle minor}.
If ${\bf A}\preceq {\bf R}_2$ with ${\bf A}\npreceq {\bf Q}$, then 
\begin{align*}
\kappa_{M_1/C}(E(G_2),E({\bf A}))&\ge \min\big(r_{\pT/C}(E({\bf A})),r_{G/C}(E(G_2))\big)\\
&\ge \min(n_{d({\bf A})},m_1)\\
&\ge n_{d({\bf A})},
\end{align*}
where the first inequality holds by applying Lemma \ref{tangle restriction} to $\pT/C$ with $J=E(G_2)$, and the last inequality holds since $m_1\ge n_{ks}$. 
Since $T\subseteq E({\bf Q})$, by Lemma \ref{hard} with $J=E(G_2)$ the matroid $M_1/C$ has a minor $N$ with $G_2$ a spanning restriction, $N|T=(M_1/C)|T=M_0|T$, and a $g$-retract  ${\bf R}_3$ of ${\bf R}_2$.
Since ${\bf R}_3$ is a retract of ${\bf R}_2$ and is thus a retract of ${\bf R}_0$, each element of $T'$ is the tip of a $g$-porcupine restriction $P$ of $N$ with $d(P)=s$.

Therefore, $M_0$ has a minor $N$ with a $\DG(r(N),\Gamma)$-restriction so that $N|T=M_0|T$ and each element of $T'$ is the tip of a $g$-porcupine restriction  $P$ of $N$ with $d(P)=s$.
Since porcupines and Dowling geometries are simple, we may assume that $N|(E(N)-T)$ is simple.
Let $C_1,D\subseteq E(M_0)$ be disjoint sets so that $N=M_0/C_1\del D$.
Since $N|(E(N)-T)$ is simple, we may assume that $X\subseteq D$.
Thus, $N$ is a minor of $M$, and the theorem holds. 
\end{proof}
\end{subsection}
\end{section}


\begin{section}{The Main Result} \label{main proof chap}
In this section we combine the main results of Sections \ref{crit chap} and \ref{conn chap} to prove Theorem \ref{main corollary}.


\begin{subsection}{The Spanning Clique Case} \label{span sec}
We first combine Corollary \ref{mader}, Theorem \ref{spanning clique}, and Theorem \ref{clique minor} to prove a version of Theorem \ref{main corollary} for matroids with a spanning clique restriction.

\begin{proposition} \label{spanning clique case}
There is a function $f_{\ref{spanning clique case}}:\bZ^3\to \bZ$ so that for all integers $\ell\ge 2$, $t\ge 1$, and $k,n\ge 3$, if $\cM$ is a minor-closed class of matroids such that $U_{2,\ell+2}\notin \cM$, then either 
\begin{itemize}
\item $\cM$ contains $\LG^+(n,\Gamma')$ for some nontrivial group $\Gamma'$, 
\item $\cM$ contains $\DG(k,\Gamma)$ for some group $\Gamma$ with $|\Gamma|\ge t$, or
\item each $M\in \cM$ with a spanning clique restriction satisfies $\elem(M)\le (t-1){r(M)\choose 2}+f_{\ref{spanning clique case}}(\ell,n,k)\cdot r(M)$.
\end{itemize}
\end{proposition}
\begin{proof}
Let $h=h_{\ref{spanning clique}}(\ell,n)$ and $\alpha=\alpha_{\ref{clique minor}}(\ell,n^{2^{\ell}})$, and define $f_{\ref{spanning clique case}}(\ell,n,k)=2\ell k^{2^{\ell}}+h \ell \alpha$.
Suppose that the first two outcomes do not hold for this value of $f_{\ref{spanning clique case}}$, and let $M\in \cM$ be simple with a spanning clique restriction. 
By Theorem \ref{spanning clique} there are disjoint sets $C_1,C_2\subseteq E(M)$ with $r_M(C_1\cup C_2)\le h$ and a $B$-clique $\hat M$ such that $\elem(\hat M)\ge \elem(M/(C_1\cup C_2))$ and $\si(\hat M)$ is isomorphic to a restriction of $M/C_1$. 

\begin{claim} \label{hat M}
$\elem(M)\le \elem(\hat M)+(h\ell\alpha)r(M)$.
\end{claim}
\begin{proof}
Since $\elem(M/(C_1\cup C_2))\le \elem(\hat M)$ it suffices to show that $\elem(M)\le \elem(M/(C_1\cup C_2))+(h\ell\alpha)r(M)$.
Assume for a contradiction that $\elem(M)-\elem(M/(C_1\cup C_2))> (h\ell\alpha) r(M)$.
Let $C\subseteq C_1\cup C_2$ be a maximal independent set so that $\elem(M)-\elem(M/C)\le r_M(C)(\ell\alpha)r(M)$, and let $e\in (C_1\cup C_2)-\cl_M(C)$.
Then 
\begin{align*}
\elem(M/C)-\elem(M/C/e)&=\big(\elem(M)-\elem(M/C/e)\big)- \big(\elem(M)-\elem(M/C)\big)\\
&>r_M(C\cup\{e\})(\ell\alpha) r(M)-r_M(C)(\ell\alpha)r(M) \\
&\ge (\ell\alpha)r(M).
\end{align*}
Since $M\in \cU(\ell)$, this implies that there are greater than $\alpha \cdot r(M)$ long lines of $M/C$ through $e$. 
Let $T_0\subseteq E(M/C)$ be the set of elements of $M/C$ on long lines through $e$, and let $N=(M/C)|T_0$.

Let $(N/e)|T$ be a simplification of $N/e$, so $\elem((N/e)|T)>\alpha \cdot r(N/e)$.
By Theorem \ref{clique minor}, the matroid $(N/e)|T$ has an $M(K_{n^{2^{\ell}}})$-minor. 
Say $(N/X/e)|T_1\cong M(K_{n^{2^{\ell}}})$ for disjoint sets $X,T_1\subseteq T$.
Let $J\subseteq E(N)$ denote the set of elements of $N$ on a line through $e$ and an element of $T_1$.
Each line of $N$ through $e$ and an element of $T_1$ has rank two in $N/X$, since $(N/X/e)|T_1$ is simple.
Thus, $(N/X)|J$ has a rank-$n^{2^{\ell}}$ doubled-clique restriction with tip $e$.
But then $N/X$ has an $\LG^+(n,\Gamma)$-minor for some nontrivial group $\Gamma$, by Proposition \ref{LG minor} with $m=n$, a contradiction.
\end{proof}

Since $\si(\hat M)$ is isomorphic to a restriction of $M/C_1$, the $B$-clique $\hat M$ has no $U_{2,\ell+2}$-minor and no $\DG(k,\Gamma)$-minor with $|\Gamma|\ge t$. Then by Lemma \ref{mader} we have $\elem(\hat M)\le (t-1){r(\hat M)\choose 2}+2\ell k^{2^{\ell}}\cdot r(\hat M)$.
By  \ref{hat M},
$$\elem(M)\le (t-1){r(M)\choose 2}+\Big(2\ell k^{2^{\ell}}+h \ell \alpha\Big)r(M),$$
as desired.
\end{proof}

\end{subsection}


\begin{subsection}{Stacks} \label{stacks cert}
We now prove that a matroid with a spanning clique restriction and a huge $(\cF\cap \cU(t))$-stack restriction has either a projective-geometry minor, or a Dowling-geometry minor over a group of size at least $t$.
Recall that a matroid $M$ is an \emph{$(\mathcal O,b,h)$-stack} if there are disjoint sets $P_1,P_2,\dots,P_h\subseteq E(M)$ such that $\cup_i P_i$ spans $M$, and for each $i\in [h]$ the matroid $(M/(P_1\cup\dots\cup P_{i-1}))|P_i$ has rank at most $b$ and is not in $\mathcal O$.
Note that for each $j\in [h]$, the matroid $M/(P_1\cup \dots \cup P_j)$ is an $(\cO,b,h-j)$-stack.
We say that a matroid $N$ is a \emph{good} $(\cO,b,h-j)$-stack-minor of $M$ if there is some $j\in [h]$ so that $N=M/(P_1\cup \dots \cup P_j)$.
It is clear from this definition that if $M_1$ is a good $(\cO,b,h-j)$-stack-minor of $M$, and $M_2$ is a good $(\cO,b,h-j-i)$-stack-minor of $M_1$ for $i\in [h-j]$, then $M_2$ is a good $(\cO,b,h-j-i)$-stack-minor of $M$.

We first prove a general lemma which shows that stacks are somewhat robust under lifts and projections.

\begin{lemma} \label{stack perturb}
Let $d,b,m\ge 0$ be integers, and let $M$ and $N$ be matroids.
If $M|S$ is an $(\mathcal O,b,m2^d)$-stack and $\dist(M,N)\le d$, then there are sets $C$ and $S'$ so that $(M/C)|S'=(N/C)|S'$ is a good $(\mathcal O,b,m)$-stack-minor of $M|S$.
\end{lemma}
\begin{proof}
Let $d$ be minimal so that the claim is false, so $d>0$. 
Let $M_0$ be a matroid so that $\dist(M,M_0)\le d-1$ and $\dist(M_0,N)=1$. 
By the minimality of $d$, there are sets $C_0$ and $S_0$ so that $(M/C_0)|S_0=(M_0/C_0)|S_0$ is a good $(\cO,b,2m)$-stack-minor of $M|S$, since $\frac{m2^d}{2^{d-1}}=2m$.
Let $M'=M_0/C_0$ and $N'=N/C_0$, and note that $\dist(M',N')\le \dist(M_0,N)=1$.
Let $(S_1,S_2)$ be a partition of $S_0$ so that $M'|S_1$ and $(M'/S_1)|S_2$ are $(\cO,b,m)$-stacks.

It suffices to show that there are sets $C_1$ and $S_3$ so that $(M'/C_1)|S_3=(N'/C_1)|S_3$ is a good $(\cO,b,m)$-stack-minor of $M'|S_0$, because then $(M/(C_0\cup C_1))|S_3=(N/(C_0\cup C_1))|S_3$ is a good $(\cO,b,m)$-stack-minor of $M|S$, and the lemma holds with $C=C_0\cup C_1$ and $S'=S_3$. 
Let $K$ be a matroid so that $\{K/f,K\del f\}=\{M',N'\}$ for some $f\in E(K)$.
First assume that $(K\del f,K/f)=(M',N')$.
If $f\notin\cl_{K}(S_1)$, then $N'|S_1$ is an $(\cO,b,m)$-stack and we take $C_1=\varnothing$.
If $f\in\cl_{K}(S_1)$, then 
$$(N'/S_1)|S_2=(K/S_1\del f)|S_2=(M'/S_1)|S_2,$$
so we take $C_1=S_1$.
Now assume that $(K/f ,K\del f)=(M',N')$.
If $f\notin\cl_{K}(S_1)$, then $N'|S_1=K|S_1=(K/f)|S_1=M'|S_1$, so we take $C=\varnothing$.
If $f\in\cl_{K}(S_1)$, then $(N'/S_1)|S_2=(K/S_1/f)|S_2=(M'/S_1)|S_2,$
so we take $C_1=S_1$.
\end{proof}

We will need the following corollary of Lemma \ref{stack perturb} in the next section.

\begin{corollary} \label{stack cor}
Let $d,b,m\ge 0$ be integers, and let $M$ be a matroid.
If $M|S$ is an $(\mathcal O,b,m2^d)$-stack and $r_M(X)\le d$, then there are sets $C$ and $S'$ so that $X\subseteq C$ and $(M/C)|S'$ is a good $(\mathcal O,b,m)$-stack-minor of $M|S$.
\end{corollary}
\begin{proof}
Let $N$ be the matroid with ground set $E(M)$ so that $N\del X=M/X$ and $r_N(X)=0$.
Then $\dist(M,N)\le d$, so the statement holds by Lemma \ref{stack perturb}.
\end{proof}


We now combine Theorem \ref{spanning clique} and Lemma \ref{stack perturb}.

\begin{proposition} \label{stack certificate}
There is a function $h_{\ref{stack certificate}}\colon \bZ^3\to \bZ$ so that for all integers $\ell, t,b\ge 2$ and $n,k\ge 3$, if $M\in \cU(\ell)$ has a spanning clique restriction and an $(\cF\cap \cU(t),b,h_{\ref{stack certificate}}(\ell,k,n))$-stack restriction $S$, then $M$ has either a rank-$n$ projective-geometry minor or a $\DG(k,\Gamma)$-minor with $|\Gamma|\ge t$. 
\end{proposition}
\begin{proof}
Let $h=h_{\ref{spanning clique}}(\ell,n)$, let $m=k^{2^{\ell}}$, and define $h_{\ref{stack certificate}}(\ell,k,n)=m2^h$.
Let $B_0$ be a frame for the spanning clique restriction of $M$.
Assume that $M$ has no rank-$n$ projective-geometry minor.
By Theorem \ref{spanning clique} there is some $B_1\subseteq B_0$ and a $B_1$-clique $N$ so that 
$\dist(M,N)\le h$ and $\si(N)$ is isomorphic to a minor of $M$.
By Lemma \ref{stack perturb}, $N$ has a contract-minor $N_1$ with a set $S_1$ so that $N_1|S_1$ is an  $(\cF\cap \cU(t),b,m)$-stack. 
Since $N_1$ is a frame matroid, the matroid $N_1|S_1$ is an $(\cU(t),b,m)$-stack.
Then each part of this stack has a $U_{2,t+2}$-minor, so by contracting elements of each part of this stack we see that $N_1$ has a contract-minor $N_2$ with $S_2\subseteq S_1$ so that $N_2|S_2$ is a spanning $(\cU(t),2,m)$-stack restriction of $N_2$.
Note that $N_2$ is a $B$-clique for some $B\subseteq B_1$, since $N_2$ is a contract-minor of $N$.

Since $N_2|S_2$ is a spanning $(\cU(t),2,m)$-stack restriction of $N_2$, there are disjoint sets $P_1,\dots, P_m\subseteq E(N_2)$ so that $\cup_i P_i$ spans $N_2$ and for each $i\in [m]$ the matroid $(N_2/(P_1\cup\dots\cup P_{i-1}))|P_i$ is isomorphic to $U_{2,t+2}$.
Since $t\ge 2$ and each line of length at least four spans two elements of $B$,
this gives a partition $(B_1,\dots,B_m)$ of $B$ so that $|B_i|=2$ and $P_i\subseteq \cl_{N_2/(P_1\cup \dots\cup P_{i-1})}(B_i)$ for each $i\in [m]$.
Then $P_1\cup \dots \cup P_{i-1}$ and $B_1\cup \dots \cup B_{i-1}$ span the same flat of $N_2$, so $$(N_2/(P_1\cup \dots\cup P_{i-1}))|P_i=(N_2/(B_1\cup \dots\cup B_{i-1}))|P_i.$$
Since $P_i\subseteq \cl_{N_2/(B_1\cup \dots\cup B_{i-1})}(B_i)$ and $|P_i|\ge 4$, this implies that there are two elements $e\in P_i$ for which the unique circuit of $N_2|(B\cup\{e\})$ contains $B_i$. 
Since $N_2$ is a frame matroid and $|B_i|=2$, this circuit contains no other elements of $B$. 
Thus, two elements of $P_i$ are spanned in $N_2$ by $B_i$ and are not parallel to either element of $B_i$, which implies that $P_i\subseteq \cl_{N_2}(B_i)$ for each $i\in [m]$.
Since $\elem(M|(B_i\cup P_i))\ge t+2$ for each $i\in [m]$, the matroid $N_2$ has a $\DG(k,\Gamma)$-minor with $|\Gamma|\ge t$ by Lemma \ref{DGminor}.
Since $\DG(k,\Gamma)$ is simple and $\si(N)$ is isomorphic to a minor of $M$, the matroid $M$ has a $\DG(k,\Gamma)$-minor.
\end{proof}
\end{subsection}


\begin{subsection}{Small Spikes with Common Tip} \label{spikes cert}
We now prove that a matroid with a spanning clique restriction and lots of nearly skew small spikes with a common tip has an $\LG^+(n,\Gamma)$-minor for some nontrivial group $\Gamma$. 
We need a continuation of Lemma \ref{stack perturb} in the case that $N$ is a $B$-clique.

\begin{lemma} \label{frame stack minor}
There is a function $f_{\ref{frame stack minor}}\colon \bZ^3\to \bZ$ so that for all integers $d,b,m\ge 0$, if $M$ and $N$ are matroids with $\dist(M,N)\le d$ so that $N$ is a $B$-clique and $M|S=N|S$ is an $(\cO,b,f_{\ref{frame stack minor}}(d,b,m))$-stack, then there are sets $C,D,B',S'$ so that 
\begin{itemize}
\item $M/C\del D=N/C\del D$,
\item $N/C\del D$ is a $B'$-clique, and
\item $(M/C\del D)|S'=(N/C\del D)|S'$ is a good $(\cO,b,m)$-stack-minor of $M|S=N|S$.
\end{itemize} 
\end{lemma}
\begin{proof}
Define $f_{\ref{frame stack minor}}(0,b,m)=1$, and inductively define 
$$f_{\ref{frame stack minor}}(d,b,m)=f_{\ref{frame stack minor}}(d-1,b,m)\big(1+2^{2b\cdot f_{\ref{frame stack minor}}(d-1,b,m)}\big),$$
for $d>0$.
Let $d$ be minimal so that the claim is false, so $d>0$, and let $m_1=f_{\ref{frame stack minor}}(d-1,b,m)$.
Let $M_1$ be a matroid so that $\dist(N,M_1)=1$ and $\dist(M_1,M)\le d-1$.
Let $(S_1,S_2)$ be a partition of $S$ so that $N|S_1$ is an $(\cO,b,m_1)$-stack and $(N/S_1)|S_2$ is an $(\cO,b,m_12^{2bm_1})$-stack.
We consider two cases, depending on whether $M_1$ is a projection or lift of $N$.
Recall that any contract-minor of $N$ is a $B_1$-clique for some $B_1\subseteq B$.

First assume that there is a matroid $K$ with an element $f$ so that $K\del f=N$ and $K/f=M_1$. 

\begin{claim}
There are sets $C_1,D_1,B_1,S_1'$ so that $N/C_1\del D_1=M/C_1\del D_1$ is a $B_1'$-clique, and $(N/C_1\del D_1)|S_1'=(M_1/C_1\del D_1)|S_1'$ is a good $(\cO,b,m_1)$-stack-minor of $N|S$.
\end{claim}
\begin{proof}
Let $B_1'$ be a minimal subset of $B$ so that $S_1\subseteq\cl_N(B_1')$, and note that $|B_1'|\le 2r_N(S_1)\le 2bm_1$.
If $f\notin\cl_N(B_1')$, then let $D_1=E(N)-\cl_N(B_1')$.
Then $N\del D_1$ is a $B_1'$-clique, and $N\del D_1=M_1\del D_1$ since $f\notin\cl_K(E(N)-D_1)$, so the claim holds with $(C_1,D_1,B_1,S_1')=(\varnothing,D_1,B_1',S_1)$.

If $f\in \cl_N(B_1')$, then by Corollary \ref{stack cor} applied with $(M,X,d,m)=(N,B_1',2bm_1,m_1)$, there are sets $C'$ and $S_1'$ so that $B_1'\subseteq C'$ and $(N/C')|S_1'$ is a good $(\cO,b,m_1)$-stack-minor of $N|S$.
Thus, the claim holds with $(C_1,D_1,S_1')=(C',\varnothing,S')$, since $f\in \cl_N(C')$ and $N/C'$ is a contract-minor of $N$.
\end{proof}

Now assume that there is a matroid $K$ with an element $f$ so that $K/f=N$ and $K\del f=M_1$.

\begin{claim} \label{C3}
There is a set $C_3\subseteq E(N)$ so that $|C_3|\le 3$ and $f\in \cl_K(C_3)$.
\end{claim}
\begin{proof}
If $B$ is a basis of $M_1$, then $r(M_1)=r(N)$.
But then $f$ is a coloop of $K$, which implies that $N=M_1$, and this contradicts that $d(N,M_1)>0$.
Thus, there is some element $e\in M_1-\cl_{M_1}(B)$.
Since $N$ is a $B$-clique there are elements $b,b'\in B$ so that $e\in \cl_N(\{b,b'\})$.
Since $r_{M_1}(\{b,b',e\})=3$ and $r_N(\{b,b',e\})=2$ it follows that $f\in \cl_K(\{b,b',e\})$.
\end{proof}

If $f\in \cl_K(S_1)$, then let $C_1=S_1$.
Then there are sets $B_1\subseteq B$ and $S_2'\subseteq S_2$ so that $N/C_1=M/C_1$ is a $B_1$-clique, and $(N/C_1)|S_2'=(M_1/C_1)|S_2'$ is a good $(\cO,b,m_1)$-stack-minor of $N|S$.
If $f\notin\cl_K(S_1)$, then let $C_3$ be the set given by \ref{C3}.
By Corollary \ref{stack cor} applied with $(M,X,d,m)=(N,C_3,3,m_1)$, there are sets $C'$ and $S_1'$ so that $C_3\subseteq C'$ and $(N/C')|S_1'$ is a good $(\cO,b,m_1)$-stack-minor of $N|S$.
Let $C_1=C'$.
Then $N/C_1=M_1/C_1$ since $f\in \cl_N(C_1)$, and $(N/C_1)|S_{1}'=(M_1/C_1)|S_{1}'$ is a good $(\cO,b,m_1)$-stack-minor of $N|S$.


In all cases, we have shown that there are sets $C_1,D_1,B_1,S_1$ so that $N/C_1\del D_1=M_1/C_1\del D_1$ is a $B_1$-clique, and $(N/C_1\del D_1)|S_1=(M_1/C_1\del D_1)|S_1$ is a good $(\cO,b,m_1)$-stack-minor of $N|S$.
Let $N'=N/C_1\del D_1$, and $M'=M/C_1\del D_1$.
Then $N'$ is a $B_1$-clique, and $\dist(N'/C_1\del D_1,M/C_1\del D_1)\le \dist(M_1,M)\le d-1$.
By the minimality of $d$, the definition of $m_1$, and the fact that $N'$ is a $B_1$-clique, there are sets $C_2,D_2,B_1',S_1'$ so that 
\begin{itemize}
\item $M'/C_2\del D_2=N'/C_2\del D_2$,
\item $N'/C_2\del D_2$ is a $B_1'$-clique, and
\item $(M'/C_2\del D_2)|S_1'=(N'/C_2\del D_2)|S_1'$ is a good $(\cO,b,m)$-stack-minor of $M'|S_1=N'|S_1$.
\end{itemize} 
Let $C=C_1\cup C_2$ and $D=D_1\cup D_2$.
Then $M/C\del D=N/C\del D$ is a $B_1'$-clique, and $(M/C\del D)|S_1'=(N/C\del D)|S_1'$  is a good $(\cO,b,m)$-stack-minor of $M|S=N|S$.
\end{proof}


To prove the following proposition we apply Theorem \ref{spanning clique} and Lemma \ref{frame stack minor}, and then apply Lemma \ref{find d clique} to find an $\LG^+(n,\Gamma)$-minor for some nontrivial group $\Gamma$.
Given a collection $\mathcal Z$ of sets, we will write $\cup\mathcal Z$ for $\cup_{Z\in\mathcal Z}Z$, for convenience.

\begin{proposition} \label{spikes certificate}
There is a function $m_{\ref{spikes certificate}}\colon \bZ^2\to \bZ$ so that for all integers $\ell\ge 2$ and $n\ge 3$, if $M\in\cU(\ell)$ has a spanning clique restriction and there is some $e\in E(M)$ and a collection $\cS$ of $m_{\ref{spikes certificate}}(\ell,n)$ mutually skew sets in $M/e$ so that for each $S\in\cS$, the matroid $M|(S\cup \{e\})$ is a spike of rank at most four with tip $e$, then $M$ has an $\LG^+(n,\Gamma)$-minor for some nontrivial group $\Gamma$. 
\end{proposition}
\begin{proof}
Let $h=h_{\ref{spanning clique}}(\ell,n)$ and $m_1=n^{2^{2\ell+1}}$, and define $m_{\ref{spikes certificate}}(\ell,n)=f_{\ref{frame stack minor}}(h,2,2m_1)2^{h}$.
Let $m=n^{2^{\ell}}$, and note that $m_1=m^{2^{\ell}+1}$.
Assume that $M$ has no $\LG^+(n,\Gamma)$-minor with $|\Gamma|\ge 2$; this implies that $M$ has no rank-$n$ projective-geometry minor.
Let $M$ with $\cS$ be a counterexample so that $M$ is minor-minimal.
Then each $S\in \cS$ satisfies $r_{M/e}(S)=2$, since every rank-4 spike has a rank-3 spike as a minor.
For each $S\in \cS$, let $X_S$ be a transversal of the parallel classes of $(M/e)|S$, and let $\cX=\{X_S\colon S\in \cS\}$.
Let $\cO=\cM-\{(M/e)|X\colon X\in\cX\}$, where $\cM$ is the class of all matroids.
Then $(M/e)|\cup\cX$ is an $(\cO,2,m_{\ref{spikes certificate}}(\ell,n))$-stack.
Note that if $M'$ is a good $(\cO,2,j)$-stack-minor of $(M/e)|\cup\cX$, then there is collection of $j$ sets in $\cX$ which are mutually skew in $M'$, by the definition of a good $(\cO,2,j)$-stack-minor.

Since $M/e$ has a spanning clique restriction, by Theorem \ref{spanning clique} there is a $B$-clique $N$ so that $\dist(M/e,N)\le h$.
By Lemma \ref{stack perturb} with $m=f_{\ref{frame stack minor}}(h,2,2m_1)$ and $d=h$ there are sets $C_1$ and $S_1$ so that $(M/e/C_1)|S_1=(N/C_1)|S_1$ is a good $(\cO,2,f_{\ref{frame stack minor}}(h,2,2m_1))$-stack-minor of $(M/e)|\cup\cX$.
Note that $\dist(M/e/C_1,N/C_1)\le \dist(M/e,N)\le h$ and that $N/C_1$ is a $B_1$-clique for some $B_1\subseteq B$.
Then by Lemma \ref{frame stack minor} with $M=M/e/C_1$ and $N=N/C_1$ and $S=S_1$, there are sets $C,D,B',S'$ so that $M/e/(C_1\cup C)\del D=N/(C_1\cup C)\del D$ is a $B'$-clique, and 
$(N/(C_1\cup C)\del D)|S'$ is a good $(\cO,2,2m_1)$-stack-minor of $(N/C_1)|S_1$.

Let $C'$ be a maximal subset of $C_1\cup C$ so that $e\notin\cl_M(C_1\cup C)$, and let $M'=M/C'\del D$.
Then $M'/e$ is a $B'$-clique, and there is a collection $\cX'\subseteq\cX$ of $2m_1$ mutually skew subsets of $M'/e$.
Let $\cS'=\{S\in \cS\colon X_S\in \cX'\}$.
For each $S\in \cS'$ we have $M'|(S\cup\{e\})=M|(S\cup\{e\})$, since $e$ is a nonloop of $M'$ and $(M'/e)|X=(M/e)|X$. 
Thus, $\cS'$ is a collection of mutually skew subsets of $M'/e$ so that for each  $S\in \cS'$, the matroid $M'|(S\cup\{e\})$ is a rank-3 spike with tip $e$.
Let $C_2\subseteq E(M')$ have minimum size so that $C_2\cup(\cup \cS')$ spans $M'$, and let $M_2$ be a simplification of $M'/C_2$.
Then $M_2|(\{e\}\cup(\cup \cS'))=M'|(\{e\}\cup(\cup \cS'))$ and $M_2/e$ is a $B_2$-clique for some $B_2\subseteq B'$, while $\cup \cS'$ spans $M_2$.

\begin{claim} \label{get Y}
There is a collection $\cY$ of $m_1$ pairwise-disjoint 2-subsets of $B_2$ so that each $Y\in\cY$ spans a nontrivial parallel class of $M_2/e$ which contains neither element of $Y$.
\end{claim}
\begin{proof}
Let $\cY$ be a maximal such collection of subsets of $B_2$, and assume for a contradiction that $|\cY|<m_1$.
Since each element of $\cup \cS'$ is in a nontrivial parallel class of $M_2/e$, the maximality of $|\cY|$ implies that each nonloop element of $(\cup\cS')-(\cup \cY)$ is parallel in $M_2/e/(\cup \cY)$ to an element of $B_2-(\cup \cY)$.
Since $r_{M_2/e}(\cup \cY)\le 2(m_1-1)<|\cS'|$ and the sets in $\cS'$ are mutually skew in $M_2/e$, there is some set $S\in \cS'$ so that $(M_2/e/(\cup \cY))|S=(M_2/e)|S$.
But since $(M_2/e)|S$ contains a rank-2 circuit, this implies that three elements of $B_2-(\cup \cY)$ are in a circuit in $(M_2/e/(\cup \cY))$, which contradicts that $B_2$ is independent in $M_2/e$.
\end{proof}

Lemma \ref{find d clique} and \ref{get Y} and the definition of $m_1$ imply that $M_2$ has a rank-$n^{2^{\ell}}$ doubled-clique minor.
But then Proposition \ref{LG minor} with $m=n$ implies that $M_2$ has an $\LG^+(n,\Gamma)$-minor for some nontrivial group $\Gamma$, a contradiction.
\end{proof}
\end{subsection}


\begin{subsection}{Porcupines} \label{porcs cert}
The following proposition shows that any matroid with a large independent set so that each element is the tip of a $g$-porcupine is not a bounded distance from a frame matroid.
We will apply this with $h=h_{\ref{spanning clique}}(\ell,n)$ to find a rank-$n$ projective-geometry minor.
Recall that if $P$ is a porcupine with tip $f$, then we write $d(P)$ for the corank of $\si(P/f)$.

\begin{proposition} \label{porcupine cshift}
For each integer $h\ge 0$, if $M$ is a matroid with a size-$(h+1)$ independent set $S$ so that each element is the tip of a $(5\cdot 2^h)$-porcupine $P$ with $d(P)=h+1$, then there are no sets $C_1,C_2\subseteq E(M)$ with $r_M(C_1\cup C_2)\le h$ and a frame matroid $N$ on ground set $E(M)$ such that 
\begin{enumerate}[$(*)$]
\item for all $X\subseteq E(M)-(C_1\cup C_2)$, if $(M/(C_1\cup C_2))|X$ is simple, then $N|X=(M/C_1)|X$.
\end{enumerate}
\end{proposition}
\begin{proof}
Let $e\in S-\cl_M(C_1\cup C_2)$, and note that $e$ is a nonloop of $N$ by $(*)$.
We will show that $e$ is the tip of a spike of rank at least five in $N$.
Let $P$ be a $(5\cdot 2^{h})$-porcupine restriction of $M$ with tip $e$ and $d(P)=h+1$, and let $(P/e)|T_0$ be a simplification of $P/e$.
By Lemma \ref{projections1} there is some $T_1\subseteq T_0$ so that $(M/e/(C_1\cup C_2))|T_1$ has corank $h+1$ and girth at least five.
In particular, $(M/e/(C_1\cup C_2))|T_1$ is simple.

Then $(M/e/C_1)|T_1$ has corank at least one and girth at least five, since $r_M(C_2)\le h$.
So there is some $T_2\subseteq T_1$ so that $(M/C_1/e)|T_2$ is a circuit of size at least five.
Let $S\subseteq E(P)$ be the union of lines of $P$ through $e$ and an element of $T_2$.

\begin{claim}
$(M/C_1)|S$ is a spike of rank at least five with tip $e$, and $(M/C_1/e)|T_2$ is a simplification of $(M/C_1/e)|(S-\{e\})$.
\end{claim}
\begin{proof}
First note that each line of $P|S$ through $e$ has rank two in $M/C_1$, or else $(M/C_1/e)|T_2$ contains a loop.
Also, each line of $P|S$ through $e$ is a line of $(M/C_1)|S$, or else two elements of $(M/C_1/e)|T_2$ are parallel.
Since $(M/C_1/e)|T_2$ is a circuit of size at least five, the claim holds.
\end{proof}

Since $(M/e/(C_1\cup C_2))|T_2$ is simple, no two elements of $S$ are parallel in $M/(C_1\cup C_2)$, and $S$ is disjoint from $\cl_{M/C_1}(C_2)$.
Thus, $(M/(C_1\cup C_2))|S$ is simple.
By $(*)$ we have $N|S=(M/C_1)|S$, so $N$ has a spike restriction of rank at least five, which contradicts that $N$ is a frame matroid.
\end{proof}

Theorem \ref{spanning clique} and Proposition \ref{porcupine cshift} combine to give the following corollary, which we will use to find an $\LG^+(n,\Gamma)$-minor for some nontrivial group $\Gamma$.

\begin{corollary} \label{porcs cor}
Let $\ell,n\ge 2$ be integers, and let $h=h_{\ref{spanning clique}}(\ell,n)$.
 If $M\in\cU(\ell)$ has a spanning $B$-clique restriction, and a size-$(h+1)$ independent set $S$ for which each element is the tip of a $(5\cdot 2^h)$-porcupine $P$ with $d(P)=h+1$, then $M$ has a rank-$n$ projective-geometry minor.
\end{corollary}
\end{subsection}


\begin{subsection}{The Main Proof}
We now prove a result which readily implies Theorem \ref{main corollary}.
The proof applies Theorems \ref{new reduction}, \ref{critical structure} and \ref{main porcupines} in that order, and then uses the previous results of this section.

\begin{theorem} \label{main}
There is a function $f_{\ref{main}}:\bZ^4\to \bZ$ so that for all integers $\ell\ge 2$, $t\ge 1$ and $n,k\ge 3$, if $\cM$ is a minor-closed class of matroids such that $U_{2,\ell+2}\notin \cM$, then either 
\begin{itemize}
\item $\cM$ contains $\LG^+(n,\Gamma')$ for some nontrivial group $\Gamma'$, or
\item $\cM$ contains $\DG(k,\Gamma)$ for some group $\Gamma$ with $|\Gamma|\ge t$, or
\item each $M\in \cM$ satisfies $\elem(M)\le (t-1){r(M)\choose 2}+f_{\ref{main}}(\ell,n,k,t)\cdot r(M)$.
\end{itemize}
\end{theorem}
\begin{proof}
Let $h_0=h_{\ref{spanning clique}}(\ell,n)$, let $m_0=m_{\ref{spikes certificate}}(\ell,n)$, and let $h_1=h_{\ref{stack certificate}}(\ell,k,n)$.
Then let $h=\max(h_0,m_0,h_1)$; we will use $h$ to apply Theorem \ref{critical structure}. 
Let $m_1=\max(4m_0, 15h_12^h)$, which will be an upper bound on the rank of a stack or collection of small spikes which we find.
Let $$n_0=\max\Big(k, r_{\ref{main porcupines}}(\ell,3,h_0+1,n,h_0+1,5\cdot 2^{h_0}), r_{\ref{Dowling minor}}(\ell,m_1,n)\Big),$$
which is the rank of a Dowling-geometry minor that we will find,
and let $$s_1=\max\Big(s_{\ref{main porcupines}}(h_0+1,h_0+1,5\cdot 2^{h_0})+h2^{h+7},m_1+2+h2^{h+7},2^{15h}\Big).$$
Define $$f_1=\max\Big(f_{\ref{spanning clique case}}(\ell,n,k),\ell^{28h},\alpha_{\ref{clique minor}}(\ell,n_0+h2^{h+7})\Big),$$ 
to handle matroids with a spanning clique restriction and the case $t=1$, and
$$r_1=r_{\ref{new reduction}}\Big(\frac{t-1}{2},f_1,0,\ell,1,s_1\Big),$$
and finally 
$$f_{\ref{main}}(\ell,t,n,k)=\max(\ell^{r_1},f_1).$$
Let $p(r)=(t-1){r\choose 2}+f_{\ref{main}}(\ell,t,n,k)\cdot r$.
Assume that the third outcome of the theorem statement does not hold, and that $\cM$ contains no $\LG^+(n,\Gamma)$ with $|\Gamma|\ge 2$.

Let $M\in \cM$ be a matroid so that $\elem(M)>p(r(M))$.
Then $r(M)\ge r_1$, by the definition of $f_{\ref{main}}(\ell,t,n,k)$ and the fact that $M\in\cU(\ell)$.
We will show that $M$ has a $\DG(k,\Gamma)$-minor with $|\Gamma|\ge t$.
By Theorem \ref{new reduction} with $r=1$ and $s=s_1$, the matroid $M$ has a minor $N$ so that $r(N)\ge 1$ and $\elem(N)>p(r(N))$, while $N$ either has a spanning clique restriction, or is vertically $s_1$-connected and has an $s_1$-element independent set $S$ so that each $e\in S$ satisfies $\elem(N)-\elem(N/e)>p(r(N))-p(r(N)-1)$.
Since $\elem(N)>p(r(N))$ and $f_1\ge f_{\ref{spanning clique case}}(\ell,n,k)$, the matroid $N$ does not have a spanning clique restriction by Theorem \ref{spanning clique case}.
Since $\elem(N)> \alpha_{\ref{clique minor}}(\ell,n_0+h2^{h+7})\cdot r(N)$ by the definition of $f_1$, the matroid $N$ has an $M(K_{n_0+1+h2^{h+7}})$-minor $G$ by Theorem \ref{clique minor}.
Since $n_0\ge k$, this proves the result in the case that $t=1$, so we may assume that $t\ge 2$.

By Theorem \ref{critical structure} with $h=\max(h_0,m_0,h_1)$ and the facts that $f_1\ge \ell^{28h}$ and $s_1\ge 2^{15h}$, there is some $C\subseteq E(N)$ with $r_N(C)\le h2^{h+7}$ so that $N/C$ has either
\begin{enumerate}[(i)]
\item an $(\cF\cap \cU(t), 15\cdot 2^h, h_1)$-stack restriction,

\item an element $e$ and a collection $\cS$ of $m_0$ mutually skew sets in $N/(C\cup\{e\})$ such that for each $R\in \cS$, the matroid $(N/C)|(R\cup\{e\})$ is a spike of rank at most four with tip $e$, or

\item a size-($h_0+1)$ independent set so that each element is the tip of a $(5\cdot 2^{h_0})$-porcupine restriction $P$ of $N/C$ with $d(P)=h_0+1$.
\end{enumerate}

Note that $N/C$ has an $M(K_{n_0+1})$-minor since $r_N(C)\le h2^{h+7}$ and $N$ has an $M(K_{n_0+1+h2^{h+7}})$-minor.
Also, $N/C$ is vertically $(s_1-h2^{h+7})$-connected since $N$ is vertically $s_1$-connected and $r_M(C)\le h2^{h+7}$.
Note that $s_1-h2^{h+7}\ge m_1+2$.

If (i) holds, then $N/C$ has a rank-$n$ projective-geometry minor or a $\DG(k,\Gamma)$-minor with $|\Gamma|\ge t$ by Theorem \ref{Dowling minor} with $m=m_1$ and Proposition \ref{stack certificate}.
If (ii) holds, then $N/C$ has an $\LG^+(n,\Gamma)$-minor for some nontrivial group $\Gamma$ by Theorem \ref{Dowling minor} with $m=m_1$, and Proposition \ref{spikes certificate}.

If (iii) holds, then by Theorem \ref{main porcupines} with $s=k=h_0+1$ and $g=5\cdot 2^{h_0}$ and $m=3$, the matroid $N/C$ has a minor with a spanning clique restriction and a size-$(h_0+1)$ independent set so that each element is the tip of an $5\cdot 2^{h_0}$-porcupine $P$ with $d(P)=h_0+1$.
By Corollary \ref{porcs cor} and the definition of $h_0$, the matroid $N/C$ has a rank-$n$ projective-geometry minor and thus an $\LG^+(n,\Gamma)$-minor for some nontrivial group $\Gamma$. 
\end{proof}


We now prove Theorem \ref{main corollary}.

\begin{statement}{Theorem \ref{main corollary}}
For each integer $\ell\ge 2$, if $\cM$ is a minor-closed class of matroids so that $U_{2,\ell+2}\notin \cM$,  then there is finite group $\Gamma$ and a constant $c_{\cM}$ so that either
\begin{enumerate}[(1)]
\item $h_{\cM}(r)\le c_{\cM}\cdot r$ for all $r\ge 0$, or

\item $\Gamma$ is nontrivial and $\cM$ contains all $\Gamma$-lift matroids, or

\item $|\Gamma|{r\choose 2}+r\le h_{\cM}(r)\le |\Gamma|{r\choose 2}+c_{\cM}\cdot r$ for all $r\ge 0$, and $\cM$ contains all $\Gamma$-frame matroids.
\end{enumerate}
\end{statement}
\begin{proof}
Assume that (1) does not hold for $\cM$, and that there is no group $\Gamma$ for which (2) holds. 

\begin{claim} \label{upper}
There is an integer $c'$ so that $h_{\cM}(r)\le c'\cdot r^2$ for all $r\ge 0$.
\end{claim}
\begin{proof}
Assume for a contradiction that there is a prime power $p^k$ for which $\cM$ contains all $\GF(p^k)$-representable matroids.
Then $\cM$ contains $\LG^+(n,\bZ_p)$ for each integer $n\ge 3$, by Theorem \ref{LG rep}.
Since the class of $\GF(p^k)$-representable matroids is closed under parallel extensions and adding loops, this implies that $\cM$ contains all $\bZ_p$-lift matroids.
But then (2) holds with $\Gamma=\bZ_p$, a contradiction.
By the Growth Rate Theorem, the claim holds.
\end{proof}

The following claim is not difficult.

\begin{claim} \label{simple}
Let $\Gamma$ be a finite group. Then
\begin{enumerate}[(i)]
\item If $\cM$ contains $\LG^+(k,\Gamma)$ for infinitely many integers $k$, then $\cM$ contains all $\Gamma$-lift matroids, and 

\item If $\cM$ contains $\DG(k,\Gamma)$ for infinitely many integers $k$, then $\cM$ contains all $\Gamma$-frame matroids.
\end{enumerate}
\end{claim}
\begin{proof}
Using the fact that if $M\cong \LG^+(r(M),\Gamma)$ and $e\ne e_0$, then $\si(M/e)\cong \LG^+(r(M)-1,\Gamma)$, it is not hard to show that $\LG^+(2n,\Gamma)$ has a minor $N$ so that $\si(N)\cong \LG^+(n,\Gamma)$, and each parallel class of $N$ has size at least two.
This can be done by contracting a set of elements which forms a maximum matching of the $\Gamma$-labeled graph associated with $\LG^+(2n,\Gamma)$.
Combined with the fact that each simple rank-$n$ $\Gamma$-lift matroid is a restriction of $\LG^+(n,\Gamma)$, this shows that each $\Gamma$-lift matroid is a minor of $\LG^+(n,\Gamma)$ for some sufficiently large integer $n$.
This proves (i), and (ii) is similarly proved.
\end{proof}

The following claim shows that we may apply Theorem \ref{main}.
Note that $\cM$ contains all graphic matroids, or else (1) holds by Theorem \ref{clique minor}.
Let $t\in \bZ$ be maximal so that $\cM$ contains a rank-$k$ Dowling geometry with group size $t$ for infinitely many integers $k$; this integer $t$ is well-defined by \ref{upper}.

\begin{claim} \label{final}
There is an integer $n\ge 3$ so that $\cM$ contains no $\LG^+(n,\Gamma)$ with $|\Gamma|\ge 2$, and no $\DG(n,\Gamma)$ with $|\Gamma|\ge t+1$.
\end{claim}
\begin{proof}
By \ref{upper}, there is an integer $n_0$ so that for all $n\ge n_0$, $\cM$ contains no $\LG^+(n,\Gamma)$ or $\DG(n,\Gamma)$ with $|\Gamma|>c'$.
For each nontrivial group $\Gamma$ of size at most $c'$, there is an integer $n_{\Gamma}$ so that 
$\cM$ contains no $\LG^+(n,\Gamma)$ with $n\ge n_{\Gamma}$; otherwise (2) holds by \ref{simple} (i).
For each group $\Gamma$ with $t+1\le |\Gamma|\le c'$, there is an integer $n'_{\Gamma}$ so that 
$\cM$ contains no $\DG(n,\Gamma)$ with $n\ge n'_{\Gamma}$, by the maximality of $t$. 
Let $n$ be the maximum of $n_0$, and all $n_{\Gamma}$ and $n'_{\Gamma}$.
\end{proof}

Since there are finitely many groups of size $t$, there is some group $\Gamma$ such that $|\Gamma|=t$ and $\cM$ contains $\DG(k,\Gamma)$ for infinitely many integers $k$.
By \ref{simple} (ii), $\cM$ contains all $\Gamma$-frame matroids.
By Theorem \ref{main} and \ref{final}, there is a constant $c_{\cM}$ so that $h_{\cM}(r)\le |\Gamma|{r\choose 2}+c_{\cM}\cdot r$ for all $r\ge 0$.
\end{proof}

\end{subsection}


\begin{subsection}{$\Delta$-Modular Matrices}
We conclude this section by using Theorem \ref{main} to prove Theorem \ref{IP}, which has an application in integer programming theory. 
Given a positive integer $\Delta$, we say that an integer matrix $A$ is \emph{$\Delta$-modular} if the determinant of each $\rank(A)\times \rank(A)$ submatrix is at most $\Delta$ in absolute value. 

\begin{statement}{Theorem \ref{IP}}
There is a function $f_{\ref{IP}}\colon \bZ\to \bZ$ so that for all positive integers $\Delta$ and $m$, if $A$ is a rank-$m$ $\Delta$-modular matrix, then the number of distinct columns of $A$ is at most $m^2+f_{\ref{IP}}(\Delta)\cdot m$.
\end{statement}

This theorem can be used to prove that any integer program with constraint matrix $A$ can be solved efficiently; see \cite{Paat} and \cite{Glanzer} for a detailed discussion of the implications of this result. 
Note that there exist rank-$m$ totally unimodular matrices with $m^2+m$ distinct columns, so this result is tight up to the function $f_{\ref{IP}}$.
For $\Delta\ge 2$, the previous best result was $\Delta^{2+2\log_2\log_2(\Delta)}\cdot m^2+1$, due to Glanzer et al. \cite{Glanzer}.
The proof of Theorem \ref{IP} relies on the following proposition.


\begin{proposition} \label{Delta}
Let $\Delta$ be a positive integer, and let $\cM_{\Delta}$ denote the class of matroids with a representation as a {$\Delta$-modular} matrix. Then
\begin{enumerate}[(i)]
\item $\cM_{\Delta}$ is minor-closed,

\item $U_{2,2\Delta+2}\notin \cM_{\Delta}$, and

\item $\DG(2\lfloor\log_2(\Delta)\rfloor+2,\{1,-1\})\notin\cM_{\Delta}$.
\end{enumerate}
\end{proposition}
\begin{proof}
Given a matrix $A\in \bZ^{R\times E}$ and a set $X\subseteq E$, we write $A[X]$ for the submatrix of $A$ consisting of the columns indexed by $X$.
We first prove a claim.

\begin{claim} \label{upper-triangle}
Let $A\in \bZ^{R\times E}$ be a $\Delta$-modular matrix with full row-rank, and let $X\subseteq E$ so that $X$ is independent in $M(A)$. 
Then there is a $\Delta$-modular matrix $A_1$ which is row-equivalent to $A$, so that $A_1[X]$ is an upper-triangular matrix after removing all-zero rows.
\end{claim}
\begin{proof}
Let $|X|$ be minimal so that the claim is false. Since it is vacuously true when $|X|=0$, we may assume that there is some $e\in X$.
By the minimality of $|X|$, there is a $\Delta$-modular matrix $A'$ which is row-equivalent to $A$, so that the submatrix $A'[X-\{e\}]$ is an upper-triangular matrix after removing all-zero rows.
Let $R_X\subseteq R$ denote the rows of $A'$ which are all-zero in $A'[X-\{e\}]$.

Let $a_1,\dots,a_k$ denote the nonzero entries of $A'[e]$ in the rows indexed by $R_X$;
by scaling rows by $-1$, we may assume that each $a_i$ is positive.
We say that the row of $A'$ which has entry $a_i$ in column $A'[e]$ is row $i$ of $A'$. 
The following three cases describe row operations to perform iteratively on $A'$, corresponding to a version of the Euclidean algorithm: 
\begin{enumerate}[(1)]
\item If $a_1=a_i$ for all $i\in \{2,\dots,k\}$, then for each $i\in \{2,\dots, k\}$, subtract row $1$ from row $i$, and then cease all row operations.

\item If $a_1<a_i$ for all $i\in \{2,\dots,k\}$, then for each $i\in\{2,\dots, k\}$, subtract row $1$ from row $i$ until $0 < a_i \le a_1$.

\item If $a_1>a_i$ for some $i\in \{2,\dots,k\}$, then subtract row $i$ from row $1$ until $0 < a_1 \le a_i$. 
\end{enumerate}

After each application of (2) or (3), the sum of the entries in the column indexed by $e$ and rows indexed by $R_X$ is strictly smaller, and each positive entry remains positive.
Thus, at some point we must perform the operations in case (1).
Therefore, we obtain a matrix $A''$ which is row-equivalent to $A'$, so that $A''[e]$ has precisely one nonzero entry some row indexed by an element of $R_X$.
This implies that after swapping rows, the submatrix $A''[X]$ is upper-triangular, after removing all-zero rows.

Since we did not scale any row, each $\rank(A)\times \rank(A)$ submatrix of $A''$ has the same determinant as the corresponding submatrix of $A'$.
Thus, the matrix $A''$ is $\Delta$-modular, so the claim holds with $A_1=A''$.
\end{proof}

We now prove (i).

\begin{claim}
$\cM_{\Delta}$ is minor-closed.
\end{claim}
\begin{proof}
Clearly $\cM_{\Delta}$ is closed under deletion.
Let $M\in \cM_{\Delta}$, and let $e\in E(M)$. We will show that $M/e\in \cM_{\Delta}$.
Let $A\in \bZ^{[r(M)]\times E(M)}$ be a {$\Delta$-modular} representation of $M$.
By \ref{upper-triangle}, there is a $\Delta$-modular matrix $A_1$ which is row-equivalent to $A$, so that the column indexed by $e$ has precisely one nonzero entry. 
Then the submatrix of $A_1$ obtained by deleting the column indexed by $e$ and the row with a nonzero entry in that column is a representation of $M/e$ as a $\Delta$-modular matrix.
\end{proof}

The following claim proves (ii).

\begin{claim}
$U_{2,2\Delta+2}\notin \cM_{\Delta}$.
\end{claim}
\begin{proof}
Let $\Delta\ge 1$ and $k\ge 3$ be integers, and let $A$ be $\Delta$-modular representation of $U_{2,k}$ with two rows.
Since $U_{2,k}$ is simple, no two columns of $A$ are parallel. 
By Chebyshev's Theorem, there is a prime $p$ so that $\Delta<p\le 2\Delta$.
By taking each entry of $A$ modulo $p$, we obtain a matrix $A'$ over $\GF(p)$.
The determinant of each square submatrix of $A'$ has the same value modulo $p$ as the corresponding submatrix of $A$.
In particular, each $2\times 2$ submatrix of $A'$ has determinant $0$ modulo $p$ if and only if the corresponding submatrix of $A$ has determinant $0$, since $A$ is $\Delta$-modular and $\Delta<p$.
Thus, $A'$ represents $U_{2,k}$ over $\GF(p)$, so $k\le p+1\le 2\Delta+1$.
\end{proof}

The final claim proves (iii).

\begin{claim}
$\DG(2\lfloor\log_2(\Delta)\rfloor+2,\{1,-1\})\notin\cM_{\Delta}$.
\end{claim}
\begin{proof}
The claim holds if $\Delta=1$, since $\DG(2,\{1,-1\})\cong U_{2,4}$ is not in $\cM_2$ by the previous claim, so we may assume that $\Delta\ge 2$.
Let $k\ge 2$ be an even integer so that $\DG(k,\{1,-1\})\in \cM_{\Delta}$; one can check that $\DG(2,\{1,-1\})\in \cM_{\Delta}$ for all $\Delta\ge 2$.
Let $A$ be a $\Delta$-modular matrix with $k$ rows whose vector matroid is isomorphic to $\DG(k,\{1,-1\})$.
Let $B$ be a frame for $\DG(k,\{1,-1\})$.
By \ref{upper-triangle}, there is a $\Delta$-modular matrix $A_1$ which is row-equivalent to $A$, so that $A_1[B]$ is an upper-triangular matrix.
Since $B$ is a frame for $\DG(k,\{1,-1\})$, each pair of elements of $B$ spans a $U_{2,4}$-restriction. 
Since $U_{2,4}$ is not a regular matroid, the matrix $A_1$ has a $k\times k$ submatrix with $2\times 2$ block submatrices on the diagonal, each with determinant at least two in absolute value, and with all zeroes below. 
Thus, the determinant of this submatrix is at least $2^{k/2}$ in absolute value. 
Since the determinant is at most $\Delta$, we find that $k\le 2\lfloor\log_2(\Delta)\rfloor$, as desired.
\end{proof}

The three claims complete the proof.
\end{proof}


We now prove Theorem \ref{IP}.

\begin{proof}[Proof of Theorem \ref{IP}]
Let $A$ be a rank-$m$ {$\Delta$-modular} matrix so that no two columns of $A$ are equal.
Let $X$ be a maximum-size subset of the columns of $A$ such that no two are scalar multiples of each other by an integer other than $-1$.
We first bound $|X|$, using the fact that the vector matroid of $A$ is $\bR$-representable.

By Theorem \ref{LG rep}, there is no group $\Gamma$ of size at least two so that the vector matroid of $A$ has a $\LG^+(3,\Gamma)$-minor.
By Theorem \ref{Dowling fields}, there is no group $\Gamma$ of size at least three so that the vector matroid of $A$ has a $\DG(3,\Gamma)$-minor.
Then by Proposition \ref{Delta},  and Theorem \ref{main} with $(\ell,n,k,t)=(2\Delta, 3, 2\lfloor\log_2(\Delta)\rfloor+2, 2)$, we see that
\begin{align*}
|X|&\le 2\Big({m\choose 2}+f_{\ref{main}}\big(2\Delta, 3, 2\lfloor\log_2(\Delta)\rfloor+2, 2\big)\cdot m\Big)\\
&=m^2+\Big(2\cdot f_{\ref{main}}\big(2\Delta, 3, 2\lfloor\log_2(\Delta)\rfloor+2, 2\big)-2\Big) \cdot m.
\end{align*}

We will write $d(\Delta)=f_{\ref{main}}(2\Delta, 3, 2\lfloor\log_2(\Delta)\rfloor+2, 2)$ for convenience.
Thus, we need only show that there are few columns in $X$ with many parallel columns.

\begin{claim} \label{ind}
There is no independent set of $\lfloor\log_2(\Delta)\rfloor+1$ columns in $X$ which all have a scalar multiple not in $X$.
\end{claim}
\begin{proof}
If so, let $A_1$ be an invertible $m\times m$ submatrix of $A$ containing all columns in $X$ which can all be scaled by some integer at least two to obtain another column of $A$; such a matrix exists because $X$ is independent.
Now, scale $|X|$ columns of $A_1$ each by a factor of at least two in absolute value to obtain a submatrix $A_1'$ of $A$.
Since we scale $|X|$ columns by magnitude at least two, we have $|\det(A_1')|\ge 2^{|X|}|\det(A_1)|$, and thus 
$$|X|\le \log_2\bigg(\frac{|\det(A_1')|}{|\det(A_1)|}\bigg)\le \log_2(\Delta),$$
since $|\det(A_1')|\le \Delta$.
\end{proof}

Since $U_{2,2\Delta+2}\notin\cM_{\Delta}$, \ref{ind} implies that at most $(2\Delta)^{\lfloor\log_2(\Delta)\rfloor}$ columns in $X$ have a scalar multiple not in $X$.
Each set of pairwise parallel columns of $A$ has size at most $2\Delta$; if we scale a column by magnitude greater than $\Delta$, then any invertible $m\times m$ submatrix containing that column has determinant greater than $\Delta$ in absolute value.
Thus, there are at most $2\Delta \cdot (2\Delta)^{\lfloor\log_2(\Delta)\rfloor} = (2\Delta)^{1+\lfloor\log_2(\Delta)\rfloor}$ columns of $A$ which are not in $X$.
Therefore, the number of columns of $A$ is at most 
\begin{align*}
|X|+(2\Delta)^{1+\lfloor\log_2(\Delta)\rfloor}&\le m^2+\Big(2\cdot d(\Delta)-2\Big) \cdot m+(2\Delta)^{1+\lfloor\log_2(\Delta)\rfloor}\\
&\le m^2+\Big(2\cdot d(\Delta)-2 + (2\Delta)^{1+\lfloor\log_2(\Delta)\rfloor}\Big) \cdot m.
\end{align*}
Thus, the theorem holds with $c = 2\cdot d(\Delta)-2+(2\Delta)^{1+\lfloor\log_2(\Delta)\rfloor}.$
\end{proof}

\end{subsection}
\end{section}

 
\begin{section}{An Exact Theorem} \label{exact chap}
In this section, we use Theorem \ref{main} to prove a result, Theorem \ref{exact main}, which will easily imply Theorem \ref{mainc}.


\begin{subsection}{A New Connectivity Reduction} \label{reduc3}
We first prove an analogue of Theorem \ref{new reduction} for matroids with no $\LG^+(n,\Gamma)$-minor with $|\Gamma|\ge 2$, where the spanning minor we find is a Dowling geometry.
The proof is essentially identical to the proof of Theorem \ref{new reduction}, except we apply Theorem \ref{main} instead of Theorem \ref{clique minor}.


\begin{theorem}\label{stronger reduction} There is a function $r_{\ref{stronger reduction}}\colon \bR^8\to \bZ$ so that for all integers $\ell,k,t\ge 2$ and $r,s\ge 1$ and any real polynomial $p(x)=ax^2+bx+c$ with $a>\frac{t-1}{2}$, if $M\in\cU(\ell)$ has no $\LG^+(k,\Gamma')$-minor with $|\Gamma'|\ge 2$ and satisfies $r(M)\ge r_{\ref{stronger reduction}}(a,b,c,\ell,t,k,r,s)$ and $\elem(M)>p(r(M))$, then $M$ has a minor $N$ with $\elem(N)>p(r(N))$ and $r(N)\ge r$ such that either
\begin{enumerate}[(1)]
\item $N$ has a $\DG(r(N),\Gamma)$-restriction with $|\Gamma|\ge t$, or
\item $N$ is vertically $s$-connected and has an $s$-element independent set $S$ so that each $e\in S$ satisfies $\elem(N)-\elem(N/e)>p(r(N))-p(r(N)-1)$.
\end{enumerate}
\end{theorem}
\begin{proof}
We first define the function $r_{\ref{stronger reduction}}$.
Let $\nu=\nu_{\ref{define nu}}(a,b,c,\ell,k,r,s)$, and define $\hat r_1$ to be an integer so that 
$$(2s+1)a(x+y)+s(\nu+b)+c-as^2\le 2axy$$
and $p(x-s)\le p(x-s+1)$ for all real $x,y\ge \hat r_1$. 
Let $f$ be a function which takes in an integer $m$ and outputs an integer $f(m)\ge \max(r,2m,2\hat r_1)$ such that $p(x)-p(x-1)\ge ax+\ell^{\max(m,\hat r_1)}$ for all real $x\ge f(m)$.
Define $r_{\lceil\nu/a\rceil}=1$, and for each $i\in \{0,1,2,\dots,\lceil\nu/a\rceil-1\}$ recursively define $r_i$ to be an integer so that 
$$p(x)>(t-1){x\choose 2}+f_{\ref{main}}\Big(\ell,k,n_{\ref{spanning Dowling}}\big(a,b,c,\ell,f(r_{i+1}),s\big),t\Big)\cdot x$$
 for all $x\ge r_i$.
Such an integer $r_i$ exists because $a>\frac{t-1}{2}$.
Finally, define $r_{\ref{stronger reduction}}(a,b,c,\ell,t,k,r,s)=r_0$.

Let $M\in \cU(\ell)$ with no $\LG^+(k,\Gamma)$-minor with $|\Gamma|\ge 2$, such that $r(M)\ge r_0$ and $\elem(M)>p(r(M))$. 
Let $\cM$ denote the class of minors of $M$.
We may assume that $h_{\cM}(n)\le p(n)+\nu n$ for all $n\ge 1$, or else (2) holds by Lemma \ref{define nu}.
The following claim essentially finds some $\nu'$ so that the coefficient of the linear term of $h_{\cM}(n)$ is in the interval $[\nu'+b-a,\nu'+b+a]$.

\begin{claim}\label{nu'}
There is some $0\le \nu'< \nu$ and $i\ge 0$ so that $h_{\cM}(n)>p(n)+\nu' n$ for some $n\ge r_i$, and $h_{\cM}(n)\le p(n)+(\nu'+a)n$ for all $n\ge r_{i+1}$. 
\end{claim}
\begin{proof}
We will break up the real interval $[0,\nu]$ into subintervals of size $a$.
Define $\nu_i=ai$ for $i\in \{0,1,2,\dots,\lceil\frac{\nu}{a}\rceil\}$.
Let $i\ge 0$ be minimal so that $h_{\cM}(n)\le p(n)+\nu_{i+1}n$ for all $n\ge r_{i+1}$. 
This choice of $i$ is well-defined, because $i=\lceil \nu/a\rceil-1$ is a valid choice since $\nu_{\lceil \nu/a\rceil}\ge \nu$ and $h_{\cM}(n)\le p(n)+\nu n$ for all $n\ge 1=r_{\lceil \nu/a\rceil}$. 

If $i>0$, then $h_{\cM}(n)>p(n)+\nu_i$ for some $n\ge r_i$ by the minimality of $i$.
If $i=0$, then $M$ certifies that $h_{\cM}(n)>p(n)$ for some $n\ge r_0$.
Thus, there is some $i\ge 0$ so that $h_{\cM}(n)>p(n)+\nu_i n$ for some $n\ge r_i$, and $h_{\cM}(n)\le p(n)+\nu_{i+1}n$ for all $n\ge r_{i+1}$. 
Since $\nu_i+a=\nu_{i+1}$, we may choose $\nu'=\nu_i$.
Note that $\nu_i=ai< \nu$ since $i\le \lceil \frac{\nu}{a}\rceil-1$.
\end{proof}

By  \ref{nu'}, $M$ has a minor $M_1$ such that $r(M_1)\ge r_i$ and $\elem(M_1)>p(r(M_1))+\nu' r(M_1)$. 
By Theorem \ref{main} and definition of $r_i$, the matroid $M_1$ has a $\DG(n_{\ref{spanning Dowling}}(a,b,c,\ell,f(r_{i+1}),s),\Gamma)$-minor with $|\Gamma|\ge t$.
Then by Lemma \ref{spanning Dowling} with $r=f(r_{i+1})$ and $q=p+\nu'$, the matroid $M_1$ has a minor $N$ such that $r(N)\ge f(r_{i+1})$ and $\elem(N)>p(r(N))+\nu'r(N)$, and $N$ either has a $\DG(r(N),\Gamma)$-restriction or an $s$-element independent set $S$ so that each $e\in S$ satisfies $\elem(N)-\elem(N/e)>p(r(N))-p(r(N)-1)+\nu'$. 
We may assume that $N$ is simple.
Since $f(r_{i+1})\ge r$ and $\nu'\ge 0$ we may assume that $N$ is not vertically $s$-connected, or else either (1) or (2) holds.

Let $(A,B)$ be a partition of $E(N)$ so that $r_N(A)\le r_N(B)< r(N)$ and $r_N(A)+r_N(B)-r(N)<s-1$. 
Let $r_N=r(N)$ and $r_A=r_N(A)$ and $r_B=r_N(B)$.
We first show that $r_A\ge \max(\hat r_1,r_{i+1})$.
If not, then $r_B\ge r_N-r_A\ge \max(r_{i+1},\hat r_1)$, using that $r_N\ge f(r_{i+1})\ge \max(2r_{i+1},2\hat r_1)$.
Also, 
\setcounter{equation}{0}
\begin{align}
|B|=|N|-|A|&>p(r_N)+\nu'r_N-\ell^{\max(\hat r_1,r_{i+1})}\\
&\ge p(r_N-1)+(\nu'+a)r_N\\
&\ge p(r_B)+(\nu'+a)r_B.
\end{align}
Line (1) holds because $r_A<\max(\hat r_1,r_{i+1})$ and $N\in \cU(\ell)$, and line (2) holds because $r_N\ge f(r_{i+1})$.
Line (3) holds because $r_B\ge \hat r_1$, so $p(r_B)\le p(r_N-1)$ since $r_B\le r_N-1$.
But then $r_B\ge r_{i+1}$ and $|B|>p(r_B)+(\nu'+a)r_B$, which contradicts  \ref{nu'} and the choice of $\nu'$.
Thus, $r_B\ge r_A\ge \max(\hat r_1,r_{i+1})$. 
Then 
\begin{align*}
p(r_A+r_B-s)+\nu'(r_A+r_B-s)&\le p(r_N)+\nu' r_N\\
&<|A|+|B|\\
&\le p(r_A)+p(r_B)+(\nu'+a)(r_A+r_B),
\end{align*}
where the first inequality holds because $r_A+r_B-s\le r_N$ and $p(x-s)\le p(x-s+1)$ for all $x\ge \hat r_1$, and the last inequality holds by  \ref{nu'} because $r_B\ge r_A\ge r_{i+1}$.
Expanding $p(x)=ax^2+bx+c$ and simplifying, we have 
$$(2s+1)a(r_A+r_B)+s(\nu'+b)+c-as^2>2ar_Ar_B,$$
which contradicts that $r_A\ge \hat r_1$, since $\nu'< \nu$.
\end{proof}

\end{subsection}


\begin{subsection}{Locally Frame Matroids} \label{local Dowling}
We now prove a lemma which will help prove the uniqueness of Dowling geometries as extremal matroids in Theorem \ref{mainc}.
We first need a straightforward lemma about frame matroids.

\begin{lemma}\label{easy1}
If $M$ is a simple frame matroid and $e\in E(M)$, then 
no element on a line through $e$ of length at least four is the tip of a spike in $M/e$.
\end{lemma}
\begin{proof}
If $f$ is on a line of $M$ through $e$ of length at least four in $M$, then $f$ is parallel to a frame element in $M/e$, by Lemma \ref{frame obvious} (i).
Thus, $f$ is not the tip of a spike in $M/e$ by Lemma \ref{frame obvious} (iii).
\end{proof}

The following proposition relies on Proposition \ref{g-star}.
We freely use the fact that if $M$ is framed by $B$ and $e\in E(M)$ is on at least two lines of $M$ of length at least four, then $e$ is parallel to an element of $B$.


\begin{proposition} \label{locally frame}
Let $t\ge 1$ be an integer, and let $M$ be a simple matroid of rank at least seven for which $|M|=t{r(M)\choose 2}+r(M)$.
If there is some $e\in E(M)$ so that $\si(M/e)\cong \DG(r(M)-1,\Gamma)$ for some group $\Gamma$ with $|\Gamma|=t$, then either 
\begin{enumerate}[(1)]
\item there is a set $X\subseteq E(M)$ for which $r_M(X)\le 15$ and $M|X$ is either not a frame matroid, or has a $U_{2,t+3}$-minor, or

\item $e$ is the tip of a spike of rank at least five, or

\item $M\cong \DG(r(M),\Gamma)$.
\end{enumerate}
\end{proposition}
\begin{proof}
Assume for a contradiction that (1), (2), and (3) do not hold for $M$.

\begin{claim} \label{no frame}
There is no $B\subseteq E(M)$ so that $M$ is framed by $B$.
\end{claim}
\begin{proof}
If the claim is false, then, since $M$ has no $U_{2,t+3}$-restriction and $\elem(M)=t{r(M)\choose 2}+r(M)$, each pair of elements of $B$ spans a $U_{2,t+2}$-restriction of $M$.
Since no rank-4 subset of $E(M)$ has a $U_{2,t+3}$-minor, Theorem \ref{dowling} implies that there is a group $\Gamma'$ so that $M\cong \DG(r(M),\Gamma')$.
Then $\si(M/e)\cong \DG(r(M)-1,\Gamma')$, which implies that $\Gamma'\cong \Gamma$, since Dowling geometries are isomorphic if and only if their groups are isomorphic.
Thus, $M\cong \DG(r(M),\Gamma)$ and (3) holds, a contradiction.
\end{proof}

We now reduce to the case that $t\ge 2$.
Note that $\elem(M)-\elem(M/e)=t(r(M)-1)+1$.

\begin{claim}
$t\ge 2$.
\end{claim}
\begin{proof}
Assume for a contradiction that $t=1$.
Since (1) does not hold, each restriction of $M$ of rank at most $15$ is binary, by a classical result of Tutte \cite{Tutte1}.
Since every binary spike has an $F_7$-minor and $F_7$ is not a frame matroid, binary spikes are not frame matroids. 
This implies that $e$ is not the tip of a spike, since (1) and (2) do not hold.
Since $\elem(M)-\elem(M/e)=r(M)$ and $M$ has no $U_{2,4}$-restriction, there are $r(M)-1$ lines of length three of $M$ through $e$.
We claim that each transversal $B$ of the long lines of $M$ through $e$ is a frame for $M/e$.
Since $e$ is not the tip of a spike in $M$ and $\si(M/e)\cong M(K_n)$, the set $B$ corresponds to a spanning tree of $M/e$.
If this tree has a path $P$ with three edges, then there is some $f$ so that $P\cup\{f\}$ is a size-4 circuit. 
But then $e$ is the tip of a rank-3 spike in $M/f$, and (1) holds since binary spikes are not frame matroids.
Thus, $B$ corresponds to a spanning star of $K_n$ and is a thus a frame for $M/e$.
In particular, this implies that each pair $F,F'$ of long lines of $M$ through $e$ satisfies $|\cl_M(F\cup F')|\ge \elem_{M/e}(\cl_{M/e}(F\cup F'))+3=6$, since $\si(M/e)\cong M(K_n)$.

We will choose a specific transversal $B'$ of the long lines of $M$ through $e$, and show that $M$ is framed by $B'\cup\{e\}$.
Let $F_1$ be a long line of $M$ through $e$, and let $b_1\in F_1-\{e\}$.
Let $F$ be a long line through $e$ other than $F_1$.
Since $|\cl_M(F_1\cup F)|=6$ and (1) does not hold, we have $M|\cl_M(F_1\cup F)\cong M(K_4)$.
Since $e$ and $b_1$ are in a triangle, the matroid $M|\cl_M(F_1\cup F)$ has a unique frame containing $e$ and $b_1$, and the third frame element $b_F$ is on $F-\{e\}$.
Then $\{b_1\}\cup\{b_F\colon F \ne F_1 \text{ is a long line through } e\}$ is a transversal $B'$ of the long lines of $M$ through $e$.

We claim that $B'\cup\{e\}$ is a frame for $M$.
Clearly $B'\cup\{e\}$ is a basis for $M$ since $B'$ is a basis for $M/e$.
Let $x\in E(M)$.
Since $M/e$ is framed by $B'$, there are elements $b,b'\in B'$ so that $x\subseteq\cl_M(\{e,b,b'\})$.
If $b_1\in \{b,b'\}$, then $B'\cup\{e\}$ contains a frame for $M|\cl_M(\{e,b,b'\})$ and $x$ is spanned by two elements of $\{e,b,b'\}$, so we may assume that $b_1\notin \{b,b'\}$.
Then $r_M(\{e,b_1,b,b'\})=4$, and $|\cl_M(\{e,b_1,b,b'\})|=\elem_{M/e}(\cl_{M/e}(\{b_1,b,b'\}))+4=10$, since $B'$ frames $M/e$.
Since (1) does not hold, this implies that $M|\cl_M(\{e,b_1,b,b'\})\cong M(K_5)$.
Since $e$ and $b_1$ are in a triangle, this matroid has a unique frame $B_1$ containing $e$ and $b_1$.
This set $B_1$ contains the unique frame for $M|\cl_M(\{e,b_1,b\})$ and $M|\cl_M(\{e,b_1,b'\})$, and so $b,b'\in B_1$.
Since $r_M(B_1)=4$ and $e,b_1,b,b'\in B_1$, we have $B_1=\{e,b_1,b,b'\}\subseteq B'\cup\{e\}$.
Thus, $x$ is spanned by two elements of $B'\cup\{e\}$.
Therefore, $M$ is framed by $B'\cup\{e\}$, which contradicts \ref{no frame}.
\end{proof}


Since $t\ge 2$, we can exploit that fact that any rank-3 frame matroid with $3t+3$ elements and no $U_{2,t+3}$-restriction has three elements each on two lines of length at least four, and thus has a unique frame.
Let $(M/e)|T_1$ be a simplification of $M$, and let $T\subseteq T_1$ be a transversal of the nontrivial parallel classes of $M/e$.
Let $B$ be a frame for $(M/e)|T_1$.
There are two distinct cases for the structure of the long lines of $M$ through $e$.

\begin{claim} \label{no spike}
If $e$ is not the tip of a spike, then $e$ is on $r(M)-1$ lines of length $t+2$ so that each contains an element of $B$.
\end{claim}
\begin{proof}
Recall that $\elem(M)-\elem(M/e)=t(r(M)-1)+1$.
If $e$ is not the tip of a spike, then $e$ is on $r(M)-1$ lines of length $t+2$ since $M$ has no $U_{2,t+3}$-restriction.
Let $F$ be a line of $M$ through $e$.
If $F\cap B=\varnothing$, then there is an element of $F$ which is the tip of a spike $S$ of rank at most four in $M/e$, since each element of $(M/e)|(T_1-B)$ is the tip of a spike of rank at most four in $M$ since $(M/e)|T_1\cong \DG(r(M/e),\Gamma)$.
But then $M|(F\cup S)$ is not a frame matroid by Lemma \ref{easy1} and (1) holds, a contradiction.
\end{proof}

The structure is a bit more complex when $e$ is the tip of a spike.

\begin{claim} \label{line F}
If $e$ is the tip of a spike, then $(M/e)|T$ has a star-partition $(\cL,\{x\})$ so that $|\cl_M(\{e,x\})|=t+2$, and $|L|=t$ for each $L\in\cL$.
Moreover, $x\in B$ and for each $b\in B-\{x\}$ there some $L\in \cL$ so that $L\subseteq \cl_M(\{e,x,b\})$.
\end{claim}
\begin{proof}
Since (1) and (2) do not hold, by Proposition \ref{g-star} with $g=5$ there is a star-partition $(X,\cL)$ of $(M/e)|T$ so that each line of $M$ through $e$ and an element of $T-X$ has length three.
Note that each $L\in \cL$ satisfies $|L|\le t+2$, or else (1) holds.
Then $|\cL|\le r(M)-1$, or else (2) holds by taking the union of lines through $e$ and each element of a transversal of $\cL$.
Then there is some $L\in \cL$ so that $|L|\ge 2$, or else $\elem(M)-\elem(M/e)\le (r(M)-1)+2t+2<t(r(M)-1)+1$, using that $r(M)\ge 7$.
This implies that $|L|\le t$ for all $L\in \cL$, or else (1) holds by Lemma \ref{m=1}. 
Letting $m=|\cl_M(\{e\}\cup X)|$, we have 
\begin{align*}
t(r(M)-1)+1&=\elem(M)-\elem(M/e)\\
&=\sum_{L\in \cL}|L|+m-1\\
&\le t(r(M)-2)+m-1,
\end{align*}
which implies that $m\ge t+2$ and each $L\in \cL$ satisfies $|L|=t$.

By the same reasoning as in \ref{no spike}, $\cl_M(\{e,x\})\cap B\ne \varnothing$. 
Since $B\subseteq T_1$ and $\cl_M(\{e,x\})\cap T_1=\{x\}$ we have $x\in B$.
Fix some $L\in \cL$.
By the definition of a star-partition, the set $\{x\}\cup L$ is a line of $M/e$ of length $t+1\ge 3$, and thus spans some element $b_L\in B-\{x\}$ in $M/e$ since $x\in B$.
Thus, $L\subseteq \cl_M(\{e,x,b_L\})$, so for each $L\in\cL$ there is some $b_L\in B-\{x\}$ so that $\{e,x,b_L\}$ spans $L$ in $M$.
If $L\ne L'$ then $b_L\ne b_{L'}$, or else $L\cup X$ is not a flat of $(M/e)|T$, which contradicts the definition of a star-partition.
Since $|\cL|=|B-\{x\}|=r(M)-2$, for each $b\in B-\{x\}$ there is a unique $L\in \cL$ so that $L\subseteq \cl_M(\{e,x,b\})$.
\end{proof}

Let $b_1\in B$ be on a line of length $t+2$ through $e$; such an element exists by \ref{no spike} and \ref{line F}.
We will show that for each $b\in B-\{b_1\}$, the matroid $N=M|\cl_M(\{e,b_1,b\})$ is a frame matroid with a unique frame.
Note that $\elem_{M/e}(\cl_{M/e}(\{b_1,b\}))=t+2$ since $\si(M/e)$ is isomorphic to a Dowling geometry.
If $e$ is not the tip of a spike of $M$, then since $B$ is a transversal of the long lines of $M$ through $e$ by \ref{no spike} we have $\elem(N)=(t+2)+2t+1=3t+3$, and thus $N$ has a unique frame since it has no $U_{2,t+3}$-restriction. 
If $e$ is the tip of a spike of $M$, then by \ref{line F} the point $e$ is on $t$ lines of length three and a line of length $t+2$ in $N$, so $\elem(N)=(t+2)+2t+1=3t+3$ and again $N$ has a unique frame.


We now define a frame for $M$.
Fix some element $b_2\in B-\{b_1\}$, and let $\{e',b_1',b_2'\}$ be the unique frame for $M|\cl_M(\{e,b_1,b_2\})$, where $\{e',b_1'\}$ spans $\{e,b_1\}$.
Define a function $f$ from $B-\{b_1\}$ to $E(M)$ which maps $b$ to the unique element $f(b)$ so that $\{e',b_1',f(b)\}$ is a frame for $M|\cl_M(\{e,b_1,b\})$.
We will show that $\hat B=\{e',b_1'\}\cup \{f(b)\colon b\in B-\{b_1\}\}$ is a frame for $M$.
Clearly $|\hat B|\le r(M)$ and $B\cup \{e\}\subseteq \cl_M(\hat B)$ by the definition of $f$, so $\hat B$ is a basis of $M$.
Let $x\in E(M)$, so $x\in \cl_M(\{e,b_i,b_j\})$ for some $b_i,b_j\in B$, since $B$ is a frame for $M/e$.
Then $x\in \cl_M(\{e',b_1',f(b_i),f(b_j)\})$, since $\{e,b_1\}$ spans $\{e',b_1'\}$.
The matroid $M|\cl_M(\{e',b_1',f(b_i),f(b_j)\})$ is a frame matroid of rank at most four, and the frame contains $\{e',b_1',f(b_i),f(b_j)\}$, as each of these elements is in the unique frame for $M|\cl_M(\{e,b_1,b_i\})$ or $M|\cl_M(\{e,b_1,b_j\})$.
Since $\{e',b_1',f(b_i),f(b_j)\}$ spans $\{e,b_1,b_i,b_j\}$ and is contained in a frame, the set $\{e',b_1',f(b_i),f(b_j)\}$ is a frame for $M|\cl_M(\{e,b_1,b_i,b_j\})$.
Therefore, $x$ is spanned by two elements of $\hat B$.
But then $M$ is framed by $\hat B$, which contradicts \ref{no frame}.
\end{proof}
\end{subsection}


\begin{subsection}{The Proof}
We now prove a general result which will easily imply Theorem \ref{mainc}.
The proof of the upper bound follows the same outline as the proof of Theorem \ref{main}.

\begin{theorem}\label{exact main}
There is a function $r_{\ref{exact main}}\colon \bZ^4\to \bZ$ so that for all integers $t\ge 1$, $\ell\ge 2$, and $k,n\ge 3$, if $\cM$ is a minor-closed class of matroids so that $U_{2,\ell+2}\notin\cM$, then either 
\begin{itemize}
\item $\cM$ contains $\LG^+(n,\Gamma')$ for some nontrivial group $\Gamma'$, or
\item $\cM$ contains a nontrivial extension of $\DG(k,\Gamma)$ with $|\Gamma|\ge t$, or
\item each $M\in \cM$ with $r(M)\ge r_{\ref{exact main}}(\ell,t,k,n)$ satisfies $\elem(M)\le t{r(M)\choose 2}+r(M)$. Moreover, if $r(M)\ge r_{\ref{exact main}}(\ell,t,k,n)$ and $\elem(M)=t{r(M)\choose 2}+r(M)$, then $\si(M)$ is isomorphic to a Dowling geometry.
\end{itemize}
\end{theorem}
\begin{proof}
We first define a sequence of large integers, ending with $r_{\ref{exact main}}(\ell,t,k,n)$.
Define $n_0=\max(r_{\ref{Dowling minor}}(\ell,k,n),r_{\ref{main porcupines}}(\ell,k,1,n,1,5))$, and
define $r_0$ to be an integer so that $t{r\choose 2}+r>(t-1){r\choose 2}+f_{\ref{main}}(\ell,t,n_0,n) r$ for all $r\ge r_0$.
This will allow us to find a $\DG(n_0,\Gamma)$-minor with $|\Gamma|\ge t$.
Define $r_1=\max(k,r_0)$, and $r_2= r_{\ref{stronger reduction}}(\frac{t}{2},1,0,\ell,t,k,r_1,200)$.

This integer $r_2$ is large enough to show that $h_{\cM}(r)\le t{r\choose 2}+r$ for all $r\ge r_2$, but we must define a larger integer to prove that each extremal matroid is a Dowling geometry.
Define $n_1\ge r_2$ to be an integer so that 
$$txy\ge 400(x+y)+200(1-t/2)-200^2$$
for all real $x,y\ge n_1$.
Define $n_2$ to be an integer so that $$\frac{t}{2}(2x-1)>(1-t/2)+\ell^{r_2}+\ell^{n_1}$$
for all real $x\ge n_2$.
Finally, define $r_{\ref{exact main}}(\ell,t,k,n)$ to be an integer so that $t{r\choose 2}+r>(t-1){r\choose 2}+f_{\ref{main}}(\ell,t,n_2,n)\cdot r$ for all $r\ge r_{\ref{exact main}}(\ell,t,k,n)$.
We will write $\DG(k,\Gamma)+e$ to denote any nontrivial extension of $\DG(k,\Gamma)$.
Assume that $\cM$ contains no $\LG^+(n,\Gamma)$ with $|\Gamma|\ge 2$, and no $\DG(k,\Gamma)+e$ with $|\Gamma|\ge t$.
The following claim provides sufficient conditions for finding a $(\DG(k,\Gamma)+e)$-minor with $|\Gamma|\ge t$.


\begin{claim} \label{find minor}
Let $N\in\cM$ be a vertically 200-connected matroid so that $N$ has a $\DG(n_0,\Gamma)$-minor with $|\Gamma|\ge t$.
If $N$ has either a restriction of rank at most $15$ which is not in $\cF\cap \cU(t+1)$ or a spike restriction of rank at least five, then $N$ has a $(\DG(k,\Gamma)+e)$-minor with $|\Gamma|\ge t$.
\end{claim}
\begin{proof}
If there is some $X\subseteq E(N)$ so that $r(N|X)\le 15$ and $N|X\notin \cF\cap \cU(t+1)$, then  by Theorem \ref{Dowling minor} with $m=k$ and $s=15$, the matroid $N$ has a minor $N_1$ with a spanning $\DG(k,\Gamma)$-restriction so that $N_1|X=N|X$. 
Since $N_1|X$ is not a restriction of a Dowling geometry with group size at most $t$, the matroid $N_1$ has a $(\DG(k,\Gamma)+e)$-restriction.

If $N$ has a spike restriction of rank at least five, then by Theorem \ref{main porcupines} with $m=k$, $s=k=1$ and $g=5$, the matroid $N$ has a minor $N_2$ of rank at least $k$ with a $\DG(r(N_2),\Gamma)$-restriction and a spike restriction of rank at least five, using that $s_{\ref{main porcupines}}(1,1,5)<200$.
Since $r(N_2)\ge k$ and spikes of rank at least five are not frame matroids, $N_2$ has a $(\DG(k,\Gamma)+e)$-minor.
\end{proof}

We first prove the upper bound on the extremal function of $\cM$.

\begin{claim}
$h_{\cM}(r)\le t{r\choose 2}+r$ for all $r\ge r_2$.
\end{claim}
\begin{proof}
If this is false, let $M\in \cM$ so that $r(M)\ge r_2$ and $\elem(M)>t{r(M)\choose 2}+r(M)$.
By Theorem \ref{stronger reduction} with $r=r_1$ and $s=200$, $M$ has a minor $N$ so that $r(N)\ge r_1$ and $\elem(N)>t{r(N)\choose 2}+r(N)$ and $N$ has either a $\DG(r(N),\Gamma)$-restriction with $|\Gamma|\ge t$, or is vertically $200$-connected and has an element $e$ so that $\elem(N)-\elem(N/e)>t(r(N)-1)+1$.
We may assume that the first outcome does not hold or else $N$ has a $(\DG(r(N),\Gamma)+e)$-restriction, and thus a $(\DG(k,\Gamma)+e)$-minor since $r(N)\ge k$.
Since $r(N)\ge r_0$ and $\elem(N)>t{r(N)\choose 2}+r(N)$, we have $\elem(N)>(t-1){r(N)\choose 2}+f_{\ref{main}}(\ell,t,n_0,n)\cdot r(N)$, so $N$ has a $\DG(n_0,\Gamma)$-minor with $|\Gamma|\ge t$ by Theorem \ref{main}.

Since $N$ has an element $e$ so that $\elem(N)-\elem(N/e)>t(r(N)-1)+1$, by Theorem \ref{critical structure} with $h=0$, $N$ either has a restriction of rank at most $15$
which is not in $\cF\cap \cU(t+1)$, or a spike restriction of rank at least five.
Here we use the fact that the union of two spikes with tip $e$ which are skew in $M/e$ is not a frame matroid, which follows from Lemma \ref{star}.
Then $N$ has a $(\DG(k,\Gamma)+e)$-minor by \ref{find minor}, a contradiction.
\end{proof}


We now prove the uniqueness of Dowling geometries as extremal matroids.
Let $M\in \cM$ be a simple matroid so that $r(M)\ge r_{\ref{exact main}}(\ell,t,k,n)$ and $|M|= t{r(M)\choose 2}+r(M)$.
Assume for a contradiction that $M$ is not isomorphic to a Dowling geometry.
Since $r(M)\ge r_{\ref{main}}(\ell,t,n_2,n)$, the matroid $M$ has a $\DG(n_2,\Gamma)$-minor $G$ with $|\Gamma|\ge t$, by Theorem \ref{main}.
Let $C\subseteq E(M)$ so that $G$ is a restriction of $M/C$.
Let $C_1$ be a maximal subset of $C$ so that $\elem(M/C_1)=t{r(M/C_1)\choose 2}+r(M/C_1)$ and $\si(M/C_1)$ is not isomorphic to a Dowling geometry.
Let $M_1$ be a simplification of $M/C_1$, and note that $C_1\ne C$.

\begin{claim} \label{200-conn}
$M_1$ is vertically $200$-connected.
\end{claim}
\begin{proof}
If not, then there is a partition $(A,B)$ of $E(M_1)$ with $r_{M_1}(A)\le r_{M_1}(B)<r(M_1)$ so that $r_{M_1}(A)+r_{M_1}(B)<r(M_1)+200$.
Let $r_A=r_{M_1}(A)$ and $r_B=r_{M_1}(B)$.
We first show that $r_A\ge n_1$.
If not, then 
\begin{align*}
|M_1|=|A|+|B|&<\ell^{n_1}+\max\Big(\ell^{r_2},t{r_B\choose 2}+r_B\Big)\\
&\le \ell^{n_1}+\ell^{r_2}+t{r_B\choose 2}+r_B\\
&<t{r(M_1)\choose 2}+r(M_1),
\end{align*}
a contradiction.
The last line holds since $r(M_1)\ge r(G)\ge n_2$ and $r_B<r(M_1)$, and by the definition of $n_2$.
Thus, $r_B\ge r_A\ge n_1\ge r_2$.
Then using that $r_A+r_B-200<r(M_1)$, we have
\begin{align*}
t{r_A+r_B-200\choose 2}+r_A+r_B-200&<t{r(M_1)\choose 2}+r(M_1)\\
&=|M_1|=|A|+|B|\\
&\le t{r_A\choose 2}+r_A+t{r_B\choose 2}+r_B.
\end{align*}
After expanding these polynomials and rearranging, this contradicts that $r_B\ge r_A\ge n_1$.
\end{proof}

Let $e\in C-C_1$.
By the maximality of $C_1$, either $\elem(M_1/e)<t{r(M_1/e)\choose 2}+r(M_1/e)$, or $\si(M_1/e)$ is isomorphic to a Dowling geometry.
If $\elem(M_1/e)<t{r(M_1/e)\choose 2}+r(M_1/e)$, then $\elem(M_1)-\elem(M_1/e)>1+t(r(M_1)-1)$.
Then by Theorem \ref{critical structure} with $h=0$, the matroid $M_1$ either has a restriction of rank at most $15$ which is not in $\cF\cap \cU(t+1)$, or a spike restriction of rank at least five.
By  \ref{find minor} and \ref{200-conn}, this implies that $M_1$ has a $(\DG(k,\Gamma)+e)$-minor.
This is a contradiction, so $M_1/e$ is isomorphic to a Dowling geometry.
Then by Proposition \ref{locally frame}, the matroid $M_1$ either has a restriction of rank at most $15$ which is not in $\cF\cap\cU(t+1)$, or a spike restriction of rank at least five.
Again, \ref{find minor} and \ref{200-conn} imply that $M_1$ has a $(\DG(k,\Gamma)+e)$-minor with $|\Gamma|\ge t$, a contradiction.
Therefore, $M$ is isomorphic to a Dowling geometry.
\end{proof}

\end{subsection}
\end{section}


\begin{section}{Proofs} \label{proofs}
In this section we prove Theorems \ref{mainc}, \ref{mainr}, and \ref{rep exact}-\ref{alg GF(p^k)}.
We rely on a class of rank-3 matroids called Reid geometries.
A \emph{Reid geometry} is a simple rank-3 matroid $R$ consisting of long lines $L_1,L_2,L_3$ with a common intersection point $x$ so that $L_3=\{x,y,z\}$ has length three. 
The \emph{incidence graph} $I(R)$ of $R$ is the bipartite graph with bipartition $(L_1-\{x\},L_2-\{x\})$ such that $a_1\in L_1-\{x\}$ is adjacent to $a_2\in L_2-\{x\}$ if and only if $\{a_1,a_2,y\}$ or $\{a_1,a_2,z\}$ is a line of $R$. 
Note that $I(R)$ has maximum degree at most two, since each $a_1\in L_1-\{x\}$ is on at most one long line of $R$ with each of $y$ and $z$.

It is not hard to show that for each integer $m\ge 2$, there is a unique Reid geometry whose incidence graph is a cycle of length $2m$; we denote this matroid by $R(m)$.
Reid geometries are useful for us because of the following result of Gordon \cite{Gordon}.

\begin{theorem}[Gordon] \label{alg Reid}
For each integer $p\ge 2$, the Reid geometry $R(p)$ is algebraic over a field $\bF$ if and only if $\bF$ has characteristic $p$.
\end{theorem}

We now show the connection between Dowling geometries and Reid geometries, which is that every nontrivial extension of a Dowling geometry with lines as long as possible has an $R(m)$-minor for some integer $m\ge 2$. 
Recall that an extension $M$ of a Dowling geometry is \emph{nontrivial} if $M$ is simple, and has no coloops.

\begin{proposition}\label{Dowling}
For all integers $\ell\ge 2$ and $k\ge 3$, if $M\in \cU(\ell)$ is isomorphic to a nontrivial extension of $\DG(k,\Gamma)$ with $|\Gamma|=\ell-1$, then $M$ has an $R(m)$-minor for some $2\le m\le \ell$. 
\end{proposition}
\begin{proof}
Let $e$ be an element of $M$ so that $M\del e\cong \DG(k,\Gamma)$ with $|\Gamma|=\ell-1$.
Let $B$ be a frame for the spanning $\DG(k,\Gamma)$-restriction of $M$, and let $B_e\subseteq B$ so that $B_e\cup\{e\}$ is the unique circuit contained in $B\cup\{e\}$ which contains $e$. 
Then $|B_e|\ge 3$, or else $\cl_M(B_e)$ has a $U_{2,\ell+2}$-restriction. 
Let $b_1,b_2,b_3$ be distinct elements of $B_e$, and let $N=M/(B-\{b_1,b_2,b_3\})$. 
Let $L_1=\cl_{N}(\{b_1,b_2\})$ and $L_2=\cl_{N}(\{b_1,b_3\})$. 
There is some element $y$ spanned by $\{b_2,b_3\}$ such that $\{b_1,y,e\}$ is a line, or else $N/b_1$ has a $U_{2,\ell+2}$-restriction. 
Let $L_3=\{b_1,y,e\}$, and $R=N|(L_1\cup L_2\cup L_3)$, so $R$ is a Reid geometry. 
Note that each long line of $R$ through $e$ or $y$ other than $L_3$ contains precisely one element from both $L_1$ and $L_2$.
Then each element of $(L_1\cup L_2)-\{b_1\}$ is on a long line with $y$,
or else $\elem(N/y)\ge \elem(R)-(\ell-1)=\ell+2$.
Similarly, each element of $(L_1\cup L_2)-\{b_1\}$ is on a long line with $e$.
Thus, each vertex of $I(R)$ has degree at least two in $I(R)$, so $I(R)$ contains a cycle. 
Thus, $R$ has an $R(m)$-restriction with $2\le m\le \ell$.
\end{proof}

The following lemma shows that if $|\Gamma|\ge 2$, then $\LG^+(k,\Gamma)$ has an $R(m)$-restriction.

\begin{proposition} \label{doubled clique}
Let $k\ge 3$ be an integer, and let $\Gamma$ be a nontrivial finite group. 
Then $\LG^+(k,\Gamma)$ has an $R(m)$-restriction for some $2\le m\le |\Gamma|$.
Moreover, if $\Gamma\cong \bZ_p$ for a prime $p$, then $\LG^+(k,\Gamma)$ has an $R(p)$-restriction.
\end{proposition}
\begin{proof}
 Let $L_1$ and $L_2$ be distinct long lines through $e_0$, and let $y,z$ be distinct elements other than $e_0$ on a third long line through $e_0$; then the restriction to $(L_1\cup L_2)\cup\{y,z\}$ is a Reid geometry $R$.
Since $\LG^+(k,\Gamma)$ has no $U_{2,|\Gamma|+2}$-minor by Theorem \ref{identify LG}, each element in $(L_1\cup L_2)-\{e_0\}$ is on a long line with $y$ and a long line with $z$.
Thus, $I(R)$ contains a cycle, so $\LG^+(k,\Gamma)$ has an $R(m)$-restriction.

Now assume that $\Gamma\cong \bZ_p$.
There is some integer $m$ so that $\LG^+(k,\bZ_p)$ has an $R(m)$-restriction.
Then $\LG^+(k,\bZ_p)$ is algebraic over a field $\bF$ only if $\bF$ has characteristic $m$, by Theorem \ref{alg Reid}.
By Theorem \ref{LG rep}, $\LG^+(k,\bZ_p)$ is representable over $\GF(p)$, and therefore $m=p$.
\end{proof}

We use this to prove the following proposition.

\begin{proposition} \label{LG alg}
Let $k\ge 3$ be an integer, let $\Gamma$ be a nontrivial finite group, and let $\bF$ be a field.
If $\LG^+(k,\Gamma)$ is algebraic over $\bF$, then $\bF$ has positive characteristic $p$, and the order of each element of $\Gamma$ is a power of $p$.
\end{proposition}
\begin{proof}
By Proposition \ref{doubled clique}, there is some integer $p$ so that $\LG^+(k,\Gamma)$ has an $R(p)$-restriction.
Since $\LG^+(k,\Gamma)$ is algebraic over $\bF$, and $R(p)$ is a restriction of $\LG^+(k,\Gamma)$, Theorem \ref{alg Reid} implies that $\bF$ has characteristic $p>0$.

Let $q$ be a prime so that some element of $\Gamma$ has order $q$.
By Proposition \ref{doubled clique}, $\LG^+(k,\Gamma)$ has an $R(q)$-restriction, and so $\LG^+(k,\Gamma)$ is only algebraic over fields of characteristic $q$.
Thus, $q=p$, since $\bF$ has characteristic $p$.
\end{proof}


We now prove our main results; as discussed in Section \ref{results}, we need only prove Theorems \ref{rep approx}, \ref{alg exact}, and \ref{alg GF(p^k)}.
We start with Theorem \ref{rep approx}.

\begin{proof}[Proof of Theorem \ref{rep approx}]
By Theorems \ref{doubled clique} and \ref{alg Reid}, $\cM$ contains no $\LG^+(3,\Gamma)$ with $|\Gamma|\ge 2$.
By Theorem \ref{Dowling fields}, $\cM$ contains no Dowling geometry with group size greater than $\alpha$.
Thus, $$h_{\cM}(r)\le \alpha{r\choose 2}+f_{\ref{main}}(\ell,\alpha+1,3^{2^{\ell}},3)\cdot r$$ for all $r\ge 0$, by Theorem \ref{main}.
By Theorem \ref{Dowling fields}, $\cM$ contains all Dowling geometries over the cyclic group of size $\alpha$, so $h_{\cM}(r)\ge \alpha{r\choose 2}+r$ for all $r\ge 0$.
\end{proof}

We now prove Theorem \ref{alg exact}.

\begin{proof}[Proof of Theorem \ref{alg exact}]
By Proposition \ref{doubled clique} and Theorem \ref{alg Reid}, $\cM$ contains no $\LG^+(3,\Gamma')$ with $|\Gamma'|\ge 2$.
By Proposition \ref{Dowling} and Theorem \ref{alg Reid}, $\cM$ contains no nontrivial extension of $\DG(3,\Gamma)$ with $|\Gamma|\ge t$.
Thus, Theorem \ref{exact main} applies.
\end{proof}

We prove one more exact result, which is stronger than Theorem \ref{alg ff1}.
It was conjectured for representable matroids in \cite{Structure Theory}.

\begin{theorem}
Let $\bF$ be a finite field of characteristic $p$.
If $r$ is sufficiently large, and $M$ is a simple rank-$r$ matroid algebraic over $\bF$ and with no $R(p)$-minor and no $U_{2,|\bF|+2}$-minor,
then $|M| \le (|\bF|-1){r\choose 2}+r$.
Moreover, if equality holds, then $M$ is isomorphic to a Dowling geometry over $\bF^{\times}$.
\end{theorem}
\begin{proof}
Assume for a contradiction that $\cM$ contains some $\LG^+(3,\Gamma')$ with $|\Gamma'|\ge 2$.
By Proposition \ref{doubled clique}, this matroid has an $R(m)$-minor $R$ with $2\le m\le |\Gamma'|$.
Since $R\in\cM$ and $\bF$ has characteristic $p$, Theorem \ref{alg Reid} implies that $m=p$.
Thus, $R(p)\in\cM$, a contradiction, so $\cM$ contains no $\LG^+(3,\Gamma')$ with $|\Gamma'|\ge 2$.
By Proposition \ref{Dowling} and similar reasoning, $\cM$ does not contain any nontrivial extension of $\DG(3,\Gamma)$ with $|\Gamma|\ge |\bF|-1$.
Thus, Theorem \ref{exact main} applies.
\end{proof}

Finally, we prove Theorem \ref{alg GF(p^k)}.

\begin{proof}[Proof of Theorem \ref{alg GF(p^k)}]
Assume for a contradiction that $\cM$ contains some $\LG^+(r(N),\Gamma')$ with $|\Gamma'|\ge 2$.
Then the order of each element of $\Gamma'$ is a power of $p$, by Proposition \ref{LG alg}.
In particular, $\Gamma'$ has a subgroup isomorphic to $\bZ_p$.
But then $\LG^+(r(N),\Gamma')$ has a restriction isomorphic to $N$, a contradiction.
Then by Theorem \ref{main}, there is a group $\Gamma$ and a constant $c$ so that $\cM$ contains all $\Gamma$-frame matroids, and $h_{\cM}(r)\le |\Gamma|{r\choose 2}+c \cdot r$, for all $r\ge 0$.
Then $|\Gamma|\le p^k-1$ since $U_{2,p^k+2}\notin\cM$, and so the first statement holds.
If $N$ is not a $(\bZ_{p^k-1})$-frame matroid, then $\cM$ contains all $(\bZ_{p^k-1})$-frame matroids by Theorem \ref{Dowling fields}, so the second statement holds.
\end{proof}


We conclude with three interesting problems regarding algebraic matroids.
Theorem \ref{alg Reid} and Proposition \ref{doubled clique} show that lift geometries are not algebraic over fields of characteristic zero, and a result of Lindstr\"om \cite{Lindstrom} shows that if a matroid is algebraic over a field $\bF$, then it is algebraic over the prime field of $\bF$.
Thus, the following conjecture would go some way towards characterizing the fields over which $\LG^+(k,\Gamma)$ is algebraic.

\begin{conjecture}
Let $\Gamma$ be a nontrivial finite group, and let $p$ be a prime.
Then $\LG^+(3,\Gamma)$ is algebraic over $\GF(p)$ if and only if there is some integer $j\ge 1$ so that $\Gamma\cong \bZ_p^j$.
\end{conjecture}

Proposition \ref{doubled clique} shows that if $\LG^+(3,\Gamma)$ is algebraic over a field of characteristic $p$, then the order of each element of $\Gamma$ is a power of $p$, but it is unclear whether $\Gamma$ must be Abelian.
We also make the analogous conjecture for Dowling geometries.

\begin{conjecture}
Let $\Gamma$ be a nontrivial finite group, and let $p$ be a prime.
Then $\DG(3,\Gamma)$ is algebraic over $\GF(p)$ if and only if there is some integer $j\ge 1$ so that $\Gamma$ is a subgroup of $\bZ_{p^j-1}$.
\end{conjecture}

The final problem is interesting in light of Theorem \ref{alg exact}.

\begin{problem}
Let $\ell\ge 2$ be an integer, and let $\cF$ be a family of fields.
When is there a matroid with no $U_{2,\ell+2}$-minor, which is algebraic over all fields in $\cF$, but not representable over all fields in $\cF$?
\end{problem}


\end{section}



\end{document}